\documentclass{amsart}

\usepackage{amsmath,amsfonts,amsthm,amssymb,amscd}
\usepackage{appendix}
\usepackage{amsaddr}
\usepackage{amsmath,bm}
\usepackage{cite}
\usepackage{geometry}
\usepackage{graphicx}
\usepackage{hyperref}
\usepackage{setspace}
\usepackage{url}
\newtheorem{remark}{Remark}
\newtheorem{theorem}{Theorem}[section]
\newtheorem{lemma}[theorem]{Lemma}
\newtheorem{proposition}{Proposition}[section]

\numberwithin{equation}{section}

\author{Chengyang Shao}
\thanks{
Department of Mathematics, Massachusetts Institute of Technology; \tt{shaoc@mit.edu}
}

\title{Long Time Behavior of a Quasilinear Hyperbolic System Modelling Elastic Membranes}

\begin{document}
\maketitle
\begin{spacing}{1.2}
\begin{abstract}
The paper studies the long time behavior of simplified model of elastic membrane driven by surface tension and inner air pressure. The system is a degenerate quasilinear hyperbolic one that involves the mean curvature, and also includes a damping term that models the dissipative nature of genuine physical systems. With the presence of damping, a small perturbation of the sphere converges exponentially in time to the sphere, and without the damping the evolution that is $\varepsilon$-close to the sphere has life span longer than $\varepsilon^{-1/6}$. Both results are proved using a new Nash-Moser-H?rmander type theorem proved by Baldi and Haus.
\end{abstract}

\section{Introduction}
\subsection{The Equation and the Main Results}
In this paper, we study a degenerate hyperbolic system that is related to small vibrations of an ideal elastic membrane governed by surface tension and air pressure. The mathematical model considered here is a toy model by its nature, but the methods used may shed light on genuine physical models.

Let $M$ be a smooth compact surface, modelling the topological configuration of the membrane. We fix an embedding $i_0:M\hookrightarrow\mathbb{R}^3$, and fix the induced metric $g_0=g(i_0)$ on $M$. We define $\mu_0$ to be the surface measure induced by $i_0$. On any time interval $[0,T)$, the motion of this membrane is described by a family of smooth embeddings $u:[0,T)\times M\to\mathbb{R}^3$. This description corresponds to Lagrangian coordinates in the realm of continuum mechanics. The geometric quantities of interest associated to $u$ consist of the induced metric $g(u)$, the induced surface measure $\mu(u)$, the outward pointing unit normal vector field $N(u)$, the volume $\text{Vol}(u)$ of the region enclosed by the surface $u(M)$, the second fundamental form $h(u)$ and the mean curvature vector field $-H(u)N(u)=\Delta_{g(u)}u$. All these geometric quantities will be regarded as mappings from $M$ to $\mathbb{R}^3$. Furthermore, we will also take into account the projection $\top_u$ to the tangent direction of $u(M)$ and the projection $\perp_u$ to the normal direction of $u(M)$.

We will let the motion of the membrane be governed by its own surface tension and a force due to pressure of the ideal gas enclosed. The model was first posed by Notz \cite{Notz2010}. Following this paper, we start with the following Lagrangian action:
\begin{equation}\label{LA}
\begin{aligned}
\mathcal{S}(u)&:=\int_0^T(\mathcal{K}(u)-\mathcal{A}(u)-\kappa\mathcal{I}(u))dt\\
&:=\int_0^T\left[\frac{1}{2}\int_M|\partial_tu|^2d\mu_0-\int_M d\mu(u)+\kappa\log\frac{\text{Vol}(u)}{\text{Vol}(u(0))}\right]dt,
\end{aligned}\tag{LA}
\end{equation}
where $\mathcal{K}(u)$, $\mathcal{A}(u)$ represent the kinetic energy and potential energy (proportional to surface area) of the membrane respectively. The material that fills the region bounded by $M$ will be assumed as ideal gas with its energy modeled by $\kappa {\mathcal I}(u)$, where $\kappa>0$ is a constant. The kinetic energy of the gas due to its macroscopic motion is neglected. The measure $\mu_0$ on $M$ is considered as a reference mass distribution, so the mass of a surface patch of area $d\mu(u)$ should be $d\mu_0/d\mu(u)$ when there is no mass transport in this idealistic system, which leads to the above expression of kinetic energy. The Euler-Lagrange equation takes the form
\begin{equation}\label{EQ0}
\frac{\partial^2u}{\partial t^2}
=\frac{d\mu(u)}{d\mu_0}\left(-H(u)+\frac{\kappa}{\mathrm{Vol}(u)}\right)N(u),\,
\left(
\begin{matrix}
u(0,x)\\
\partial_tu(0,x)
\end{matrix}
\right)=
\left(
\begin{matrix}
u_0(x)\\
u_1(x)
\end{matrix}
\right).
\tag{EQ0}
\end{equation}
This is the equation studied by Notz in \cite{Notz2010}, who proved its local well-posedness and obtained initial stability results. The equation falls into the genealogy of the equation
$$
\frac{\partial^2 u}{\partial t^2}=-H(u)N(u)
$$
suggested by S.T. Yau \cite{Yau2000}, who pointed out its relation to a vibrating membrane.

We would also like to introduce a damping term to the original action (\ref{LA}) to simulate dissipative features of a genuine physical system. The simplest way is to imitate the action of a damped harmonic oscillator:
$$
\int_0^Te^{bt}\left(\frac{m\dot x^2}{2}-\frac{kx^2}{2}\right)dt,\,b\geq0,
$$
and modify (\ref{LA}) as
\begin{equation}\label{LA'}
\begin{aligned}
\mathcal{S}(u)=\int_0^Te^{bt}\left[\frac{1}{2}\int_M|\partial_tu|^2d\mu_0-\int_M d\mu(u)+\kappa\log\frac{\text{Vol}(u)}{\text{Vol}(u(0))}\right]dt,\,b\geq0,
\end{aligned}\tag{Action'}
\end{equation}
whence the Euler-Lagrange equation becomes
\begin{equation}\label{EQWD}
\frac{\partial^2u}{\partial t^2}+b\frac{\partial u}{\partial t}
=\frac{d\mu(u)}{d\mu_0}\left(-H(u)+\frac{\kappa}{\mathrm{Vol}(u)}\right)N(u),\,
\left(
\begin{matrix}
u(0,x)\\
\partial_tu(0,x)
\end{matrix}
\right)=
\left(
\begin{matrix}
u_0(x)\\
u_1(x)
\end{matrix}
\right).
\tag{EQWD}
\end{equation}
We may view (\ref{EQ0}) as a special case of (\ref{EQWD}) with $b=0$, and deal with (\ref{EQWD}) unless the non-damped equation (\ref{EQ0}) needs to be investigated specifically. The term $b\partial_tu$ is usually referred as a ``weak damping" compared to strong damping of the form $-\Delta\partial_tu$.

We now set the stage to state the main results of this paper. Roughly speaking, they are concerned with the  evolution of spherical membranes under (\ref{EQWD}).

We know that the only compact embedded $C^2$ hypersurfaces of constant mean curvature in Euclidean spaces are spheres; this is known as the \emph{Alexandrov sphere theorem}. See \cite{MP2019} for a concise review of its proofs and refinements. Consequently, the only static embedded solutions of system (\ref{EQWD}) for compact $M$ must be Euclidean spheres. We thus assume, throughout the paper, that $M=S^2$, and without loss of generality that the static configuration is the unit sphere. So we fix in (\ref{EQWD})
$$\kappa=\frac{8\pi}{3}.$$
We may also fix $i_0:S^2\hookrightarrow\mathbb{R}^3$ as the standard embedding. It thus seems reasonable to conjecture that the evolution under (\ref{EQWD}) converges to a unit sphere if it starts from a small perturbation of the unit sphere. Nevertheless, the set of \emph{static solutions} possesses an obvious symmetry group:: if $\varphi:S^2\to S^2$ is a diffeomorphism and $a\in\mathbb{R}^3$ is any position vector, then $u(t)\equiv i_0\circ\varphi+a$ is still a static solution for (\ref{EQWD}). If $\varphi$ is close to the identity and $a$ is close to the origin, then $i_0\circ\varphi+a$ should be legitimately regarded as a perturbative solution near $i_0$, but it converges to (indeed, identically equals) \emph{another} embedding of $S^2$ different from $i_0$. The question thus arises: \emph{when $b>0$, which embedding do we expect the perturbative solution of (\ref{EQWD}) to converge to?}

The following theorem answers this question. It roughly states that a perturbation of the unit sphere converges to a ``nearby" unit sphere:

\begin{theorem}\label{Thm1}
Fix $M=S^2$ in equation (\ref{EQWD}). For a given $b>0$, set
$$
\beta=\left\{
\begin{matrix}
\displaystyle{{\frac{b-\sqrt{b^2-4}}{3}}}, & b\geq2\\
b/3, & b<2
\end{matrix}
\right.
$$
There is an $\varepsilon_0>0$ depending on $\beta$, such that if $\|u_0-i_0\|_{H^{41}}+\|u_1\|_{H^{41}}\leq\varepsilon<\varepsilon_0$, the solution $u$ to (\ref{EQWD}) exists globally in time, is of class $C^{3-j}([0,+\infty);C^{31+j}(S^2;\mathbb{R}^3))$, and is always a $C^2$ embedding. There is a $C^{31}$-diffeomorphism $\varphi:S^2\to S^2$ close to the identity mapping, a vector $a\in\mathbb{R}^3$ close to the origin, such that for some constant $C=C(b)$,
$$
\|u(t)-(i_0\circ\varphi+a)\|_{H^{33}_x}
\leq Ce^{-\beta t}\varepsilon,
\quad
|a|+\|\varphi-\mathrm{id}\|_{C^{31}}
\leq C\varepsilon.
$$
If in addition $u_0,u_1\in H^{n}$ for $n>41$, then
$$
\|u(t)-(i_0\circ\varphi+a)\|_{H^{n-8}_x}
\leq C_ne^{-\beta t}\|u_0-i_0\|_{H^{n}}+\|u_1\|_{H^{n}},$$
$$
|a|+\|\varphi-\mathrm{id}\|_{C^{n-10}}
\leq C_n\|u_0-i_0\|_{H^{n}}+\|u_1\|_{H^{n}}.
$$
\end{theorem}

The next natural question is: \emph{when the damping coefficient $b$ becomes zero, what is the lifespan of the vibration before singularities appear?} This question is partially answered by the following theorem:
\begin{theorem}\label{Thm2}
Fix $M=S^2$ in equation (\ref{EQ0}). There are positive numbers $\varepsilon_0>0$ and $C>0$ such that if $\|u_0-i_0\|_{H^{24}}+\|u_1\|_{H^{24}}\leq\varepsilon<\varepsilon_0$, the solution $u$ to (\ref{EQ0}) exists in the time interval $[0,T_\varepsilon)$ with $T_\varepsilon\sim{{\varepsilon}^{-1/6}}$, is of class $C^3([0,T_\varepsilon);C^{19}(S^2;\mathbb{R}^3))$, and is always a $C^2$ embedding throughout this time interval. Furthermore, there holds the following estimate:
$$
\|u(t)-i_0\|_{H_x^{21}}\leq (1+t)^3\varepsilon,\,
t\in[0,T_\varepsilon).
$$
If in addition $u_0,u_1\in H^{n}$ for $n>24$, then $u\in C^3([0,T_\varepsilon);H^{n-3}(S^2;\mathbb{R}^3))$.
\end{theorem}

\begin{remark}
Sobolev and H?lder norms in the above statements are all taken with respect to the metric $g_0$ induced by $i_0$. The methodology of \cite{Notz2010} applies to both (\ref{EQ0}) and (\ref{EQWD}), so they are locally-well-posed.
\end{remark}

\begin{remark}
Theorem \ref{Thm2} gives a longer lifespan estimate compared to the original $\log1/\varepsilon$ obtained in \cite{Notz2010}. The original one was obtained for the perturbation around a general static solution, i.e. a constant mean curvature hypersurface in a general ambient Riemannian manifold. However, for $S^2$ in $\mathbb{R}^3$ we obtain a much better lifespan estimate. This corresponds to the fact that $S^2$ is the \emph{stable} critical point of the area functional defined for surfaces enclosing a given volume.
\end{remark}

\subsection{Idea of the Proof}\label{Subsec1.2}
Before outlining the proof, we shall compare (\ref{EQ0}) and (\ref{EQWD}) with some well-studied geometric flows. Equation (\ref{EQ0}) seems very similar to the hyperbolic mean curvature flow (HMCF) introduced by LeFloch and Smoczyk in \cite{LeFlochSmoczyk2008}. Another way of introducing an HMCF is to modify the time derivative in the well-known mean curvature flow (MCF) to second-order derivative:
\begin{equation}\label{MCF}
\frac{\partial u}{\partial t}=-H(u)N(u)\quad
\Rightarrow\quad
\frac{\partial^2 u}{\partial t^2}=-H(u)N(u),
\end{equation}
as suggested by Yau in \cite{Yau2000}. However, both HMCFs in \cite{LeFlochSmoczyk2008} and \cite{Yau2000} are geometric evolution equations, just as the parabolic MCF, in the sense that they are invariant under diffeomorphisms of the underlying manifold, while the dynamical equation (\ref{EQWD}) is, as noticed by Notz in \cite{Notz2010}, not a geometric one. We shall explain this difference and reveal the resulting difficulty by linearizing (\ref{EQWD}). Write
$$
\Psi(u)=\frac{\partial^2u}{\partial t^2}+b\frac{\partial u}{\partial t}
-\frac{d\mu(u)}{d\mu_0}\left(-H(u)+\frac{\kappa}{\mathrm{Vol}(u)}\right)N(u),
$$
i.e., the nonlinear differential operator. The linearization of this operator around a given motion $u$ is already calculated in \cite{Notz2010}: for any $u,v\in C^\infty([0,T]\times S^2;\mathbb{R}^3)$,
\begin{equation}\label{LEQ1}
\begin{aligned}
\Psi'(u)v&=\frac{\partial^2v}{\partial t^2}+b\frac{\partial v}{\partial t}\\
&\quad-\frac{d\mu(u)}{d\mu_0}\left(\Delta_{g(u)}[v\cdot N(u)]+|h(u)|^2[v\cdot N(u)]
-\frac{\kappa}{\mathrm{Vol}(u)^2}\int_{S^2}[v\cdot N(u)] d\mu(u)\right)N(u)\\
&\quad-\frac{d\mu(u)}{d\mu_0}(\nabla^{g(u)} H(u)\cdot\top_uv)N(u)\\
&\quad-\frac{d\mu(u)}{d\mu_0}\left(-H(u)+\frac{\kappa}{\mathrm{Vol}(u)}\right)(\text{div}^{g(u)}{\top_u}v+H(u)[v\cdot N(u)])N(u)\\
&\quad+\frac{d\mu(u)}{d\mu_0}\left(-H(u)+\frac{\kappa}{\mathrm{Vol}(u)}\right)(\nabla^{g(u)}[v\cdot N(u)]-h^{kl}(u)[v\cdot\partial_lu]\partial_ku).
\end{aligned}
\end{equation}
This clearly shows that (\ref{EQWD}) is a \emph{highly degenerate} hyperbolic system: the principal symbol of the right-hand-side depends only on the \emph{normal} direction of $u(M)$ and vanishes for all tangent directions.

A similar degeneracy also occurs for the MCFs or the Ricci flow, which are evolutionary equations of obvious geometric significance. Hamilton \cite{Hamilton19822} and Gage and Hamilton \cite{GageHamilton1986} initially used the Nash-Moser technique to resolve this degeneracy problem to obtain local well-posedness results of Ricci flow or MCF Cauchy problems. In \cite{Hamilton19822}, Hamilton developed a scheme to deal with such degenerate parabolic systems whose principal symbol of linearization satisfies some certain integrability conditions. Fortunately, it was noticed by several authors that the geometric invariance of MCF allows one to reformulate the problem into a non-degenerate parabolic one and thus avoid using the complicated Nash-Moser scheme. One way is to transform through the DeTurck trick \cite{DeTurck1983}, by including the evolution of the gauge itself into the equation. A review of these techniques can be found in \cite{Mantegazza2011}.

However, as pointed out by Notz in \cite{Notz2010}, the right-hand-side of (\ref{EQWD}) is not invariant under diffeomorphisms due to the factor $d\mu(u)/d\mu_0$, whence he was forced to use the Nash-Moser inverse function theorem for the local well-posedness problem since the DeTurck trick does not apply. This argument can be further justified: the Lagrangian (\ref{LA}) is not invariant under a general diffeomorphism $\varphi\in\mathfrak{Diff}(S^2)$, but the only problematic term is the kinetic energy $\mathcal{K}(u)$, since the area and volume are invariant under a mere ``coordinate change" of the underlying surface. To keep $\mathcal{K}(u)$ unchanged, the diffeomorphism has to preserve the reference measure $\mu_0$ induced by $i_0$. Such diffeomorphisms form a closed subgroup of the diffeomorphism group with infinite codimension. Even if one applies the DeTurck trick, it is only possible to capture the symmetry governed by this subgroup, and there are still ``infinitely many degrees of freedom" that cannot be canceled. To summarize, the difficulty is degeneracy due to symmetry possessed by the space of static solutions, and impossibility to cancel this degeneracy due to lack of symmetry for the system itself.

We thus still employ the Nash-Moser technique. The linearized system (\ref{LEQ1}) is a so-called \emph{weakly hyperbolic linear system} (WHLS), which is a natural generalization of Hamilton's notion of weakly parabolic system in \cite{Hamilton19822}. The Cauchy problem of WHLS was studied in \cite{Notz2010}, where the author obtained a tame estimate for the inverse of the linearized differential operator $\Psi$ on the tame Fréchet  space $C^\infty([0,T]\times S^2;\mathbb{R}^3)$. We will sketch his results in Appendix \ref{B}. According to the general framework for the Nash-Moser category proposed by Hamilton \cite{Hamilton19821}, as long as such tame estimate holds, the original nonlinear problem is solvable. This is how local well-posedness was established. We shall continue to work along this outline.

In order to deal with the problem of symmetry, it is necessary to separate out all the possible symmetries. We shall follow the basic idea employed by Hintz and Vasy in \cite{HintzVasy2018} when investigating perturbations of Kerr-de Sitter spacetime: \emph{consider the unknown as a decaying perturbation of the eventual geometric configuration, which in turn is a part of the unknown}. In \cite{HintzVasy2018}, the evolution is governed by Einstein's equation of positive cosmological constant, and starts from a perturbation of an initial data set of a Kerr-de Sitter spacetime. The initial data set consists of a Cauchy surface, a Riemannian metric and a lapse tensor satisfying geometric constraints (Gauss-Codazzi equations). It is uniquely determined by four real blackhole parameters (a scalar mass and a vector angular momentum). Hintz and Vasy were able to show that if the angular momentum is small, then the evolution converges to a Kerr-de Sitter spacetime (to be precise, the difference between the solution metric with the Kerr-de Sitter metric decays exponentially fast in time) with \emph{possibly different} blackhole parameters. The Lorentz metric to be solved was decomposed to be a tuple of unknowns, including the eventual blackhole parameters and the exponentially decaying perturbation.

For (\ref{EQWD}), the eventual configuration in general takes the form $i_0\circ\varphi+a$, where $\varphi\in\mathfrak{Diff}(S^2)$ is a diffeomorphism and $a$ a vector in $\mathbb{R}^3$ signifying spatial shift. As commented above, there is no gauge invariance for (\ref{EQWD}), so it is not helpful to solve a geometric gauge. These two facts mark the major differences compared to \cite{HintzVasy2018}: the space of ``geometric parameters" in our problem is infinite dimensional, and since the system is not gauge-invariant, \emph{we do not introduce additional evolving gauge}.

The unknown will thus become (with a change in the meaning of symbol $u$)
$$
(X,a,u)\in\mathfrak{X}\times\mathbb{R}^3\times\mathbf{E},
$$
where $\mathfrak{X}$ is the Fréchet  space of all smooth tangent vector fields on $S^2$, and $\mathbf{E}$ is the Fréchet  space of time-dependent mappings from $S^2$ to $\mathbb{R}^3$ decaying exponentially in time. The space $\mathfrak{X}\times\mathbb{R}^3\times\mathbf{E}$ will be considered as the tangent space of the Fréchet  manifold
$$
\mathfrak{Diff}(S^2)\times\mathbb{R}^3\times\mathbf{E}.
$$
The triple $(X,a,u)$ shall correspond to a time-dependent embedding in a unique manner. The technical heart of this paper is the proof that the linearized equation (\ref{LEQ1}) has a solution in the Fréchet  space $\mathfrak{X}\times\mathbb{R}^3\times\mathbf{E}$, satisfying tame estimates required by the Nash-Moser scheme.

Let us briefly explain how a triple $(X,a,u)$ is obtained as a solution to the linearized problem. We shall write $\Xi(X,a,u)$ for the evolving embedding determined by $(X,a,u)$, whose precise form will be indicated later. For simplicity, let us linearize around the standard embedding $i_0$. This gives a linear system
\begin{equation}\label{LEQSimple}
\begin{aligned}
\frac{\partial^2\phi}{\partial t^2}+b\frac{\partial\phi}{\partial t}&=(\Delta_{g_0}+2)\phi-\frac{6}{4\pi}\int_{S^2}\phi d\mu_0,\\
\frac{\partial^2\psi}{\partial t^2}+b\frac{\partial\psi}{\partial t}&=0,
\end{aligned}
\end{equation}
where $\phi$ is a scalar function and $\psi$ is a tangent vector field along $i_0$. Although generally normal and tangent components are not decoupled in (\ref{LEQ1}), they actually behave as if they were decoupled under the energy norm introduced in \cite{Notz2010} for a WHLS, so this illustrative example can still be used.

The mapping $\phi N(i_0)+\psi$ is the linearized $\Xi(X,a,u)$. To re-obtain the vector field $X$, the shift vector $a$ and the decaying perturbation $u$, we notice two crucial geometric facts:
\begin{itemize}
\item All eigenvalues of the elliptic operator acting on $\phi$ in (\ref{LEQSimple}) are non-positive. This corresponds to the fact that the sphere is a \emph{stable} critical point of the area functional for surfaces enclosing a given volume (see Darbosa and do Carmo \cite{BD2012} for the definition), and ensures that in (\ref{LEQSimple}) most of the modes will be exponentially decaying.

\item The null space of this operator is exactly spanned by the three components of $N(i_0)$, or equivalently spherical harmonics with lowest eigenvalue. This corresponds to the differential geometric identity $\Delta N+|h|^2N+\nabla H=0$ for any embedded orientable surface, where $\Delta$ is the Laplacian on the surface, $N$ is the outward normal vector field, $h$ is the second fundamental form and $H$ is the mean curvature. Note that we do not distinguish between the outward normal vector field and the Gauss map.
\end{itemize}
Equation (\ref{EQ0}) was derived for a generic evolving submanifold in \cite{Notz2010}. However, for a generic constant mean curvature hypersurface in a generic ambient Riemannian manifold, no stability result is valid; for example, the section of $S^{n}\subset\mathbb{R}^{n+1}$ with an $n$-dimensional hyperplane gives a sphere of dimension $n-1$, which has constant mean curvature in $S^n$ but is not stable. Furthermore, for a generic constant mean curvature hypersurface in a generic ambient Riemannian manifold, the null space of the second variation of area functional may be larger than the span of the components of the normal. Hence the above two facts are specific for $S^n\subset\mathbb{R}^{n+1}$. Using the language of \cite{HintzVasy2018}, we point out that most of the modes of the linearized problem are decaying, and the zero modes are well-understood and do not destroy this decay.

Thus the solution to (\ref{LEQSimple}) takes the following form:
$$
\left(\begin{matrix}
\phi(t)\\
\psi(t)
\end{matrix}
\right)
=\sum_{k=1}^3\left(\begin{matrix}
\langle\phi(0)+b^{-1}\phi'(0),N^k(i_0)\rangle_{L^2(g_0)}N^k(i_0)\\
\psi(0)+b^{-1}\psi'(0)
\end{matrix}
\right)+\text{exponentially decaying terms}.
$$
This gives a $\mathbb{R}^3$-valued mapping
$$
\begin{aligned}
\sum_{k=1}^3\langle\phi(0)+b^{-1}\phi'(0),N^k(i_0)\rangle_{L^2(g_0)}&N^k(i_0)N(i_0)
+\psi(0)+b^{-1}\psi'(0)
+\text{exponentially decaying terms}.
\end{aligned}
$$
We immediately notice that $\sum_{k=1}^3\langle\phi(0)+b^{-1}\phi'(0),N^k(i_0)\rangle_{L^2(g_0)}N^k(i_0)N(i_0)$ is nothing but the projection of a constant vector $a$ in $\mathbb{R}^3$ along $N(i_0)$, with $a^k=\langle\phi(0)+b^{-1}\phi'(0),N^k(i_0)\rangle_{L^2(g_0)}$. We may then subtract from $\psi(0)+b^{-1}\psi'(0)$ the tangent projection of $a$ to obtain a tangent vector field $X$. This realignment then gives the desired $(X,a,u)$. This is the key ingredient for proving Theorem \ref{Thm1}.

It is of technical interest which version of Nash-Moser scheme should be chosen. We can certainly choose the simplest version, for example Saint-Raymond's account \cite{SR1989} (the one used by Hintz and Vasy \cite{HintzVasy2018}), or the ``structuralist" version, for example Hamilton's account \cite{Hamilton19821}, if we only care about well-posedness. However, we choose the Nash-Moser-H?rmander theorem presented by Baldi and Haus in \cite{BaldiHaus2017}. The advantage is that it gives an explicit estimate of the solution in terms of initial data, thus enabling us to estimate the lifespan in Theorem \ref{Thm2}. \cite{BaldiHaus2017} also provides, as noted by the authors, a sharp regularity result: under this version, ``the nonlinear problem reaches exactly the same regularity given by the linearized one". So the optimal regularity bound can be obtained once the optimal result is obtained for the linearized system.

With the aid of this Nash-Moser-H?rmander theorem, the energy estimates obtained in establishing tame estimates for Theorem \ref{Thm1} will automatically give the lifespan estimate in Theorem \ref{Thm2}. If the initial data is $\varepsilon$-close to the static solution, the lifespan will be approximately some negative power of $\varepsilon$. This is an almost global result, and it of course corresponds to the fact that the elliptic operator acting on $\phi$ in (\ref{LEQSimple}) is non-positive, which in turn corresponds to the stability of the sphere. In the original argument \cite{Notz2010}, the lifespan estimate around a given constant mean curvature hypersurface was $\log1/\varepsilon$, which was a weaker generic result compared to ours since no stability assumption was posed. Furthermore, we shall explain why this generic method cannot provide a lifespan estimate better than $\log1/\varepsilon$ in Appendix \ref{B}.

To get a heuristic about how the power is obtained, it is helpful to consider a prototype ODE problem
$$
u''(t)=A(t)u(t),
\quad
\text{with }A(t)\in M_n(\mathbb{R}^n),\,\|A(t)\|\leq\lambda\sim0.
$$
If $A(t)\equiv\lambda$, then the optimal growth estimate of $|u(t)|$ uniform in $\lambda$ is $|u(t)|\simeq te^{\sqrt{\lambda}t}$, which is obtained by explicitly solving this problem. For generic $A(t)$, it is necessary to consider a weighted energy norm
$$
E(t):=\left(|u'(t)|^2+\lambda|u(t)|^2\right)^{1/2}.
$$
Differentiating, using Young's inequality $ab\leq (c^ra^2+c^{2-r}b^2)/2$,
$$
\frac{1}{2}\frac{d}{dt}E(t)^2
=A(t)u\cdot u'+\lambda u'\cdot u
\leq\sqrt{\lambda}|u'|^2+\lambda^{3/2}|u|^2
\leq\sqrt{\lambda} E(t)^2.
$$
This gives the estimate $|u'(t)|\simeq e^{c\sqrt{\lambda}t}$, hence the estimate $|u(t)|\simeq te^{c\sqrt{\lambda}t}$ (uniform in $\lambda$), which is optimal in terms of the rate of exponential growth as illustrated by the example $A(t)\equiv\lambda$.

Returning to the original problem, we need to choose the length of the time interval to keep control of the growth of the norms. We thus need to obtain energy estimate for the perturbation of a system which has three non-trivial zero modes (as discussed above). The idea is to separate the non-growing modes and these three almost zero modes, and investigate them separately. The technical difficulty of this scheme was pointed out by H?rmander in section 6.5 of \cite{Hormander1997}: it is similar to the difficulty encountered when dealing with
$$
\square u=F(u,\partial u,\partial^2u),\,u:\mathbb{R}_t\times\mathbb{R}_x^n\to\mathbb{R},
$$
whence the energy estimate has to gain a factor of positive power in $t$ due to the presence of $u$ on the right-hand-side. For our case, this positive power is $t^3$. This factor undermines the expectation that the lifespan should be approximately $\varepsilon^{-1/2}$, but still gives an almost global result $\varepsilon^{-1/6}$.

To summarize the proof, we point out that several geometric facts enter into the analysis and play crucial roles in establishing stability results.

\subsection{Outline of the Paper}
The paper contains four sections devoted to the proof of Theorem \ref{Thm1} and \ref{Thm2} together with three appendices. The proof follows a Nash-Moser scheme, and the technical heart is to obtain tame estimates for the linearized operator (\ref{LEQ1}).
\begin{itemize}
\item In section 2, we introduce the Fréchet spaces that will be used, together with several functional theoretic lemmas concerning the tame estimates for several geometric mappings. There will essentially be two Fréchet spaces of time-dependent smooth mappings from $S^2$ to $\mathbb{R}^3$, one defined on a finite time interval $[0,T]$ and one consisting of exponentially decaying mappings.
\item In section 3, we first study the elliptic operator that acts on the normal component in (\ref{LEQ1}). We will obtain some spectral properties and some tame elliptic estimates for this operator, and use these properties to obtain a decay estimate for (\ref{LEQ1}) with $u=0$.
\item In section 4, we will study the full linearized operator (\ref{LEQ1}), based on the spectral properties and elliptic estimates obtained. We first derive an energy estimate which holds uniformly for $b\geq0$, and explicitly determine the dependence of the constants on the time interval $[0,T]$. Next, we will fix a $b>0$ and obtain the decay estimate for the full linearized problem in the space of exponentially decaying mappings. The stability of $S^2$ in $\mathbb{R}^3$ enters into these tame estimates.
\item Section 5 will be devoted to the proof of our main theorems. Using the tame estimates in section 4, the Nash-Moser scheme ensures the solvability of equation (\ref{EQWD}), and the refined Nash-Moser-H?rmander theorem in \cite{BaldiHaus2017} will help us determine the lifespan for (\ref{EQ0}) with initial data $\varepsilon$-close to the unit sphere.
\item Appendix \ref{A} briefly describes general weakly hyperbolic linear systems that are studied in \cite{Notz2010}. Appendix \ref{B} closely studies the methodology of \cite{Notz2010} in obtaining the lifespan for (\ref{EQ0}) and explains why the original method cannot give the optimal result.
\end{itemize}

\subsection{Physical Appropriateness of the Model}\label{Subsec1.4}
The action (\ref{LA}) looks similar to the well-known
\begin{equation}\label{LA''}
\int_0^T\left[\frac{1}{2}\int_{\mathbb{R}_x^2}|\partial_t u|^2dx -\int_{\mathbb{R}^2_x}\left(\sqrt{1+|\nabla u|^2}dx-1\right)dx\right]dt.
\end{equation}
Here $u$ is a scalar function, and the addendum $-1$ in the area term is to ensure integrability. It's Euler-Lagrange equation is
\begin{equation}\label{EQClassical}
\frac{\partial^2 u}{\partial t^2}=\frac{1}{\sqrt{1+|\nabla u|^2}}\sum_{i,j=1}^2\left(\delta_{ij}-\frac{\partial_iu\partial_ju}{1+|\nabla u|^2}\right)\partial_i\partial_ju.
\end{equation}
This equation is known as the equation of nonlinear vibrating elastic membrane, whose potential energy is proportional to its area. If the vibration is of small slope, the area term in (\ref{LA''}) is then approximately $\int_{\mathbb{R}^2_x}|\nabla u|^2$, and the Euler-Lagrange equation becomes the linear wave equation $\partial_t^2u-\Delta_xu=0$. This is a usual way how linear wave equation is derived (see, for example, Courant and Hilbert \cite{CouHi2008}, Volume I, Section IV.10). It is thus natural to generalize this argument to curved surfaces, and one possibility for such generalization is (\ref{EQ0}), with a curved manifold $M$ in place of $\mathbb{R}^2$ and $\mu_0$ in place of the Lebesgue measure, and the additional term relating to pressure. Note that the action (\ref{LA''}) is defined for vibrations \emph{assumed to be} normal to the static configuration $u\equiv0$ and has no tangential degree of freedom at all. This marks the aforementioned difference from the aspect of degeneracy between our problem and the classical problem (\ref{EQClassical}) (usually referred as ``quasilinear perturbation of linear wave equation").

It should still be emphasized that the model (\ref{EQ0}) (and its damped version (\ref{EQWD})) is merely an \emph{idealistic model} for micro-vibrations of membranes, just as the classical equation (\ref{EQClassical}), the linear wave equation $\partial_t^2u-\Delta_xu=0$, or any other PDE modelling genuine physical phenomenon. In fact, the model is suitable for small oscillations that are ``almost" normal to a standard sphere (although tangential deformation is allowed as explained), and with all but one degree of freedom of the surrounded (ideal) gas neglected. If no such smallness assumption is posed, the model may fail to be appropriate: the potential energy in (\ref{LA}) is invariant under \emph{any} diffeomorphism,
but if $u$ is composed with a diffeomorphism not close to the identity, it is unphysical to imagine that the potential energy remains unchanged. Thus we should restrict to vibrations of not very large magnitude. If we are also interested in taking more degrees of freedom of the enclosed gas into account (hence closer to a genuine scenario), we might need to replace the ideal gas law with other equation of status, or even consider the gas as Eulerian fluid, and the membrane as a free boundary.

Furthermore, the artificial damping term $b\partial_tu$ in (\ref{EQWD}) may not represent any physical dissipation, although it gives the model a form quite similar to the well-known damped harmonic oscillator $m\ddot{x}+b\dot x+kx=0$, or damped string vibration $\partial_t^2u+b\partial_tu-\partial_x^2u=0$. In fact, to take actual physical dissipation into account, it should be clarified what underlying mechanism results in the dissipation -- viscosity, dispersion, drag force due to friction, etc.. This obviously is not reflected in (\ref{EQWD}). Since appropriate description of any of the above mechanisms is too complicated and will significantly change the original equation (\ref{EQ0}), we limit to consider this idealized weak damping for simplicity.

On the other hand, the damping term $b\partial_tu$ is not a mere compromise on complexity to produce a mathematically acceptable object. The reason that we are interested in (\ref{EQWD}) is that, it is a simple prototype of nonlinear evolutionary problems with the following features:
\begin{itemize}
    \item The problem itself possesses degeneracy, either may or may not be removed by fixing a gauge;
    \item The linearized problem exhibits a uniform decay for most modes, with non-decaying modes controlled;
    \item The set of terminal configurations is non-trivial, so that only \emph{orbital stability} (instead of asymptotic stability) can be expected for the original problem.
\end{itemize}
As explained in subsection \ref{Subsec1.2}, the problem of stability of perturbed Kerr-de Sitter spacetime falls into this category (with the second feature requiring most technicality). The work of Hintz and Vasy \cite{HintzVasy2018} enlightened ours: they used a Nash-Moser scheme that takes the blackhole parameters as part of the unknown, and the technical heart is to understand how the non-decaying modes behave. We thus consider the damped equation (\ref{EQWD}) as a \emph{simple toy model} that best illustrates how a Nash-Moser scheme solving eventual configuration can resolve a nonlinear evolutionary problem with the above features. The artificial damping in (\ref{EQWD}) might be non-physical, but the means we treat it still has a chance to generalize to more physical scenarios.

To close this subsection, we turn to discuss whether (\ref{EQ0}) and (\ref{EQWD}) are suitable for describing soap bubbles, and point out that the method we develop for the toy model (\ref{EQWD}) might be applicable for a better model. Since soap films obey the Young-Laplace law ``$\text{Pressure difference}\propto\text{Mean curvature}$", it is not surprising that the (parabolic) mean curvature flow is usually referred as describing a moving soap film. In \cite{Notz2010}, the author considered the hyperbolic mean curvature flow (\ref{EQ0}) as ``idealized model for a moving soap bubble". In \cite{IYAH2017},  Ishida, Yamamoto, Ando and Hachisuk derived a hyperbolic mean curvature flow from two-dimensional Euler equation as a model for the motion of fluid films, and used it for numerical simulation of soap films. Their equation is similar to the prototype suggested by Yau \cite{Yau2000}.

However, fluid physicits have also worked hard to describe fluid films, and derived several models directly from the general principles of continuum mechanics. For example, Ida and Miksis \cite{IdaMiksis19981} \cite{IdaMiksis19982} proposed a model, and did some initial linear stability analysis. Their model involves the following equations (precise forms omitted due to limited space) :
\begin{itemize}
    \item General Navier-Stokes equation involving van der Waals force, since the liquid film is very thin;
    \item Diffusion-reaction equations of surfactant on the two liquid-gas interfaces;
    \item Boundary conditions on the two liquid-gas interfaces involving the surface tension (hence mean curvature).
\end{itemize}
A similar model was also obtained by Chomaz in \cite{Chomaz2001}. These models all seem to be very different from any mean curvature flow for surfaces, whether parabolic or hyperbolic. This raises the question whether mean curvature flow type equations, including (\ref{EQ0}) and (\ref{EQWD}), are suitable for soap bubbles.

Nevertheless, we do observe, in any real bubble-blowing game, that spherical bubbles are stable configurations for evolving soap films: if the original bubble is not of very large scale, then it stablizes to a round sphere within a few seconds, and remains so until vaporizes; if the bubble is of very large scale, say of 2 meters' diameter, then it keeps oscillating before it vaporizes. These facts do resemble Theorem \ref{Thm1}. Of course, dissipation plays an important role here (reflected in the viscosity term of the Navier-Stokes equation and diffusion term of the surfactant equation), since otherwise the configuration may not converge at all. But the mathematics behind this is extremely complicated. Fortunately, it is still not hard to conclude that the (very thin) co-centric spherical annulus is a static configuration of the Ida-Miksis system. It is thus very interesting to ask whether the annulus is \emph{orbitally stable} \footnote{Should it exist, such a stable thin annulus itself is always uniquely determined. However, if worked under Lagrange coordinate, the problem involves, at least implicitly, determining a Lebesgue measure-preserving diffeomorphism of the annulus to itself. That is, it should be determined how fluid particles redistribute eventually.} for the Ida-Miksis system. If the answer is positive, it certainly justifies the Navier-Stokes equation in describing thin fluid films. There is no wonder that the Ida-Miksis system is a highly complicated free boundary problem, but it exhibits, at least partially, the features mentioned above. From this aspect, (\ref{EQWD}) is not completely a mathematical amusement. The author will keep working on this aspect to verify whether the techniques developed in this paper -- Nash-Moser scheme solving eventual configuration -- applies to this more genuine dissipative system.

\subsection{Comparison with other Evolutionary Problems}
As already discussed above, equation (\ref{EQ0}) and (\ref{EQWD}) are different from the well-known mean curvature flow due to their non-geometric nature. In this subsection, we will compare these equations with other evolutionary problems in literature.

Damped hyperbolic systems have been objects of interest mathematics and physics for a long-time, although it sometimes inevitably introduces oversimplification of actual existing dissipation, as pointed out in subesection \ref{Subsec1.4}. Some classical results were collected in O. Ladyzhenskaya's monograph \cite{Ladyzhenskaya1991}, where the author treated damped equations of the following form in a Hilbert space $\mathcal{H}$:
$$
\frac{\partial^2 u}{\partial t^2}+\nu\frac{\partial u}{\partial t}+Au+f(u)=h,
$$
where $\nu\geq0$ is a constant, $A$ is assumed to be a \emph{positive-definite} operator on $\mathcal{H}$ with compact resolvent, and $f$ is a nonlinear and unbounded operator. The prototype of $A$ is, of course, the negative Laplacian $-\Delta$. Ladyzhenskaya developed a complete global well-posedness theory for $\nu>0$ by introducing the scales $\mathcal{H}_s(A):=\text{Range}(A^{s/2})$ (i.e. by imitating the classical Sobolev space) and assuming that the regularity loss under the nonlinear perturbation $f$ does not exceed that resulting from $A$; transferring to the usual Sobolev space setting, the perturbation $f(u)$ should consist of only first order derivatives of $u$. The monograph also proved exponential decay of the solution for $\nu>0$, and obtained results on global attractors of this system under further assumptions on $f$. However, equation (\ref{EQWD}) \emph{cannot} be reduced to the form indicated above, as can be illustrated by  linearizing around the static solution $i_0$; in fact, if we replace $u$ by $i_0+u$ and regard $u$ as the perturbation, the equation may be re-written as
$$
\begin{aligned}
\frac{\partial^2u}{\partial t^2}+b\frac{\partial u}{\partial t}
&=\Delta_{g_0}[u\cdot N(i_0)]N(i_0)+2\perp_{i_0}u\\
&\quad+\left(\frac{3}{2\pi}\int_{S^2}[u\cdot N(i_0)] d\mu_0\right)N(i_0)+(\text{Quadratic order terms in }u),
\end{aligned}
$$
where the principal part of the remainder is
$$
\frac{d\mu(i_0+u)}{d\mu_0}\Delta_{g(i_0+u)}(i_0+u)-\Delta_{g_0}[ u\cdot N(i_0)] N(i_0).
$$
It consists of \emph{both} normal and tangential (with respect to the fixed embedding $i_0$) second order spatial derivatives of $u$, while the would-be operator $A$ consists of only normal derivative of $u$; this reflects the \emph{degeneracy} of the system, which appears, as commented above, impossible to be gauged and counteracted. Theorem \ref{Thm1} shows that the set of spherical configurations attracts at least a neighbourhood of itself under (\ref{EQWD}), but the existence of globally attracting sets, as studied in \cite{Ladyzhenskaya1991}, is not known.

On the other hand, hyperbolic equations relating to mean curvature problems have been studied by numbers of mathematicians. (\ref{EQClassical}) is a prototype of such equations. Besides the issue of degeneracy addressed in last section, there is yet another significant difference between (\ref{EQ0}) and (\ref{EQClassical}). Klainerman was the first mathematician to obtain global-in-time results for equation (\ref{EQClassical}) with $n\geq6$: in \cite{Klainerman1980}, Klainerman used a modified Nash-Moser technique. Afterwards, Klainerman proved a similar result without using Nash-Moser technique \cite{Klainerman1982}, around the same time as Shatah \cite{Shatah1982}. Finally, Klainerman developed the well-known vector field method for $n\geq4$ in the groundbreaking work \cite{Klainerman1985}. Nevertheless, all the above results strongly rely on the dispersive properties possessed by wave equations on $\mathbb{R}^n$: the original Nash-Moser technique in \cite{Klainerman1980} and later works \cite{Klainerman1982} \cite{Shatah1982} made use of the following dispersive estimate for solution to the Cauchy problem $\partial_t^2u-\Delta_xu=0$, $u(0)=0$, $u'(0)=g$ in $\mathbb{R}_t\times\mathbb{R}_x^n$:
$$
|\nabla u(x,t)|\leq\frac{1}{(1+|t|)^{(n-1)/2}}\|g\|_{W^{1,n}}.
$$
The Klainerman-Sobolev inequality in \cite{Klainerman1985} can be considered as a refined version of this estimate. The global well-posedness results obtained for quasilinear wave equations on $\mathbb{R}^n$ required these dispersive estimates, but \emph{no such inequality can hold when the underlying space is compact}. Intuitively, we imagine that a wave should disperse in the Euclidean space if there are more than one spatial direction; although the total energy is conserved, the amplitude of the wave around a given location should decay because the dispersion takes energy away. But if the underlying space is compact, the wave should travel back within a finite time and keeps recurring, so there cannot exist any decay unless some dissipative feature enters into the system. This is a general difficulty for dispersive systems defined on compact spaces, and the corresponding global-in-time theory is thus expected to be very different from the known results for the Euclidean setting.

Finally, Theorem \ref{Thm2} shows that the damping-free equation (\ref{EQ0}) has lifespan $\simeq\varepsilon^{-1/6}$ for initial data $\varepsilon$-close to the unit sphere. The question of \emph{what happens next} is an exciting one. For equations of type (\ref{EQClassical}) with $n=2$, small-data blow-up results were initially obtained by Alinhac \cite{Alinhac1999}, and it is very reasonable to conjecture that (\ref{EQWD}) also develops a ``cusp" type blow-up point, just as its Euclidean prototype (\ref{EQClassical}) does. A heuristic is that the absence of dispersion on compact spaces should reinforce the intensity of blow-up, but this is yet beyond our knowledge.

On the other hand, since (\ref{EQ0}) is a Hamiltonian system, it is also reasonable ask whether some specific global solutions could exist. As pointed out by Notz \cite{Notz2010}, spheres with radius $r(t)$ solving the ODE
$$
r''=-2r+\frac{2}{r}
$$
form a family of global periodic solutions of (\ref{EQ0}), foliating the phase space $\{(r,z)\in\mathbb{R}^2:r>0\}$. These solutions could be referred as \emph{breathers}. Some natural questions then arises:
\begin{itemize}
    \item Do these breather solutions persist under perturbation, i.e. are the perturbed solutions still global and periodic?
    \item Does there exist any quasi-periodic solution close to a solution to the linearized equation around a breather solution?
\end{itemize}

For some nonlinear dispersive Hamiltonian PDEs defined on line segments or tori, these questions have been partially answered. Typical examples include the semilinear string vibration equation, periodic water wave equation or semi-linear Schr\"{o}dinger equation. For these systems, either a scheme of finding periodic solutions is valid (see e.g. Berti's monograph \cite{Berti2007}), or a KAM type theory could be established (see e.g. recent work of Baldi, Berti, Haus and Montalto \cite{BBHM2018}) and a Birkhoff normal form reduction is possible (see e.g. \cite{BFP2018}), which extends the estimate on the lifespan. However, to the author's knowledge, all such results are concerned with either segments or tori -- in fact, mostly the circle  $S^1$, except for works following Bourgain's paper \cite{Bourgain1998}. Either extension or disproof of these results for Hamiltonian PDEs on compact curved manifolds, for example (\ref{EQ0}) on sphere, has not been discovered yet.

\textbf{Acknowledgment.} This paper will be a part of the author's doctoral thesis. The author would like to thank professor Gigliola Staffilani for supervision and inspiration of the work, and would also like to thank professor Peter Hintz and Andrew Lawrie for providing useful ideas, suggestions and comments. The author is also grateful for the comments from the referees on physical aspects of the problem.

\section{Functional Settings}\label{Sec2}
\subsection{Notation}\label{Subsec2.1}
Throughout the paper, if $w:S^2\to\mathbb{R}^3$ is a smooth embedding, we shall use $N(w)$ to denote its outward pointing normal vector field, $g(w)$ to denote the induced Riemannian metric on $S^2$ (and $\mu(w)$ to denote the induced surface measure), $h(w)$ to denote the second fundamental form of the embedding and $H(w)$ to denote the (scalar) mean curvature. We shall also write, for a general vector $A\in\mathbb{R}^3$,
$$
\perp_{w}A:=N(w)\cdot A,
$$
$$
\top_wA:=A-[N(w)\cdot A]N(w).
$$
They are respectively the (scalar) normal projection and (vectorial) tangent projection along $w$, and $\top_w$ is in fact a bundle map from the trivial bundle $\mathbb{R}^3$ to the (extrinsic) tangent bundle $\bigsqcup_{x}T_{w(x)}S^2$.

Unless specified otherwise, all tensorial and functional norms on $S^2$ will be taken with respect to a fixed atlas, i.e. the most commonly used two-disk covering $\{B_1,B_2\}$ of $S^2$, and the standard spherical metric $g_0=g(i_0)$. For example, the $C^n$ norm of a vector field $X$ on $S^2$ is defined to be
$$
|X|_{C^n}
:=\sum_{j=1}^2\sum_{\alpha:|\alpha|\leq n}\left|D^\alpha_x X|_{B_j}\right|_{L^\infty_x(B_j)},
$$
and the Sobolev norm of a smooth function $f$ is defined to be
$$
\|f\|_{W^{n,p}_x(g_0)}
:=\left(\int_{S^2}|f|^pd\mu_0\right)^{1/p}+\sum_{j=1,2}\left(\sum_{\alpha:1\leq|\alpha|\leq n}\int_{B_j}|D^\alpha_x f|^pd\mu_0\right)^{1/p}.
$$
When $p=\infty$, the modification to Lipschitz or $C^n$ norms is obvious. These are not intrinsically defined tensorial norms, but it is easy to verify that they are equivalent to the intrinsic norms defined via the Riemannian connection of $g_0$. On the other hand, we also define norms of time-derivatives on a given time-slice as follows: for $f\in C^k([0,T];W^{n,p}_x(g_0))$, we set
$$
\|f(t)\|_{W^{n,p}_x(g_0)}^{(k)}
:=\sum_{l=0}^k\|\partial_t^lf(t)\|_{W^{n,p}_x(g_0)}.
$$
Obviously
$$
\|f(t)\|_{W^{n,p}_x(g_0)}^{(k)}
\leq \|f(t_0)\|_{W^{n,p}_x(g_0)}^{(k)}
+\int_{t_0}^t\|f(s)\|_{W^{n,p}_x(g_0)}^{(k+1)}ds.
$$
Furthermore, we quote a lemma of Hamilton on time derivatives:
\begin{lemma}[Hamilton's trick, see e.g. lemma 2.1.3. of \cite{Mantegazza2011}]
Let $f$ be a real Lipschitz function on $[0,T]\times S^2$. Define $f_{\max}(t):=\max_{x\in S^2} f(t,x)$. Then $f_{\max}(t)$ is Lipschitz in $t$, and for almost all $t\in[0,T]$,
$$
\frac{d}{dt}f_{{\max}}(t)=\frac{\partial f}{\partial t}(t,x),
$$
where $x\in S^2$ is any point such that $f_{\max}(t)=f(t,x)$.
\end{lemma}
With the aid of this lemma, we find in fact $\|f(t)\|_{C_x^n}^{(k)}$ is Lipschitz continuous in $t$, and
\begin{equation}\label{HamiltonMax}
\left|\frac{d}{dt}\|f(t)\|_{C_x^n}^{(k)}\right|
\leq \|f(t)\|_{C_x^n}^{(k+1)},\,\text{a.e. }t.
\end{equation}

The graded space to be considered will be as follows. Let $\mathfrak{X}$ be the Fréchet  space of all smooth tangent vector fields on $S^2$, with the grading given by H?lder norms $|\cdot|_{C^s}$. For a fixed number $\beta>0$ and a fixed integer $k\geq0$, define $\mathbf{E}_{\beta,k}^n$ to be the collection of all $C^3$ mappings from $[0,\infty)$ to $C^\infty(S^2;\mathbb{R}^3)$ such that for any $n\geq0$, the norm
$$
\|u\|_{\beta,k;n}
:=\sup_{t\geq0}e^{\beta t}\|u\|_{H^n_x(g_0)}^{(k)}
=\sup_{t\geq0}e^{\beta t}\sum_{l=0}^k\|\partial^l_tu(t)\|_{H^n_x(g_0)}
$$
is finite, and further set $\mathbf{E}_{\beta,k}=\cap_{n\geq0}\mathbf{E}_{\beta,k}^n$, equipped with the natural Fréchet space topology.  Define the Fréchet  space $\mathbf{F}_\beta=\mathfrak{X}\oplus\mathbb{R}^3\oplus\mathbf{E}_{\beta,3}$, with grading
$$\|(X,a,u)\|_{n}:=\|X\|_{H^{n}}+|a|+\|u\|_{\beta,3;n},$$
where the norms are taken with respect to the fixed metric $g_0=g(i_0)$.

The smoothing operator on $\mathbf{F}_\beta$ is defined as follows: if $f$ maps $[0,T]\times S^2$ to $\mathbb{R}^3$, then define
\begin{equation}\label{smoothing0}
S_\theta f:=\sum_{\lambda\in\sigma[-\Delta_{g_0}]:\lambda\leq\theta}\mathcal{Q}_\lambda f,
\end{equation}
where $\mathcal{Q}_\lambda$ is the eigenprojection corresponding to eigenvalue $\lambda$ of $\Delta_{g_0}$ on $L^2(g_0)$. The operator only acts on spatial variables and commutes with $\partial_t$. If $X\in\mathfrak{X}$, then we define
\begin{equation}\label{smoothing1}
S_\theta X
:=\sum_{\lambda\in\sigma(D_{g_0}^*D_{g_0}),\lambda\leq\theta}\mathcal{Q}^{(1)}_\lambda X,
\end{equation}
where $D_{g_0}^*D_{g_0}$ is the trace Laplacian operator corresponding to $g_0$, and $\mathcal{Q}^{(1)}_\lambda$ is the spectral projection to the eigenspace corresponding to $\lambda$ in the space of $L^2$ vector fields (with respect to $g_0$). It is well-known that $D_{g_0}^*D_{g_0}$ is self-adjoint and positive definite (there is no zero eigenvalue since there is no parallel vector field on $S^2$).

The space in which solutions live in will be modelled on this $\mathbf{F}_\beta$. In particular, the diffeomorphism group $\mathfrak{Diff}(S^2)$ is a Fréchet  Lie group modelled on $\mathfrak{X}$; however, as pointed out in \cite{KM1997}, the exponential map (time 1 flow map) from $\mathfrak{X}$ to $\mathfrak{Diff}(S^2)$ is not locally surjective, and the local diffeomorphism from $\mathfrak{X}$ to $\mathfrak{Diff}(S^2)$ has to be constructed under a given Riemannian metric. We thus follow section 42 of \cite{KM1997}: choose the metric to be $g_0$ and define
\begin{equation}\label{EX}
\mathcal{E}_X(x):=\exp_x^{g_0}(X(x)).
\end{equation}
Given $X\in\mathfrak{X}$, the diffeomorphism $\mathcal{E}_X\in\mathfrak{Diff}(S^2)$ can be computed explicitly by considering the distorted embedding $i_0\circ\mathcal{E}_X:S^2\to\mathbb{R}^3$. We know that the geodesic on the embedded unit sphere starting at a point $p\in S^2\subset\mathbb{R}^3$ along direction $v$ is parameterized by its arc length as
$$
t\to\cos(|v|t)p+\sin(|v|t)\frac{v}{|v|},
$$
so if $X\in\mathfrak{X}$, then
\begin{equation}\label{IEX}
(i_0\circ\mathcal{E}_X)(x)=(\cos|X(x)|_{g_0})i_0(x)+\frac{\sin|X(x)|_{g_0}}{|X(x)|_{g_0}}(di_0)_xX(x),
\end{equation}
where $|X(x)|_{g_0}$ is the norm of $X(x)$ with respect to the metric $g_0$. Consequently, $i_0\circ\mathcal{E}_X$ is smooth (in fact, analytic) in $X$ only and does not depend on any derivative of $X$. We note that the outward normal vector field $N(i_0\circ\mathcal{E}_X)=i_0\circ\mathcal{E}_X$, so $N(i_0\circ\mathcal{E}_X)$ has the same regularity property as $X$ does. It is thus easily verified that this gives a diffeomorphism from some $C^2$-neighbourhood of $0\in\mathfrak{X}$ to some neighbourhood of $\text{id}\in\mathfrak{Diff}(S^2)$, and for $X$ in that $C^2$-neighbourhood of $0\in\mathfrak{X}$, we have
$$
\|i_0\circ\mathcal{E}_X-\text{id}\|_{W^{n,p}(g_0)}\leq C_{n,p}\|X\|_{W^{n,p}(g_0)}.
$$
For a diffeomorphism $\varphi\in\mathfrak{Diff}(S^2)$, we shall write $i_\varphi$ for $i_0\circ\varphi$, and let $g_\varphi:=g(i_\varphi)$, $\mu_\varphi:=\mu(i_\varphi)$.  Throughout the paper, unless otherwise noted, we write $\varphi=\mathcal{E}_X$ for some $X\in\mathfrak{X}$.

We then define a mapping $\Xi:\mathbf{F}_\beta\to C^\infty([0,\infty)\times S^2;\mathbb{R}^3)$ by
$$
\Xi(X,a,u)=i_0\circ\mathcal{E}_X+a+u
=i_\varphi+a+u.
$$
Intuitively this means that we will be dealing with perturbations near the static (spherical) configuration, where the diffeomorphisms of the static configuration itself is reduced by the term $\mathcal{E}_X$ (``fixing a gauge", or fixing the terminal embedding). The vector $a\in\mathbb{R}^3$ represents the shift of the center, and it has no effect on the configuration. We also calculate
$$
\begin{aligned}
\Xi'(X,a,u)(Y,c,v)(t,x)
&=(di_0)_{\exp_x^{g_0}(X(x))}\left(d\exp_x^{g_0}\right)_{X(x)}Y(x)+c+v(t,x).
\end{aligned}
$$
We define, for $\varphi=\mathcal{E}_X$,
\begin{equation}\label{SigmaPhi}
(\Sigma_\varphi Y)(x):=(di_0)_{\exp_x^{g_0}(X(x))}\left(d\exp_x^{g_0}\right)_{X(x)}Y(x).
\end{equation}
Note that $\Sigma_\varphi Y$ is a tangent vector field along the embedding $i_\varphi:S^2\to\mathbb{R}^3$, i.e., for any $x\in S^2$, we have $(\Sigma_\varphi Y)(x)\in T_{i_\varphi(x)}S^2$. In other words, $\Sigma_\varphi $ is a bundle isomorphism from the (intrinsically defined) tangent bundle $T(S^2)$ to  $\bigsqcup_{x}T_{i_\varphi(x)}S^2$, the (exterior) tangent bundle of the embedding $i_\varphi$.

We will be working with a neighbourhood $\mathfrak{U}=\mathfrak{U}_0\times U\times \mathbf{V}$ of $0$ in $\mathbf{F}_\beta$, such that for any $(X,a,u)\in\mathfrak{U}$, $\varphi=\mathcal{E}_X$ is a diffeomorphism of $S^2$, and $i_\varphi+u:S^2\to\mathbb{R}^3$ is always a smooth embedding. These neighbourhoods will be specified in the context.

Now write $\Phi=\Psi\circ\Xi$, and substitute $u\to\Xi(X,a,u)$, $v\to\Psi(Y,c,v)$ in (\ref{LEQ1}). Write
$$
w:=\Xi(X,a,u)=i_\varphi+a+u,
$$
and
$$
\eta:=\Xi'(X,a,u)(Y,c,v)=\Sigma_\varphi Y+c+v.
$$
Then $\Phi:\mathbf{F}_\beta\to\mathbf{E}_{\beta,1}$ is a tame mapping, and
\begin{equation}\label{LEQ2}
\begin{aligned}
\Phi'(X,a,u)&(Y,c,v)\\
&=\Psi'(\Xi(X,a,u))[\Xi'(X,a,u)(Y,c,v)]\\
&=\frac{\partial^2\eta}{\partial t^2}+b\frac{\partial \eta}{\partial t}\\
&\quad-\frac{d\mu(w)}{d\mu_0}\left(\Delta_{g(w)}[\eta\cdot N(w)]+|h(w)|^2[\eta\cdot N(w)]\right)N(w)\\
&\quad-\frac{d\mu(w)}{d\mu_0}\left(\frac{\kappa}{\mathrm{Vol}(w)^2}\int_{S^2}[\eta\cdot N(w)] d\mu(w)+[\nabla^{g(w)} H(w)\cdot\top_{w}\eta]\right)N(w)\\
&\quad+\frac{d\mu(w)}{d\mu_0}\left(-H(w)+\frac{\kappa}{\mathrm{Vol}(w)}\right)
(\text{div}^{g(w)}{\top_{w}}\eta+H(w)[\eta\cdot N(w)])N(w)\\
&\quad+\frac{d\mu(w)}{d\mu_0}\left(-H(w)+\frac{\kappa}{\mathrm{Vol}(w)}\right)(\nabla^{g(w)}[\eta\cdot N(w)]-h^{kl}(w)[ \eta\cdot\partial_l(w)]\partial_k(w)).
\end{aligned}\tag{LEQ}
\end{equation}

Following a general Nash-Moser scheme, this paper will be mostly devoted to the study of two linearized equations. The first one is
\begin{equation}\label{LS1}
\Psi'(w)\eta=f,\tag{LEQ1}
\end{equation}
where $w\in C^{3}([0,T];C^\infty(S^2;\mathbb{R}^3))$ and $f\in C^{1}([0,T];C^\infty(S^2;\mathbb{R}^3))$ are known and $\eta$ is what needs to be solved. The second one is
\begin{equation}\label{LS2}
\Phi'(X,a,u)(Y,c,v)=f,\tag{LEQ2}
\end{equation}
where $(X,a,u)\in\mathbf{F}_\beta$ and $f\in\mathbf{E}_{\beta,1}$ are known and $(Y,c,v)$ is what needs to be solved.

\subsection{Function Theoretic Lemmas}
We begin with a multi-linear interpolation inequality for derivatives.
\begin{lemma}\label{Interpolation}
Fix integers $n,k,s\geq1$ and $k_1,\cdots,k_s$ with $k<n$, $k=k_1+\cdots+k_s$. For any $v_1,\cdots,v_s,w\in C_0^\infty(\mathbb{R}^l)$, any $\varepsilon>0$ and any $p\in(1,\infty]$, there is a constant $C_{n,s,p}$ depending on $n,s,p$ only such that
$$
\|D^{k_1}_xv_1\cdots D^{k_s}_xv_s\cdot D^{n-k}_x w\|_{L^p_x}
\leq \varepsilon\|D^n_xw\|_{L^p_x}\prod_{j=1}^s|v_j|_{L^\infty_x}^{(n-k_j)/(n-k)}
+C_{n,s,p}\varepsilon^{-(n-k)/k}|D^n_xv|_{L^\infty_x}\|w\|_{L^p_x}.
$$
Here $v=(v_1,...,v_s)$.
\end{lemma}
\begin{proof}
We quote a pointwise interpolation inequality on derivatives from \cite{MS1999} by Maz'ya and Shaposhnikova: for any test function $v\in C_0^\infty(\mathbb{R}^l)$ and any $x\in\mathbb{R}^l$,
$$
|D^k_xv(x)|\leq C_n|Mv(x)|^{1-k/n}|MD^n_xv(x)|^{k/n},\,0\leq k\leq n.
$$
where $M$ is the Hardy-Littlewood maximal function. Thus, by Young's inequality,
$$
\begin{aligned}
|D^{k_1}_xv_1(x)&\cdots D^{k_s}_xv_s(x)\cdot D^{n-k}_x w(x)|\\
&\leq C_{n,s}\prod_{j=1}^s|Mv_j(x)|^{1-k_j/n}\cdot
\prod_{j=1}^s|MD^n_xv_j(x)|^{k_j/n}\cdot
|Mw(x)|^{k/n}|MD^n_xw(x)|^{1-k/n}\\
&\leq\varepsilon\left(\prod_{j=1}^s|v_j|_{L^\infty_x}^{1-k_j/n}\right)^{n/(n-k)}|MD^n_xw(x)|
+C_{n,s}\varepsilon^{-(n-k)/k}\left(\prod_{j=1}^s|D^n_xv_j|_{L^\infty_x}^{k_j/n}\right)^{n/k}|Mw(x)|\\
&\leq\varepsilon\left(\prod_{j=1}^s|v_j|_{L^\infty_x}^{(n-k_j)/(n-k)}\right)|MD^n_xw(x)|
+C_{n,s}\varepsilon^{-(n-k)/k}\left(\sum_{j=1}^s|D^n_xv_j|_{L^\infty_x}\right)|Mw(x)|.
\end{aligned}
$$
The desired result is obtained by integrating over $\mathbb{R}^l$ and applying the strong $(p,p)$ ($p\in(1,\infty]$) property of the Hardy-Littlewood maximal function.
\end{proof}

The following lemma is a minor modification of lemma 5.1. in Klainerman's paper \cite{Klainerman1980}, and is proved using Fa\`{a} di Bruno's formula:
\begin{lemma}\label{MoserComposition}
Let $F:\mathbb{R}^q\to\mathbb{R}^r$ be a smooth mapping in its arguments. For any $n\geq1$ and any $v\in C_0^\infty(\mathbb{R}^l;\mathbb{R}^q)$ supported in the unit ball $B_1$, we have
$$
\|\nabla^n(F\circ v)\|_{L^p_x}\leq C(n,p,F,|v|_{L^\infty_x})\|\nabla^nv\|_{L^p_x},
$$
where the constant $C(n,p,F,|v|_{L^\infty_x})$ depends on $n,p$, the differentials of $F$ up to order $n$, and is monotonically increasing in the norm $|v|_{L^\infty_x}$.
\end{lemma}

Combining \ref{Interpolation} and \ref{MoserComposition}, we obtain the following estimate on composition of differential operators:
\begin{lemma}\label{OperatorComposition}
Let $k=k_1+...+k_N$. Suppose $u\in C^{l+k}(\mathbb{R}^l)$, and suppose
$$
A^{(i)}
=\sum_{\alpha:|\alpha|\leq k_i}A^{(i)}_\alpha(x,u,D_xu,...,D_x^lu)D_x^\alpha,\,
i=1,..,N
$$
are linear differential operators, with coefficients smooth functions in $x$. Then for any $v\in C_0^\infty(\mathbb{R}^l)$ supported in the unit ball $B_1$, we have
$$
\left\|A^{(1)}\circ...\circ A^{(N)}v\right\|_{L^2_x}
\leq C(k,l,[A^{(i)}_\alpha]_{i,\alpha},|u|_{C^l_x})\left(
\|v\|_{H^k_x}+(1+|u|_{C^{l+k}_x})\|v\|_{L^2_x}
\right).
$$
\end{lemma}

Next we present a lemma on perturbation of differential operators; the proof relies on the first order Taylor formula and previous lemmas.
\begin{lemma}\label{PerturbOperator}
If $A(x;u,v)$ is any $k$-th order differential operator whose coefficients depend smoothly on
$$
\left(\partial_t^jD_x^lu,D_x^mv\right),\,
j\leq j_0,\,l\leq l_0,\,m\leq m_0,
$$
and vanishes if $u=0$, then for any $f\in C_0^\infty(\mathbb{R}^l;\mathbb{R}^q)$ supported in the unit ball $B_1$, we have
$$
\begin{aligned}
\|A(x;u,v)f\|_{H^n_x}
&\leq C\left(n,|u|_{C_x^{l_0}}^{(j_0)},|v|_{C^{m_0}}\right)|u|_{C_x^{n+l_0}}^{(j_0)}\|f\|_{L^2_x}\\
&\quad+C\left(n,|u|_{C_x^{l_0}}^{(j_0)},|v|_{C^{m_0}}\right)|u|_{C_x^{l_0}}^{(j_0)}\left[\left(1+|u|_{C_x^{n+l_0}}^{(j_0)}+|v|_{C^{n+m_0}}\right)\|f\|_{L^2_x}
+\|f\|_{H^{n+k}_x}\right].
\end{aligned}
$$
\end{lemma}

Now we derive some function theoretic results on evolving surfaces.

Suppose $\eta:[0,\infty)\times S^2\to\mathbb{R}^3$ is a smooth mapping. For a fixed $(X,a,u)\in\mathbf{F}_\beta$ such that $w=i_\varphi+a+u$ is a time-dependent smooth embedding, we decompose
$$\eta=\perp_w\eta+\top_w\eta=\phi N(w)+\psi.$$

If $(X,a,u)\in\mathbf{F}_\beta$ is such that $w:=i_\varphi+a+u$ is still a smooth embedding from $S^2$ to $\mathbb{R}^3$, then at each $x\in S^2$, the tangent space is spanned by
$$
\partial_iw(x)=\partial_ii_\varphi(x)+\partial_iu(t,x),\,i=1,2,
$$
so the tangent bundle depends smoothly on first order derivatives of $X$ and $u$, and
$$
N(w)=\frac{\partial_1w\times\partial_2w}{|\partial_1w\times\partial_2w|}
$$
is smooth in first order derivatives of $X$ and $u$. The following lemma, whose proof is a mere application of Lemma \ref{Interpolation} and Lemma \ref{MoserComposition}, provides the regularity information of geometric quantities on an evolving surface:

\begin{lemma}\label{GeometricQuantities}
Write $w=i_\varphi+a+u$, where $\varphi=\mathcal{E}_X$. There is a constant $\delta_0>0$ such that if $|X|_{C^2}+|u|_{C^2_x}<\delta_0$, then the following tame estimates hold:

(1) Define $\Sigma_\varphi$ as in (\ref{SigmaPhi}). For a tangent vector field $\psi$ along the embedding $i_\varphi$, we have
$$
|\Sigma_\varphi ^{-1}\psi|_{C^n}\leq C_n\left(|\psi|_{C^n_x}+|X|_{C^{n+1}}|\psi|_{L^\infty_x}\right).
$$

(2) Given any smooth mapping $\eta:S^2\to\mathbb{R}^3$, we have, for $p\in(1,\infty]$,
$$
\|\perp_w\eta\|_{W^{n,p}_x}+\|\top_w\eta\|_{W^{n,p}_x}
\leq C_{n,p}\left[\|\eta\|_{W^{n,p}_x(g_0)}
+\left(|X|_{C^{n+1}}+\|u\|_{W^{n+1,p}_x(g_0)}\right)|\eta|_{L^\infty_x}\right].
$$

(3) Under the decomposition $\eta=\perp_w\eta+\top_w\eta=\phi N(w)+\psi$, we have, for $p\in(1,\infty]$,
$$
\begin{aligned}
\|\eta\|_{W^{n,p}_x(g_0)}
\leq\|\psi\|_{W^{n,p}_x(g_0)}
+C_{n,p}\left[\|\phi\|_{W^{n,p}_x(g_0)}
+\left(|X|_{C^{n+1}}+\|u\|_{W^{n+1,p}_x(g_0)}\right)|\phi|_{L^\infty_x}\right].
\end{aligned}
$$

(4) The induced metric $g(w)$ and the second fundamental form $h(w)$ satisfy, for $p\in(1,\infty]$,
$$
\|g(w)\|_{W^{n,p}_x(g_0)}\leq C_{n,p}\left(1+|X|_{C^{n+1}}+\|u\|_{W^{n+1,p}_x(g_0)}\right),
$$
$$
\|h(w)\|_{W^{n,p}_x(g_0)}\leq C_{n,p}\left(1+|X|_{C^{n+2}}+\|u\|_{W^{n+2,p}_x(g_0)}\right).
$$
\end{lemma}

\section{Tame Estimate I: Elliptic Operators and Decay Estimates}\label{Sec3}
\subsection{Spectral Properties}
We need to investigate some elliptic operators on an evolving surface. For a vector field $X\in\mathfrak{X}$ that is $C^2$-close to 0, we consider the corresponding diffeomorphism $\varphi\in\mathfrak{Diff}(S^2)$ given by $\varphi=\mathcal{E}_X$. Define
$$
L_\varphi\phi:=\frac{d\mu_\varphi}{d\mu_0}\left(\Delta_{g(i_\varphi)}\phi+2\phi-\frac{6}{4\pi}\int_{S^2}\phi d\mu_\varphi\right);
$$
still set $w=i_\varphi+a+u$, and define a more general operator
$$
\begin{aligned}
&L(\varphi,u)\phi\\
&:=
\frac{d\mu(w)}{d\mu_0}\left[\Delta_{g(w)}\phi+|h(w)|^2\phi+\left(-H(w)+\frac{\kappa}{\mathrm{Vol}(w)}\right)H(w)\phi
-\frac{\kappa}{\text{Vol}(w)^2}\int_{S^2}\phi d\mu(w)\right]-|\partial_tN(w)|^2\phi.
\end{aligned}
$$
Further, for a section $\theta$ of $T^*(S^2)$ (the cotangent bundle), define the elliptic operator
$$
L^1(\varphi,u)\theta
:=\frac{d\mu(w)}{d\mu_0}\left[-\Delta_{g(w)}^1\theta+|h(w)|^2\theta+\left(-H(w)+\frac{\kappa}{\mathrm{Vol}(w)}\right)H(w)\theta\right]-|\partial_tN(w)|^2\theta,
$$
where $\Delta_{g(w)}^1$ is the (positive) Hodge-Laplacian with respect to $g(w)$ acting on 1-forms. It is easily seen that $L(\varphi,u)$ is self-adjoint on $L^2(g_0)$ with domain $H^2(g_0)$, and in particular, $L_\varphi$ is a non-positive self-adjoint operator on $L^2(g_0)$. Furthermore, $L^1(\varphi,u)$ is self-adjoint on $H^2$-sections of 1-forms with respect to the metric $g(w)$.

We have the following proposition concerning the spectral properties of the above elliptic operators:
\begin{proposition}\label{SpectralProp}
(A) Set
$$
M_\varphi=\sup \left|\frac{d\mu_\varphi}{d\mu_0}\right|,\,m_\varphi=\inf \left|\frac{d\mu_\varphi}{d\mu_0}\right|,
$$
then $m_\varphi\leq 1\leq M_\varphi$, and the point $0$ belongs to $\sigma[L_\varphi]$ with multiplicity $3$, and the corresponding eigenfunctions are $\{N^k(i_\varphi)\}_{k=1}^3$. Furthermore, all non-zero eigenvalues of $L_\varphi$ are all less than $-4m_\varphi$.

(B) Suppose $\varphi=\mathcal{E}_X$. There exist constants $A,\delta_0>0$ such that if, at a given time slice there holds $|X|_{C^2}+|u|_{C_x^2}^{(1)}<\delta_0$, then there are three spectral points of $L(\varphi,u)$ (counting multiplicity) with magnitude less than $A|u|_{C_t^1C_x^2}$, and all other spectral points are in the interval $\left(-\infty,-2\right]$. If we denote by $\mathcal{P}_0(\varphi,u)$ the spectral projection corresponding to the eigenvalues close to 0 and write $\mathcal{P}_\infty(\varphi,u):=\text{Id}-\mathcal{P}_0(\varphi,u)$, then both $\mathcal{P}_0(\varphi,u),\mathcal{P}_\infty(\varphi,u)$ are analytic in $u$ and $\varphi$ in the sense that they can be represented as norm-convergent operator-valued power series of $[\partial^j_t(\nabla^{g_0})^lu]_{0\leq j\leq1}^{0\leq l\leq2}$ and $[(\nabla^{g_0})^lX]_{0\leq l\leq2}$.

(C) Under the similar assumption of (B), there are three spectral points of $L^1(\varphi,u)$ with magnitude less than $A|u|_{C_x^2}^{(1)}$, and all other spectral points are in the interval $\left(-\infty,-2\right]$, and the corresponding properties for spectral projections are still valid.
\end{proposition}
\begin{proof}
(A) The assertion $m_\varphi\leq 1\leq M_\varphi$ is a direct consequence of the equality
$$
\int_{S^2}d\mu_0=\int_{S^2}d\mu_\varphi=4\pi.
$$

Let us now denote by $\mathcal{P}_\lambda^\varphi$ the eigenprojection of $-L_\varphi$ corresponding to eigenvector $\lambda$.

Geometric identities
$$
\Delta_{g(i_\varphi)} N(i_\varphi)+2N(i_\varphi)=0
$$
and
$$
\int_{S^2}N(i_\varphi) d\mu_\varphi=0
$$
ensure that $\{N^k(i_\varphi)\}_{k=1}^3$ are eigenfunctions of $L_\varphi$ with eigenvalue 0. Taking $L^2(g_0)$ inner product, applying the standard spectral theory of the Laplacian $\Delta_{g(i_\varphi)}$ (which is non-positive and self-adjoint on $L^2(\mu_\varphi)$), we obtain
$$
\begin{aligned}
-\langle L_\varphi\phi,\phi\rangle_{L^2(g_0)}
&=-\langle\Delta_{g(i_\varphi)}\phi,\phi\rangle_{L^2(\mu_\varphi)}-2\langle\phi,\phi\rangle^2_{L^2(\mu_\varphi)}+\frac{6}{4\pi}|\langle\phi,1\rangle_{L^2(\mu_\varphi)}|^2\\
&=-\sum_{\lambda\in\sigma[\Delta_{g(i_\varphi)}]}(\lambda+2)\|\mathcal{Q}^\varphi_\lambda\phi\|_{L^2(\mu_\varphi)}^2+6\|\mathcal{Q}^\varphi_0\phi\|_{L^2(\mu_\varphi)}^2\\
&\geq4\left\|(1-\mathcal{Q}^\varphi_{-2})\phi\right\|_{L^2(\mu_\varphi)}^2.
\end{aligned}
$$
Here $\mathcal{Q}^\varphi_\lambda$ denotes the spectral projection of $\Delta_{g(i_\varphi)}$ on $L^2(\mu_\varphi)$ corresponding to eigenvalue $\lambda$. So $\mathcal{Q}^\varphi_{-2}$ is a projection operator onto $\text{span}\{N^k(i_\varphi)\}_{k=1}^3$. The right-hand-side is thus strictly greater than zero unless $\phi\in\text{span}\{N^k(i_\varphi)\}_{k=1}^3$. Then Rayleigh's formula ensures that all non-zero spectral points of $L_\varphi$ are negative. Furthermore, Rayleigh's formula gives that the greatest non-zero eigenvalue of $L_\varphi$ is
$$
\begin{aligned}
\inf_{\phi\in L^2(g_0)}\frac{-\langle L_\varphi\phi,\phi\rangle_{L^2(g_0)}}{\left\|(1-\mathcal{P}^\varphi_0)\phi\right\|_{L^2(g_0)}^2}
&\geq4\inf_{\phi\in L^2(g_0)}\frac{\left\|(1-\mathcal{Q}^\varphi_{-2})\phi\right\|_{L^2(\mu_\varphi)}^2}{\left\|(1-\mathcal{P}^\varphi_0)\phi\right\|_{L^2(g_0)}^2}.
\end{aligned}
$$
Let us calculate the quotient on the right-hand-side. The numerator is estimated as
$$
\begin{aligned}
\left\|(1-\mathcal{Q}^\varphi_{-2})\phi\right\|_{L^2(\mu_\varphi)}^2
&=\int_{S^2}|(1-\mathcal{Q}^\varphi_{-2})\phi|^2d\mu_\varphi\\
&\geq\inf\left|\frac{d\mu_\varphi}{d\mu_0}\right|\cdot\int_{S^2}|(1-\mathcal{Q}^\varphi_{-2})\phi|^2d\mu_0
=m_\varphi\left\|(1-\mathcal{Q}^\varphi_{-2})\phi\right\|_{L^2(g_0)}^2.
\end{aligned}
$$
Since $(1-\mathcal{P}^\varphi_0)$ is in fact the orthogonal projection in $L^2(g_0)$ onto $\text{span}\{N^k(i_\varphi)\}_{k=1}^3$, it follows from the extremal property of orthogonal projections that $\left\|(1-\mathcal{Q}^\varphi_{-2})\phi\right\|_{L^2(g_0)}^2/\left\|(1-\mathcal{P}^\varphi_0)\phi\right\|_{L^2(g_0)}^2\geq1$. This finishes the proof of (A).

(B) The proof actually follows from a standard perturbation argument for type (B) analytic family of self-adjoint operators in the sense of \cite{Kato2013}; we shall modify the methodology presented in \cite{Kato2013}, Chapter VII, Section 4, since the result is not directly applicable due to some technical reasons.

Let $y=(y_1,y_2)$ be a real parameter of small magnitude. First note that there is a small constant $c_1$ such that if $|y|<c_1(\|u\|_{C^1_x}+|X|_{C^1})^{-1}$, then $i_0\circ\mathcal{E}_{y_1X}+y_2u$ is a $C^1$ embedding. Let us then consider the sesquilinear form associated to $L(\mathcal{E}_{y_1X},y_2u)$: with $w_y=i_0\circ\mathcal{E}_{y_1X}+y_2u$, we define
$$
\begin{aligned}
B(y_1X,y_2u)[\phi]
&:=\int_{S^2}\left(|\nabla^{g(w_y)}\phi|^2-|h(w_y)|^2|\phi|^2\right)d\mu(w_y)
+\frac{\kappa}{\text{Vol}(w_y)^2}\left|\int_{S^2}\phi d\mu(w_y)\right|^2\\
&\quad+\int_{S^2}\left(H(w_y)-\frac{\kappa}{\mathrm{Vol}(w_y)}\right)H(w_y)|\phi|^2d\mu(w_y)+\int_{S^2}|\partial_tN(w_y)|^2|\phi|^2 d\mu_0.
\end{aligned}
$$
Write $B_0$ for $B(0,0)$. It is very easily verified that $B(y_1X,y_2u)$ gives rise to an analytic family of forms of type (A) in the variable $y$: the above is a real analytic function in $y$, and each coefficient in the power series expansion in $y$ is a symmetric form. A direct calculation also gives that there is a constant $c_2$ independent of $\varphi$ such that
$$
B(y_1X,y_2u)[\phi]
\geq-c_2|y|\left(|X|_{C^2}+|u|_{C_x^2}^{(1)}\right)\|\phi\|_{H^1_x(g_0)}^2,
$$
so
$$
\sigma[L(\mathcal{E}_{y_1X},y_2u)]
\subset\left(-\infty,\,
c_2|y|\left(|X|_{C^2}+|u|_{C_x^2}^{(1)}\right)\right];
$$
a further calculation implies
$$
\begin{aligned}
\left|\left.D_y^kB(y_1X,y_2u)[\phi]\right|_{y=0}\right|
&\leq k!c_2^k\left(|X|_{C^2}+|u|_{C_x^2}^{(1)}\right)^k\|\phi\|^2_{H^1_x(g_0)}\\
&\leq k!c_2^k\left(|X|_{C^2}+|u|_{C_x^2}^{(1)}\right)^k\left(B_0[\phi]+\|\phi\|^2_{L^2_x(g_0)}\right).
\end{aligned}
$$
Applying theorem VII.4.9 of \cite{Kato2013}, we find that the resolvent $\mathcal{R}(\zeta,L(\mathcal{E}_{y_1X},y_2u))$ exists and is analytic in $y$ as long as
$$
|y|
<\frac{1}{c_2\left(|X|_{C^2}+|u|_{C_x^2}^{(1)}\right)}
\frac{1}{1+\|(1-L_{\text{Id}})\mathcal{R}(\zeta,L_{\text{Id}})\|}.
$$
Setting $y=1$, applying theorem VI.3.9 of \cite{Kato2013}, we conclude the following: there is a constant $c_3$ such that if
\begin{equation}\label{small}
|X|_{C^2}+|u|_{C_x^2}^{(1)}<c_3,
\end{equation}
then for $|\zeta|=3$, the resolvent $\mathcal{R}(\zeta-L(\varphi,u))$ is a convergent power series of $[\partial^j_t(\nabla^{g_0})^lu]_{0\leq j\leq1}^{0\leq l\leq2}$ and $[(\nabla^{g_0})^lX]_{0\leq l\leq2}$, and furthermore
$$
\|\mathcal{R}(\zeta-L(\varphi,u))-\mathcal{R}(\zeta,L_\text{Id})\|
\leq4\left(|X|_{C^2}+|u|_{C_x^2}^{(1)}\right).
$$
If (\ref{small}) holds, then the spectral projection
$$
\mathcal{P}_0(\varphi,u)
:=\frac{1}{2\pi i}\int_{|\zeta|=3}\mathcal{R}(\zeta-L(\varphi,u))d\zeta
$$
is a convergent power series of $[\partial^j_t(\nabla^{g_0})^lu]_{0\leq j\leq1}^{0\leq l\leq2}$ and $[(\nabla^{g_0})^lX]_{0\leq l\leq2}$, hence by lemma I.4.10 of \cite{Kato2013}, we know that $\text{dim}[\mathrm{Ran}P(\varphi,u)]\equiv 3$, and
$$
\|\mathcal{P}_0(\varphi,u)-\mathcal{Q}_0^\text{Id}\|\leq8\pi\left(|X|_{C^2}+|u|_{C_x^2}^{(1)}\right).
$$

(C) We notice that $L^1(\varphi,u)$ is a perturbation of the operator $-\Delta_{g_0}^1+2$, acting on sections of $T^*(S^2)$ (the cotangent bundle). Now if $\lambda$ is an eigenvalue of $\Delta_{g_0}^1$ and the 1-form $\theta$ is an eigenvector, then $\Delta_{g_0}^1\theta=\lambda\theta$, so applying the co-differential operator $\delta_{g_0}$ we find
$$
-\Delta_{g_0}(\delta_{g_0}\theta)=\lambda\delta_{g_0}\theta.
$$
Thus $\lambda$ must be (the negative of) an eigenvalue of the scalar Laplacian and $\delta_{g_0}\theta$ must be a scalar eigenfunction, i.e. a sphere harmonic. Since the first cohomology group of $S^2$ vanishes, there is no non-trivial harmonic 1-form on $S^2$, so $\theta$ is completely determined by $\Delta_{g_0}^1\theta=d\delta_{g_0}\theta$. Consequently, a 1-form is an eigenvector of $\Delta_{g_0}^1$ if and only if it is the differential of a sphere harmonic. Thus the spectrum of $-\Delta_{g_0}^1+2$ is $\{0,-4,...\}$, and the eigenspace corresponding to 0 is three dimensional, which is spanned by the exterior differentiation of first three sphere harmonics.

Perturbation theory applied to $L^1(\varphi,u)$ under the inner product
$$
\langle\theta_1,\theta_2\rangle_{L^2_x(g(w))}
:=\int_{S^2}\langle\theta_1,\theta_2\rangle_{g(w)}d\mu(w)
$$
then gives the desired result, just as in (B). We only need to verify that the perturbation really is controlled in terms of $|X|_{C^2}$ and $|u|_{C_x^2}^{(1)}$. In fact, still with $w_y=i_0\circ\mathcal{E}_{y_1X}+y_2u$, the sesquilinear form for a section $\theta$ of $T^*(S^2)$ to be investigated is
$$
\begin{aligned}
B^1(y_1X,y_2u)[\theta]
&:=\int_{S^2}\left(|d\theta|_{g(w_y)}^2+|\delta_{g(w_y)}\theta|_{g(w_y)}^2-|h(w_y)|^2|\theta|_{g(w_y)}^2\right)d\mu(w_y)\\
&\quad+\int_{S^2}\left(H(w_y)-\frac{\kappa}{\mathrm{Vol}(w_y)}\right)H(w_y)|\theta|_{g(w_y)}^2d\mu(w_y)+\int_{S^2}|\partial_tN(w_y)|^2|\theta|_{g(w_y)}^2 d\mu_0.
\end{aligned}
$$
The only term that needs additional investigation is $|\delta_{g(w)}\theta|_{g(w)}^2$. By definition,
$$
\delta_{g(w_y)}\theta=-\star_{g(w_y)}^{-1}d\star_{g(w_y)}\theta,
$$
where $\star_{g(w_y)}$ is the Hodge star operator associated to $g(w_y)$, whose local coordinate representation is
$$
\star_{g(w_y)}\theta=-\sqrt{\det g(w_y)}\theta_idx^i.
$$
So in local coordinates, the coefficients of $\delta_{g(w_y)}\theta$ depend on first order spatial derivatives of $g(w_y)$, hence second order spatial derivatives of $X$ and $u$.
\end{proof}

As a corollary, if $|X|_{C^2}+|u|_{C_x^2}^{(1)}<\delta_0$, then
$$
\langle-L(\varphi,u)\phi,\phi\rangle_{L^2_x(g_0)}
$$
is equivalent to the usual $H^1$-norm on the subspace $\mathrm{Ran}[\mathcal{P}_\infty(\text{id}),0]$. Similar result holds for
$$
\langle-L^1(\varphi,u)\theta,\theta\rangle_{L^2_x(g_0)}.
$$

\subsection{Tame Elliptic Estimate for $L_\varphi$}
Let us still follow the notation of last subsection; in particular, we set $P_\lambda^\varphi$ to be the eigenprojection of $L_\varphi$ corresponding to eigenvalue $\lambda$. The coefficients of $L(\varphi,u)$ depends smoothly on derivatives of $\varphi$ and $u$ up to order 2. We know from Proposition \ref{SpectralProp} that the only zero eigenmodes of $L_\varphi$ are components of $N(i_\varphi)$, and
$$
\mathcal{P}_0^\varphi\phi
=\sum_{k=1}^3\left(\frac{1}{\|N^k(i_\varphi)\|_{L^2(g_0)}}\int_{S^2}N^k(i_\varphi)\phi d\mu_0\right)N^k(i_\varphi).
$$

Let us prove some estimate on $L_\varphi$. The idea is to consider it as the perturbation of $L_{\text{Id}}$. Note that
$$
\sum_{\lambda\in\sigma[-L_\varphi]\setminus\{0\}}\lambda\|\mathcal{P}_\lambda^\varphi\phi\|_{L^2(g_0)}^2
=-\langle L_\varphi(1-\mathcal{P}_0^\varphi)\phi,(1-\mathcal{P}_0^\varphi)\phi\rangle_{L^2(g_0)},
$$
while
$$
\begin{aligned}
|\langle L_\varphi(1-\mathcal{P}_0^\varphi)\phi,(1-\mathcal{P}_0^\varphi)\phi\rangle_{L^2(g_0)}&-\langle L_{\text{Id}}(1-\mathcal{P}_0^\varphi)\phi,(1-\mathcal{P}_0^\varphi)\phi\rangle_{L^2(g_0)}|\\
&\leq C|X|_{C^3}\|(1-\mathcal{P}_0^\varphi)\phi\|_{H^1(g_0)}^2.
\end{aligned}
$$
But on the other hand, if $|X|_{C^3}$ is sufficiently small, then $1-\mathcal{P}_0^\text{Id}$ is an isomorphism from $\mathrm{Ran}(1-\mathcal{P}_0^\text{Id})$ to $\mathrm{Ran}(1-\mathcal{P}_0^\varphi)$, also by the standard perturbation theory as presented in \cite{Kato2013}. Thus there is a $C^2$-neighbourhood $\mathfrak{U}_0\subset\mathfrak{X}$ of 0 such that if $X\in\mathfrak{U}_0$, then
$$
\begin{aligned}
-\langle L_{\text{Id}}(1-\mathcal{P}_0^\varphi)&\phi,(1-\mathcal{P}_0^\varphi)\phi\rangle_{L^2(g_0)}\\
&=-\langle L_{\text{Id}}(1-\mathcal{P}_0^{\varphi})(1-\mathcal{P}_0^{\varphi})\phi,(1-\mathcal{P}_0^{\varphi})(1-\mathcal{P}_0^{\varphi})\phi\rangle_{L^2(g_0)}\\
&\geq\|(1-\mathcal{P}_0^{\varphi})\phi\|_{H^1(g_0)}^2.
\end{aligned}
$$
Thus by elliptic regularity theory, we find there is a universal constant $C>0$ such that
$$
\begin{aligned}
(1-C|X|_{C^3})\|(1-\mathcal{P}_0^{\varphi})\phi\|_{H^1(g_0)}^2
&\leq\sum_{\lambda\in\sigma[-L_\varphi]\setminus\{0\}}\lambda\|\mathcal{P}_\lambda^{\varphi}\phi\|_{L^2(g_0)}^2\\
&\leq(1+C|X|_{C^3})\|(1-\mathcal{P}_0^{\varphi})\phi\|_{H^1(g_0)}^2.
\end{aligned}
$$

In order to estimate higher Sobolev norms, we consider the function $[-L_\varphi]^{n/2}\phi$, $n\geq1$, which is defined via spectral calculus:
$$
[-L_\varphi]^{n/2}\phi:=\sum_{\lambda\in\sigma[-L_\varphi]}\lambda^{n/2}\mathcal{P}_\lambda^{\varphi}\phi.
$$
We are at the place to estimate $\|\phi\|_{H_x^n(g_0)}$ in terms of $\|[-L_\varphi]^{n/2}\phi\|_{L^2(g_0)}$. The following proposition establishes the fundamental tame elliptic estimate that will be used in the next section.

\begin{proposition}\label{TameEllipticEstimates1}
There is a $\delta_1>0$ such that if $|X|_{C^2}<\delta_1$, then with $\varphi=\mathcal{E}_X$,

(A) For any integer $n$,
$$
\|\phi\|_{H^{n}(g_0)}\leq C_n\left(\|[-L_\varphi]^{n/2}\phi\|_{L^2(g_0)}+(1+|X|_{C^{n+2}})\|\phi\|_{L^2(g_0)}\right),
$$
and
$$
\|[-L_\varphi]^{n/2}\phi\|_{L^2(g_0)}\leq C_n\left(\|\phi\|_{H^{n}(g_0)}+|X|_{C^{n+2}}\|\phi\|_{L^2(g_0)}\right).
$$

(B) For any integer $n$,
$$
\|\phi\|_{H^{n+1}(g_0)}\leq C_n\left(\|[-L_\varphi]^{n/2}\phi\|_{H^1(g_0)}+(1+|X|_{C^{n+3}})\|\phi\|_{L^2(g_0)}\right),
$$
and
$$
\|[-L_\varphi]^{n/2}\phi\|_{H^1(g_0)}
\leq C_n\left(\|\phi\|_{H^{n+1}(g_0)}+|X|_{C^{n+3}}\|\phi\|_{L^2(g_0)}\right),
$$
\end{proposition}
\begin{proof}
(A) First we consider $n=2m$. Take the covering of $S^2$ by two disks $B_1,B_2$ as in the beginning of subsection \ref{Subsec2.1}, and let $\{\zeta_1,\zeta_2\}$ be a corresponding smooth partition of unity. Each disk is then a coordinate patch of $S^2$. Thus it suffices to establish all the estimates on $B_1$, with $\phi_1=\zeta_1\phi$. Under this local coordinate, the operator $L_\varphi$ has a representation
$$
\begin{aligned}
L_\varphi\phi_1
&=\left(\frac{d\mu_\varphi}{d\mu_0}\Delta_{g(i_\varphi)}+2\frac{d\mu_\varphi}{d\mu_0}\right)\phi_1-\frac{6}{4\pi}\frac{d\mu_\varphi}{d\mu_0}\int_{S^2}\phi_1d\mu_\varphi\\
&=a_X^{ij}(x)\frac{\partial^2\phi_1(x)}{\partial x^i\partial x^j}+p_X^k(x)\frac{\partial\phi_1(x)}{\partial x^k}+q_X(x)\phi_1(x)+r_X(x)\int_{B_1}\phi_1(y)\varrho_X(y)dy.
\end{aligned}
$$
Here $\varrho_X(x)dx$ is the local coordinate representation of the measure $d\mu_\varphi$, and the matrix
$$
[a_X^{ij}(x)]=\left[\frac{d\mu_\varphi}{d\mu_0}(x)g^{ij}(i_\varphi)(x)\right]
$$
is uniformly positive definite if $|X|_{C^2}$ is small. Furthermore, by Lemma \ref{GeometricQuantities}, the functions $a_X^{ij}(x)$, $\mathcal{P}_X^k(x)$, $q_X(x)$, $r_X(x)$, $\varrho_X(x)$ depend smoothly on $D_x^lX,0\leq l\leq2$. We write
$$
[L_\varphi]^m
=\sum_{|\alpha|\leq2m}A_\alpha D_x^\alpha.
$$
The principal symbol of $[L_\varphi]^m$ reads
$$
\sum_{|\alpha|=2m}A_\alpha\xi^\alpha=\left(\sum_{i,j=1}^2a_X^{ij}\xi^i\xi^j\right)^m=\sum_{1\leq i_l,j_l\leq 2}a_X^{i_1j_1}\cdots a_X^{i_mj_m}\xi^{i_1}\xi^{j_1}\cdots \xi^{i_m}\xi^{j_m}.
$$

We shall estimate
$$
\|[-L_\varphi]^m\phi_1\|^2_{L^2_x(g_0)}=\int_{S^2}L_\varphi^m\phi_1\cdot L_\varphi^m\phi_1d\mu_0.
$$

For the principal symbol, it suffices to apply the standard G?rding inequality. We quote the following simplified version of Theorem 6.5.1. from Morrey's monograph \cite{Morrey1966}:

\textbf{G?rding's inequality}
\emph{Let $B_1$ be the unit ball in $\mathbb{R}^l$, and suppose $\{a_{\alpha\beta}\}_{|\alpha|,|\beta|=m}$ are continuous real functions defined on $B_1$ such that}

\emph{(1) $|a_{\alpha\beta}|\leq M\,\forall\alpha,\beta$;}

\emph{(2) $\sum_{|\alpha|,|\beta|=m}a_{\alpha\beta}\xi^\alpha\xi^\beta\geq\mu|\xi|^{2m},\,\forall\xi$.}

\emph{Then for any $f\in H^m_0(B_1)$, there holds
$$
\int_{B_1}\sum_{|\alpha|,|\beta|=m}
a_{\alpha\beta}D^\alpha f D^\beta f
\geq
\frac{\mu}{2}\|f\|_{H^m_0}^2-C\|f\|_{L^2}^2,
$$
where the constant $C$ depends on $l$, $m$, $M$, $\mu$ and the modulus of continuity of the coefficients $a_{\alpha\beta}$.}

When $|X|_{C^2}$ is sufficiently small, the bounds of the coefficients of the principal symbol are controlled in terms of $\inf|{d\mu_\varphi}/{d\mu_0}|$, hence in terms of $|X|_{C^1}$. The modulus of continuity of the coefficients of the principal symbol is controlled in terms of $m$ and $\text{Lip}[g(i_\varphi)d\mu_\varphi/d\mu_0]$, hence in terms of $m$ and $|X|_{C^2}$. Thus, if $|X|_{C^2}$ is so small that $\inf|d\mu_\varphi/d\mu_0|>1/2$, $\sup|d\mu_\varphi/d\mu_0|<2$,  $\text{Lip}[g(i_\varphi)d\mu_\varphi/d\mu_0]<10$, we have, by G?rding's inequality,
$$
\int_{S^2}\left|\sum_{|\alpha|=2m}A_\alpha D^\alpha_x\phi_1\right|^2d\mu_0
\geq 4^{-m}\|\phi_1\|_{H^{2m}_x(g_0)}^2-C_m\|\phi_1\|_{L^2_x(g_0)}^2,
$$
where $C_m$ is a constant depending on $m$. An elementary argument also gives
$$
\int_{S^2}\left|\sum_{|\alpha|=2m}A_\alpha D^\alpha_x\phi_1\right|^2d\mu_0
\leq C_m\left(\|\phi_1\|_{H^{2m}_x(g_0)}^2+\|\phi_1\|_{L^2_x(g_0)}^2\right).
$$

For lower order terms, we simply apply Lemma \ref{OperatorComposition}, since each $A_\alpha D^\alpha_x$ is the linear combination of composition of the operators
$$
a_{X}^{ij}\partial_i\partial_j\phi,\,
p^k_X\partial_k\phi,\,
q_X\phi,\,
r_X\int_{B_1}\phi\varrho_X.
$$
The integral operator does not affect the validity of Lemma \ref{OperatorComposition}. Hence for each multi-index $\alpha$ we have
$$
\|A_\alpha D^\alpha\phi_1\|_{L^2_x(g_0)}
\leq C_\alpha\left(
\|\phi_1\|_{H^{|\alpha|}_x(g_0)}+|X|_{C^{|\alpha|+2}}\|\phi_1\|_{L^2_x(g_0)}
\right).
$$
Thus
$$
\|[-L_\varphi]^m\phi_1\|^2_{L^2_x(g_0)}
\leq C_n\left(\|\phi_1\|_{H^{n}(g_0)}+|X|_{C^{n+2}}\|\phi_1\|_{L^2(g_0)}\right).
$$
Furthermore, we have
$$
\begin{aligned}
\|[-L_\varphi]^m\phi_1\|^2_{L^2_x(g_0)}
&\geq 2^{-m}\|\phi_1\|_{H^{2m}_x(g_0)}^2
-C_m\|\phi_1\|_{L^2_x(g_0)}^2\\
&\quad-\sum_{|\beta|\leq2m-1}\sum_{|\alpha|\leq 2m}\|A_\alpha D_x^\alpha\phi_1\|_{L^2_x(g_0)}
\|A_\beta D_x^\beta\phi_1\|_{L^2_x(g_0)},
\end{aligned}
$$
and we may apply Young's inequality $ab\leq\varepsilon a^2+b^2/\varepsilon$ to estimate the sum; taking a sufficiently small $\varepsilon>0$ depending on $n$ along, gluing back, we obtain the following estimate: there is a $C^2$-neighbourhood $\mathfrak{U}_0\subset\mathfrak{X}$ of $0$ such that if $X\in\mathfrak{U}_0$, then for any even $n$,
$$
\|\phi\|_{H^{n}(g_0)}\leq C_n\left(\|[-L_\varphi]^{n/2}\phi\|_{L^2(g_0)}+(1+|X|_{C^{n+2}})\|\phi\|_{L^2(g_0)}\right),
$$
and
$$
\|[-L_\varphi]^{n/2}\phi\|_{L^2(g_0)}\leq C_n\left(\|\phi\|_{H^{n}(g_0)}+|X|_{C^{n+2}}\|\phi\|_{L^2(g_0)}\right).
$$

Next we deal with the case $n=2m+1$. Note that
$$
\begin{aligned}
\|[-L_\varphi]^{n/2}\phi\|_{L^2(g_0)}^2
&=-\langle L_\varphi[-L_\varphi]^{m}\phi,[-L_\varphi]^{m}\phi\rangle_{L^2(g_0)}\\
&\sim_{|X|_{C^2}}\|[-L_\varphi]^{m}\phi\|_{H^1(g_0)}^2\\
&=\|\nabla^{g_0}[-L_\varphi]^{m}\phi\|_{L^2(g_0)}^2+\|[-L_\varphi]^{m}\phi\|_{L^2(g_0)}^2.
\end{aligned}
$$
By a similar interpolation argument as above, combining the results, we finally obtain (A).

(B) To estimate $\|\phi\|_{H_x^{n+1}(g_0)}$ in terms of $\|[-L_\varphi]^{n/2}\phi\|_{H^1(g_0)}$, we notice that when $n$ is even this is done similarly as (A), and when $n=2m+1$ is odd, we have
$$
\begin{aligned}
\|[-L_\varphi]^{n/2}\phi\|_{H^1(g_0)}^2
&\sim_{|X|_{C^2}}-\langle L_\varphi[-L_\varphi]^{n/2}\phi,[-L_\varphi]^{n/2}\phi\rangle_{L^2(g_0)}\\
&=-\langle [-L_\varphi]^{m+1}\phi,[-L_\varphi]^{m+1}\phi\rangle_{L^2(g_0)}.
\end{aligned}
$$
Hence (B) follows.
\end{proof}

Finally, as in proposition \ref{SpectralProp}, we still define the perturbed eigenprojection
$$
\mathcal{P}_0(\varphi,u)
:=\frac{1}{2\pi i}\int_{|\zeta|=3}\mathcal{R}(\zeta-L(\varphi,u))d\zeta
$$
and $\mathcal{P}_\infty(\varphi,u)=\text{Id}-\mathcal{P}_0(\varphi,u)$. Adapting proposition \ref{SpectralProp} and applying G{\aa}rding's inequality to $L(\varphi,u)$ similarly, we also obtain
\begin{proposition}\label{TameEllipticEstimates2}
There is a constant $\delta_1>0$ such that if $|X|_{C^2}+|u|_{C^2_x}^{(1)}<\delta_1$, then
$$
\|\phi\|_{H^1_x(g_0)}^2
\leq C\left(\langle-L(\varphi,u)\mathcal{P}_\infty(\varphi,u)\phi,\mathcal{P}_\infty(\varphi,u)\phi\rangle_{L^2(g_0)}
+\|\mathcal{P}_0(\varphi,u)\phi\|_{L^2(g_0)}^2
\right).
$$
\end{proposition}

\subsection{Decay of the Linearized Equation with $u=0$}
Let us now investigate the Cauchy problem for the scalar equation
\begin{equation}\label{LEQN}
\frac{\partial^2\phi}{\partial t^2}+b\frac{\partial\phi}{\partial t}-L_\varphi\phi=\gamma(t).
\end{equation}
Here we fix a $b>0$, a $\gamma\in C^2([0,\infty); C^\infty(S^2;\mathbb{R}^3))$, and a $\varphi\in\mathfrak{Diff}(S^2)$ close to the identity and represent it as $\varphi=\mathcal{E}_X$ as in (\ref{EX}). The Cauchy problem of (\ref{LEQN}) is explicitly solved as
\begin{equation}\label{Phi}
\begin{aligned}
\phi(t)&=\mathcal{P}_0^\varphi\phi(0)+\frac{1-e^{-bt}}{b}\mathcal{P}_0^\varphi\phi'(0)+\int_0^t\frac{1-e^{-b(t-s)}}{b}\mathcal{P}_0^\varphi\gamma(s)ds\\
&\quad+\sum_{\lambda\in\sigma[-L_\varphi]\setminus\{0\}}
\frac{\omega_b^+(\lambda)e^{\omega_b^-(\lambda)t}-\omega_b^-(\lambda)e^{\omega_b^+(\lambda)t}}{\omega_b^+(\lambda)-\omega_b^-(\lambda)}\mathcal{P}_\lambda^\varphi\phi(0)
+\frac{e^{\omega_b^+(\lambda)t}-e^{\omega_b^-(\lambda)t}}{\omega_b^+(\lambda)-\omega_b^-(\lambda)}\mathcal{P}_\lambda^\varphi\phi'(0)\\
&\quad+\sum_{\lambda\in\sigma[-L_\varphi]\setminus\{0\}}
\int_0^t\frac{e^{\omega_b^+(\lambda)(t-s)}-e^{\omega_b^-(\lambda)(t-s)}}{\omega_b^+(\lambda)-\omega_b^-(\lambda)}\mathcal{P}_\lambda^\varphi\gamma(s)ds,
\end{aligned}
\end{equation}
where
$$
\omega_b^\pm(\lambda)=\frac{-b\pm\sqrt{b^2-4\lambda}}{2}.
$$
Note that by Proposition \ref{SpectralProp}, if $|X|_{C^2}$ is sufficiently small, then all the non-zero eigenvalues of $-L_\varphi$ will be strictly greater than $2$, so with
\begin{equation}\label{OmegaB}
\beta:=\left\{
\begin{matrix}
\displaystyle{{\frac{b-\sqrt{b^2-4}}{3}}}, & b\geq2\\
b/3, & b<2
\end{matrix}
\right.
\end{equation}
we have
$$
\mathrm{Re}\left(\omega_b^\pm(\lambda)\right)
<-\beta<0,\quad
\forall\lambda\in\sigma[-L_\varphi]\setminus\{0\},
$$
and $-b<-\beta<0$.

As in section 2, we define a norm that captures the exponential decay for scalar functions as follows:
$$
\|\gamma\|_{\beta,0;n}
:=\sup_{t\geq0}e^{\beta t}\|\gamma(t)\|_{H^n_x(g_0)}.
$$
Recall that $\mathbf{E}_{\beta,0}^n$ is the space of all scalar functions $\gamma\in C^0([0,\infty);H^n(S^2;\mathbb{R}))$ such that $\|\gamma\|_{\beta,0;n}<\infty$ for any $n$. With the condition $\gamma\in\mathbf{E}_{\beta,0}^n$, the Duhamel integrals in (\ref{Phi}) can be estimated as follows. For example,
$$
\begin{aligned}
\left\|\int_0^te^{-b(t-s)}\mathcal{P}_0^\varphi\gamma(s)ds\right\|_{L^2_x(g_0)}
&\leq e^{-bt}\int_0^{t}e^{(b-\beta)s}\|\mathcal{P}_0^\varphi\gamma(s)\|_{\beta;0}ds\\
&\leq\frac{e^{-\beta t}}{b-\beta}\|\mathcal{P}_0^\varphi\gamma\|_{\beta;0}.
\end{aligned}
$$
\emph{Thus the solution operator to (\ref{LEQN}) does not give rise to loss of decay.}

On the other hand, when $|X|_{C^2}$ is small, Proposition \ref{SpectralProp} shows that there is a constant $C$ independent of $X$ and $\lambda\in\sigma[-L_\varphi]$ such that
$$
\left|\frac{\sqrt\lambda}{\omega_b^+(\lambda)-\omega_b^-(\lambda)}\right|\leq C,
\,\left|\frac{\omega_b^\pm(\lambda)}{\omega_b^+(\lambda)-\omega_b^-(\lambda)}\right|\leq C,
$$
and as $\lambda\to\infty$, there holds
$$
\left|\frac{\sqrt\lambda}{\omega_b^+(\lambda)-\omega_b^-(\lambda)}\right|\to\frac{1}{2},
\,\left|\frac{\omega_b^\pm(\lambda)}{\omega_b^+(\lambda)-\omega_b^-(\lambda)}\right|\to\frac{1}{2}.
$$
As a result of these, the terms in (\ref{LEQN}) involving $\phi'(0)$ and $\gamma$ gains one spatial derivative.

Hence we obtain the following proposition:
\begin{proposition}\label{DecayH1}
Fix $b>0$, and define $\beta$ as in (\ref{OmegaB}). There is a $C^2$-neighbourhood $\mathfrak{U}_0\subset\mathfrak{X}$ of $0$ and a constant $C=C(\mathfrak{U}_0,b)$ with the following properties. For any $X\in\mathfrak{U}_0$ and any $\gamma\in \mathbf{E}_{\beta,1}$, the solution $\phi(t)$ to the Cauchy problem for (\ref{LEQN}) has a limit $\phi(\infty)$ in $H^1_x(g_0)$, which is equal to
$$\mathcal{P}_0^\varphi\phi(0)+b^{-1}\mathcal{P}_0^\varphi\phi'(0)+b^{-1}\int_0^\infty \mathcal{P}_0^\varphi\gamma(s)ds,$$
and furthermore,
$$
\begin{aligned}
\|\phi(t)&-\phi(\infty)\|_{\beta,1}+\|\phi'(t)\|_{\beta,0}\\
&\leq C\left(\|\phi(0)\|_{H^1(g_0)}+\|\phi'(0)\|_{L^2(g_0)}
+\|\gamma(t)\|_{\beta,0}\right).
\end{aligned}
$$
\end{proposition}

To estimate Sobolev norm $\|\phi\|_{H^n_x(g_0)}$, we just have to consider $[-L_\varphi]^{n/2}\phi$ apply the above results to the Cauchy problem
$$
\frac{\partial^2}{\partial t^2}[-L_\varphi]^{n/2}\phi+b\frac{\partial}{\partial t}[-L_\varphi]^{n/2}\phi-L_\varphi[-L_\varphi]^{n/2}\phi=[-L_\varphi]^{n/2}\gamma(t).
$$
With the aid of Proposition \ref{TameEllipticEstimates1} (tame elliptic estimates) and Proposition \ref{DecayH1} (decay estimate), we obtain the following proposition:
\begin{proposition}\label{DecayHn}
Fix $b>0$, and define $\beta$ as in (\ref{OmegaB}). There is a $\delta_2>0$ such that if $|X|_{C^2}<\delta_2$, then for any $\gamma\in \mathbf{E}_{\beta,1}$, the solution $\phi(t)$ to the Cauchy problem for (\ref{LEQN}) has a limit $\phi(\infty)$ in $\mathrm{Ran}[\mathcal{P}_0^\varphi]$, which is equal to
$$
\mathcal{P}_0^\varphi\phi(0)+b^{-1}\mathcal{P}_0^\varphi\phi'(0)+b^{-1}\int_0^\infty \mathcal{P}_0^\varphi\gamma(s)ds,
$$
and furthermore,
$$
\begin{aligned}
\|\phi(t)&-\phi(\infty)\|_{\beta;n+1}+\|\phi'(t)\|_{\beta;n}\\
&\leq C_n\left(\|\phi(0)\|_{H^{n+1}(g_0)}+\|\phi'(0)\|_{H^n(g_0)}
+\|\gamma(t)\|_{\beta,n}\right)\\
&\quad+C_n|X|_{C^{n+3}}\left(\|\phi(0)\|_{H^{1}(g_0)}+\|\phi'(0)\|_{L^2(g_0)}
+\|\gamma(t)\|_{\beta;0}\right).
\end{aligned}
$$
The constants $C_n$ depend only on $n,b,\beta$.
\end{proposition}

We now have all the ingredients to solve $(Y,c,v)$ out of the linearized equation $\Phi'(X,a,0)(Y,c,v)=f$. The symbols in the statement of the following proposition are defined in equations (\ref{EX})-(\ref{SigmaPhi}).
\begin{proposition}\label{Linearized(u=0)}
Fix $b>0$, and define $\beta$ as in (\ref{OmegaB}). Consider the Cauchy problem of the lienarized equation
$$
\Phi'(X,a,0)(Y,c,v)=f,\,f\in\mathbf{E}_{\beta,1},
$$
where with $\varphi=\mathcal{E}_X$,
$$
\eta(t)
:=\Sigma_\varphi Y+c+v(t),\quad
\phi(t) N(i_\varphi)+\psi(t)
:=[\eta(t)\cdot N(i_\varphi)] N(i_\varphi)+\top_{i_\varphi}\eta(t),
$$
the Cauchy data $\eta[0]=(\eta(0),\eta'(0))$ is given. Then there is a $\delta_3>0$ such that if $|X|_{C^2}<\delta_3$ and $\beta<\omega_b$, this Cauchy problem has a unique solution $(Y,c,v)\in\mathbf{F}_\beta$, satisfying the following tame estimates:
$$
|c|\leq C\left(\|\eta(0)\|_{L^2(g_0)}
+\|\eta'(0)\|_{L^2(g_0)}+\|f\|_{\beta,0;0}\right),
$$
$$
\begin{aligned}
\|Y\|_{H^n}
&\leq C_n\left(\|\eta(0)\|_{H^{n}(g_0)}
+\|\eta'(0)\|_{H^{n}(g_0)}+\|f\|_{\beta,0;n}\right)\\
&\quad+C_n|X|_{C^{n+1}}\left(\|\eta(0)\|_{H^{2}(g_0)}
+\|\eta'(0)\|_{H^{2}(g_0)}+\|f\|_{\beta,0;2}\right).
\end{aligned}
$$
$$
\begin{aligned}
\|v\|_{\beta,3;n}
&\leq C_n\left(\|\eta(0)\|_{H^{n+2}(g_0)}
+\|\eta'(0)\|_{H^{n+1}(g_0)}+\|f\|_{\beta,1;n+1}\right)\\
&\quad+C_n|X|_{C^{n+3}}\left(\|\eta(0)\|_{H^{2}(g_0)}
+\|\eta'(0)\|_{H^{2}(g_0)}+\|f\|_{\beta,1;2}\right).
\end{aligned}
$$
The constants $C_n$ depend only on $n,b,\beta$.
\end{proposition}
\begin{proof}
Rewrite the equation as the following weakly linear hyperbolic system as Notz did in \cite{Notz2010}:
\begin{equation}\label{LEQu=0}
\begin{aligned}
\frac{\partial^2\phi}{\partial t^2}+b\frac{\partial\phi}{\partial t}&=L_\varphi\phi+f_\perp(t),\\
\frac{\partial^2\psi}{\partial t^2}+b\frac{\partial\psi}{\partial t}&=f_\top(t),\\
\left(
\begin{matrix}
\phi(0)\\
\psi(0)
\end{matrix}
\right)=
\left(
\begin{matrix}
\eta(0)\cdot N(i_\varphi),\\
\top_{i_\varphi}\eta(0)
\end{matrix}
\right)&,\,
\left(
\begin{matrix}
\phi'(0)\\
\psi'(0)
\end{matrix}
\right)=\left(
\begin{matrix}
\eta'(0)\cdot N(i_\varphi),\\
\top_{i_\varphi}\eta'(0)
\end{matrix}
\right).
\end{aligned}
\end{equation}
It is easily seen that the Cauchy problem for this system has a unique global solution $(\phi(t),\psi(t))$ for $f=f_\perp N(i_\varphi)+f_\top\in \mathbf{E_\beta}$, whose components are given by
$$
\begin{aligned}
\phi(t)&=\mathcal{P}_0^\varphi\phi(0)+\frac{1-e^{-bt}}{b}\mathcal{P}_0^\varphi\phi'(0)+\int_0^t\frac{1-e^{-b(t-s)}}{b}\mathcal{P}_0^\varphi f_\perp(s)ds\\
&\quad+\sum_{\lambda\in\sigma[-L_\varphi]\setminus\{0\}}
\frac{\omega_b^+(\lambda)e^{\omega_b^-(\lambda)t}-\omega_b^-(\lambda)e^{\omega_b^+(\lambda)t}}{\omega_b^+(\lambda)-\omega_b^-(\lambda)}\mathcal{P}_\lambda^\varphi\phi(0)
+\frac{e^{\omega_b^+(\lambda)t}-e^{\omega_b^-(\lambda)t}}{\omega_b^+(\lambda)-\omega_b^-(\lambda)}\mathcal{P}_\lambda^\varphi\phi'(0)\\
&\quad+\sum_{\lambda\in\sigma[-L_\varphi]\setminus\{0\}}
\int_0^t\frac{e^{\omega_b^+(\lambda)(t-s)}-e^{\omega_b^-(\lambda)(t-s)}}{\omega_b^+(\lambda)-\omega_b^-(\lambda)}\mathcal{P}_\lambda^\varphi f_\perp(s)ds,\\
\psi(t)&=\psi(0)+\frac{1-e^{-bt}}{b}\psi'(0)+\int_0^t\frac{1-e^{-b(t-s)}}{b}f_\top(s)ds.
\end{aligned}
$$
To recover $(Y,c,v)\in\mathbf{F}_\beta$ from $\phi,\psi$, we notice that $\Sigma_\varphi  Y+c=\phi(\infty)N(i_\varphi)+\psi(\infty)$, so $(Y,c,v)$ is explicitly solved as
$$
c^k=\frac{1}{\|N^k(i_\varphi)\|_{L^2(g_0)}}\left[\int_{S^2}N^k(i_\varphi)\left([ \eta(0)+b^{-1}\eta'(0)]\cdot N(i_\varphi)+b^{-1}\int_0^\infty f_\perp(s)ds\right)d\mu_0\right],
$$
$$
Y=\Sigma_\varphi ^{-1}\left[-\top_{i_\varphi}c+\top_{i_\varphi}\left(\eta(0)+b^{-1}\eta'(0)\right)+b^{-1}\int_0^\infty f_\top(s)ds\right],
$$
$$
\begin{aligned}
v(t)&=-b^{-1}e^{-bt}\mathcal{P}_0^\varphi\phi'(0)N(i_\varphi)-b^{-1}\left(\int_0^te^{-b(t-s)}\mathcal{P}_0^\varphi f_\perp(s)ds\right)N(i_\varphi)\\
&\quad+N(i_\varphi)\sum_{\lambda\in\sigma[-L_\varphi]\setminus\{0\}}
\frac{\omega_b^+(\lambda)e^{\omega_b^-(\lambda)t}-\omega_b^-(\lambda)e^{\omega_b^+(\lambda)t}}{\omega_b^+(\lambda)-\omega_b^-(\lambda)}\mathcal{P}_\lambda^\varphi\phi(0)
+\frac{e^{\omega_b^+(\lambda)t}-e^{\omega_b^-(\lambda)t}}{\omega_b^+(\lambda)-\omega_b^-(\lambda)}\mathcal{P}_\lambda^\varphi\phi'(0)\\
&\quad+N(i_\varphi)\sum_{\lambda\in\sigma[-L_\varphi]\setminus\{0\}}
\int_0^t\frac{e^{\omega_b^+(\lambda)(t-s)}-e^{\omega_b^-(\lambda)(t-s)}}{\omega_b^+(\lambda)-\omega_b^-(\lambda)}\mathcal{P}_\varphi f_\perp(s)ds\\
&\quad-b^{-1}e^{-bt}\psi'(0)-b^{-1}\int_0^te^{-b(t-s)}f_\top(s)ds.
\end{aligned}
$$
By Proposition \ref{DecayHn}, Lemma \ref{GeometricQuantities} and the Sobolev embedding $H^2\hookrightarrow C^{1-\varepsilon}$, the estimate for $|c|$, $|Y|$ and $\|v\|_{\beta,1;n}$ follows. For higher derivatives in time, it suffices to take the equation (\ref{LEQu=0}) itself into account, and differentiate it once with respect to time; this gives the estimate for $\|v\|_{\beta,3;n}$ with a loss of two more spatial derivatives.
\end{proof}

\section{Tame Estimate II: the Full Linearized Problem}
Following the general Nash-Moser scheme, we will obtain tame estimates for the solutions of the linearized equations (\ref{LS1}) or (\ref{LS2}). Note that (\ref{LS2}) can be considered as special version of (\ref{LS1}), restricted to a smaller space. Thanks to Notz \cite{Notz2010}, we know that (\ref{LS1}) admits a unique smooth solution, since it reduces to a weakly hyperbolic linear system. The corresponding energy estimate is also established. We need to deduce more refined energy estimates and decay estimates for this weakly hyperbolic linear system. It will be clear from the proof that our refinement essentially reflects the stability of $S^2$.

We will study (\ref{LS2}) more delicately. Still set $\varphi=\mathcal{E}_X$, $w=i_\varphi+a+u$, and define, for a scalar function $\phi$,
$$
\begin{aligned}
&L(\varphi,u)\phi:=\\
&\frac{d\mu(w)}{d\mu_0}\left[\Delta_{g(w)}\phi+|h(w)|^2\phi+\left(-H(w)+\frac{\kappa}{\mathrm{Vol}(w)}\right)H(w)\phi-\frac{6}{4\pi}\int_{S^2}\phi d\mu(w)\right]-|\partial_tN(w)|^2\phi;
\end{aligned}
$$
for a differential 1-form $\theta$ on $S^2$,
$$
\begin{aligned}
L^1(\varphi,u)\theta:=
\frac{d\mu(w)}{d\mu_0}\left[\Delta_{g(w)}^1\theta+|h(w)|^2\theta+\left(-H(w)+\frac{\kappa}{\mathrm{Vol}(w)}\right)H(w)\theta\right]-|\partial_tN(w)|^2\theta.
\end{aligned}
$$
We still transfer the unknown from $(Y,c,v)$ to $(\phi,\psi)\in C^\infty(S^2)\oplus C^\infty([0,\infty)\times S^2;\mathbb{R}^3)$, via $$
\phi N(w)+\psi=(\perp_w\eta)\cdot N(w)+\top_w\eta,\,
\eta:=\Sigma_\varphi Y+c+v.
$$
The linearized equation $\Phi'(X,a,u)(Y,c,v)=f$ or $\Psi'(w)\eta=f$ is written as, in the local coordinate we fixed on $S^2$ at the beginning,
$$
\begin{aligned}
\frac{\partial^2\phi}{\partial t^2}&+b\frac{\partial\phi}{\partial t}-L(\varphi,u)\phi
+\frac{d\mu(w)}{d\mu_0}[\nabla^{g(w)} H(w)\cdot\psi]\\
&-\frac{d\mu(w)}{d\mu_0}\left(-H(w)+\frac{\kappa}{\mathrm{Vol}(w)}\right)\mathrm{div}^{g(w)}\psi
+2\frac{\partial\psi^k}{\partial t}[\partial_t\partial_kw\cdot N(w)]
+\psi^k[\partial_t^2\partial_kw\cdot N(w)]
=f_\perp,
\end{aligned}
$$
$$
\begin{aligned}
\frac{\partial^2\psi^k}{\partial t^2}&+b\frac{\partial\psi^k}{\partial t}+\frac{d\mu(w)}{d\mu_0}\left(-H(w)+\frac{\kappa}{\mathrm{Vol}(w)}\right)[(\nabla^{g(w)})^{;k}\phi-h_i^k(w)\psi^i]\\
&+b\psi^l[\partial_t\partial_kw\cdot\partial_lw]
+2\frac{\partial\phi}{\partial t}[\partial_tN(w)\cdot\partial_jw] g^{jk}(w)\\
&+\phi[\partial_t^2N(w)\cdot\partial_jw]g^{jk}(w)
+2\frac{\partial\psi^l}{\partial t}[\partial_t\partial_lw\cdot\partial_jw] g^{jk}(w)
+\psi^l[\partial_t^2\partial_lw\cdot\partial_jw] g^{jk}(w)
=f_\top^k,
\end{aligned}
$$
Note that we write $\psi=\psi^k\partial_kw$. Since time derivatives of tangent vector fields along $w$ are not necessarily tangent along $w$, the above equation in the components $\psi^k$ is derived to ensure the tangential property.

Using terser symbols, we may also write
\begin{equation}\label{LinearizationInComponents}
\begin{aligned}
\frac{\partial^2\phi}{\partial t^2}+b\frac{\partial\phi}{\partial t}-L(\varphi,u)\phi
+I_0(\varphi,u)\psi+I_1(\varphi,u)D\psi
&=f_\perp,\\
\frac{\partial^2\psi}{\partial t^2}+b\frac{\partial\psi}{\partial t}
+J_0(\varphi,u)(\psi,\partial_t\psi)+K_0(\varphi,u)\phi+K_1(\varphi,u)D\phi
&=f_\top+Q_0(\varphi,u)f_\perp.
\end{aligned}
\end{equation}
Here we set
$$
D\phi=(\partial_t\phi,d\phi),\,
D\psi=(\partial_t\psi,d\psi),
$$
where $d$ is the exterior differential with respect to $x$. Note that we consider $\psi$ as a $\mathbb{R}^3$-valued function, so $d\psi$ is a well-defined section of $T^*(S^2)\otimes\mathbb{R}^3$. The reason that we contract the operators to $I_{0,1},J_0,K_{0,1},Q_0$ is simple: under a fixed local coordinate, they are all scalar, vector or matrix-valued smooth functions in
$$
(x,\partial_t^jD_x^lu,D_x^mX);\,
0\leq j\leq2;\,0\leq l,m\leq3,
$$
and are of order $O\left([\partial_t^jD_x^lu,D_x^mX]_{j\leq2;l,m\leq3}\right)$ when $(X,u)\simeq0$, and in fact vanishes for $u=0$. By proposition 2.18 of \cite{Notz2010}, the Cauchy problem of (\ref{LinearizationInComponents}) admits a unique smooth solution $\eta\in C^\infty([0,\infty)\times S^2;\mathbb{R}^3)$ if the initial data and right-hand-side are smooth. We shall take this as granted, and refine the argument by obtaining a preciser energy estimate.

\subsection{Tame Energy Estimate}
A certain energy estimate holds for a solution of the linearized problem, and we shall state it in this subsection. It gives the rate of exponential growth of the Sobolev norms $\|\eta(t)\|_{H^n_x(g_0)}$. The derivation employs differential calculus on Riemannian manifolds, since localization to a coordinate patch will result in loss of information of growth. This estimate does not depend on any lower bound of $b$, so it applies to either the damped or damping-free equation.
\begin{proposition}\label{EnergyFull}
Fix a real number $T>0$. Let $b\geq0$. Suppose $\varphi=\mathcal{E}_X$ for some $X\in\mathfrak{X}$, and suppose $u\in C^3_tC^\infty_x([0,T]\times S^2;\mathbb{R}^3)$. Consider the Cauchy problem of (\ref{LinearizationInComponents}). There is a constant $\delta_4>0$ such that if $|X|_{C^3}+|u|_{C^3_tC^4_x}<\delta_4$, the solution $\eta=\phi N(w)+\psi$ satisfies the following energy estimate: if we set
$$
E_{n}[\eta]=\|\partial_t\eta\|_{H^{n}_x(g_0)}
+\|\eta\|_{H^{n}_x(g_0)},
$$
then for $t\leq T$,
$$
\begin{aligned}
E_{n}[\eta](t)
&\leq (1+t)^2Q_{u;n}(t)\left[E_{n+1}[\eta](0)+
\left(|X|_{C^{n+3}}+|u|_{C^3_tC^{n+3}_x}\right)E_2[\eta](0)\right]\\
&\quad+(1+t)^2Q_{u;n}(t)\left[\|f(0)\|_{H^n_x(g_0)}+\left(|X|_{C^{n+3}}+|u|_{C^3_tC^{n+3}_x}\right)\|f(0)\|_{H^{2}_x(g_0)}\right]\\
&\quad+(1+t)^2Q_{u;n}(t)\left[\int_0^t\|f(s)\|_{H^{n}_x(g_0)}^{(1)}ds
+\left(|X|_{C^{n+3}}+|u|_{C^3_tC^{n+3}_x}\right)
\int_0^t\|f(s)\|_{H^2_x(g_0)}^{(1)}ds\right],
\end{aligned}
$$
where for an increasing sequence of positive numbers $(M_n)$,
$$
Q_{u;n}(t)=M_n\exp\left(M_n
\int_0^t\sqrt{\sup_{\tau\geq s}|u(\tau)|_{C_x^4}^{(3)}}ds
\right).
$$
The constants $C_n,M_n$ do not depend on $T$.
\end{proposition}
\begin{proof}
All computations below, unless otherwise noted, will be done on a fixed time slice $t$, so dependence on time will be abbreviated.

We first introduce an auxiliary energy norm
\begin{equation}\label{EXun}
E_n^{X,u}[\eta]=\|\partial_t^2\phi\|_{H^{n}_x(g_0)}
+\|\partial_t\phi\|_{H^{n+1}_x(g_0)}
+\|\phi\|_{H^{n+1}_x(g_0)}
+\|\partial_t\psi\|_{H^{n+1}_x(g_0)}
+\|\psi\|_{H^{n+1}_x(g_0)}.
\end{equation}
We suppose that $|X|_{C^2}+|u|_{C_x^2}^{(1)}$ satisfies the requirement of Proposition \ref{SpectralProp} and \ref{TameEllipticEstimates2}. For simplicity we write $\lambda(t)=\sup_{t\geq s}|u(s)|_{C_x^4}^{(3)}$. Then $\lambda(t)$ is non-increasing. We assume further $\lambda(0)=|u|_{C_t^3C_x^4}<1$. The spectrum of $L(\varphi,u)$ is then a perturbation of that of $L_\varphi$. We also write
$$
\begin{aligned}
\gamma&=-I_0(\varphi,u)\psi-I_1(\varphi,u)D\psi
+f_\perp,\\
\xi&=-J_0(\varphi,u)(\psi,\partial_t\psi)-K_0(\varphi,u)\phi-K_1(\varphi,u)D\phi+f_\top+Q_0(\varphi,u)f_\perp.
\end{aligned}
$$
Using Lemma \ref{Interpolation}- \ref{GeometricQuantities}, we find
\begin{equation}\label{gammaHn}
\begin{aligned}
\|\gamma\|_{H^n_x(g_0)}
&\leq C_n\left[\|f\|_{H^n_x(g_0)}+
\left(|X|_{C^{n+3}}+|u|_{C_x^{n+3}}^{(2)}\right)\|f\|_{L^2_x(g_0)}
\right]\\
&\quad+C_n\lambda\left(\|\partial_t\psi\|_{H^{n}_x(g_0)}+\|\psi\|_{H^{n+1}_x(g_0)}\right)\\
&\quad+C_n\left(|X|_{C^{n+3}}+|u|_{C_x^{n+3}}^{(2)}\right)\left(\|\partial_t\psi\|_{L^2_x(g_0)}+\|\psi\|_{L^2_x(g_0)}\right),
\end{aligned}
\end{equation}

\begin{equation}\label{gammaHn'}
\begin{aligned}
\|\partial_t\gamma\|_{H^n_x(g_0)}
&\leq C_n\left[\|f\|_{H^n_x(g_0)}^{(1)}+
\left(|X|_{C^{n+3}}+|u|_{C_x^{n+3}}^{(3)}\right)\|f\|_{L^2_x(g_0)}^{(1)}
\right]\\
&\quad+C_n\lambda\left(\|\partial_t\psi\|_{H^{n}_x(g_0)}+\|\psi\|_{H^{n+1}_x(g_0)}\right)\\
&\quad+C_n\left(|X|_{C^{n+3}}+|u|_{C_x^{n+3}}^{(3)}\right)\left(\|\partial_t\psi\|_{L^2_x(g_0)}+\|\psi\|_{L^2_x(g_0)}\right),
\end{aligned}
\end{equation}

\begin{equation}\label{xiHn}
\begin{aligned}
\|\xi\|_{H^n_x(g_0)}
&\leq C_n\left[\|f\|_{H^n_x(g_0)}+
\left(|X|_{C^{n+3}}+|u|_{C_x^{n+3}}^{(2)}\right)\|f\|_{L^2_x(g_0)}
\right]\\
&\quad+C_n\lambda\left(\|\partial_t\phi\|_{H^n_x(g_0)}+\|\phi\|_{H^{n+1}_x(g_0)}
+\|\partial_t\psi\|_{H^n_x(g_0)}+\|\psi\|_{H^{n}_x(g_0)}\right)\\
&\quad+C_n\left(|X|_{C^{n+3}}+|u|_{C_x^{n+3}}^{(2)}\right)\left(\|\partial_t\phi\|_{L^2_x(g_0)}+\|\phi\|_{L^2_x(g_0)}+\|\partial_t\psi\|_{L^2_x(g_0)}+\|\psi\|_{L^2_x(g_0)}\right).
\end{aligned}
\end{equation}

\textbf{Step 1: Estimate of the velocities $\partial_t\phi$, $\partial_t\psi$.}

We write $\mathcal{P}_0^{\varphi,0}$ for the spectral projection of $d\mu_\varphi/d\mu_0(\Delta_{g_\varphi}+2)$ correspongding to eigenvalue 0, and $\mathcal{P}_\infty^{\varphi,0}=\text{id}-\mathcal{P}_0^{\varphi,0}$. Set $\phi=\phi_\infty+\phi_0$, with $\phi_\infty=\mathcal{P}_{\infty}^{\varphi,0}\phi$, $\phi_0=\mathcal{P}_0^{\varphi,0}\phi$. By a direct calculation, the evolution of $\phi_j$ is described by
\begin{equation}\label{EnergyNormal0}
\frac{\partial^2\phi_j}{\partial t^2}+b\frac{\partial\phi_j}{\partial t}
=L(\varphi,u)\phi_j
-[L(\varphi,u),\mathcal{P}^{\varphi,0}_j]\phi
+\mathcal{P}_j^{(0)}\gamma,
\end{equation}
where $j=0,\infty$. We find more explicitly
$$
\begin{aligned}
\left(L(\varphi,u)\mathcal{P}^{\varphi,0}_0\right) \phi
&=\sum_{k=1}^3\frac{1}{\|N^k(i_0)\|_{L^2(g_0)}}\left(\int_{S^2}N^k(i_0)\cdot\phi d\mu_0\right)L(\varphi,u)[N^k(i_0)],\\
\left(\mathcal{P}^{\varphi,0}_0(\varphi,u)\right) \phi
&=\sum_{k=1}^3\frac{1}{\|N^k(i_0)\|_{L^2(g_0)}}\left(\int_{S^2} L(\varphi,u)[N^k(i_0)]\cdot\phi d\mu_0\right)N^k(i_0).
\end{aligned}
$$
Thus both $[L(\varphi,u),\mathcal{P}^{\varphi,0}_j]$, $[\partial_tL(\varphi,u),\mathcal{P}^{\varphi,0}_j]$ are zeroth order operators, and since $L_{\varphi}$ commutes with the projections, an explicit calculation gives
$$
\|[L(\varphi,u),\mathcal{P}^{\varphi,0}_j]\phi\|_{L^2_x(g_0)}
\leq C\lambda\|\phi\|_{L^2_x(g_0)}.
$$

We next define several weighted energy norms:
$$
F^\perp_0[\phi]
=\left(\|\partial_t\phi\|_{L^2_x(g_0)}^2
-\langle\phi_\infty, L(\varphi,u)\phi_\infty \rangle_{L^2_x(g_0)}
+\lambda\|\phi_0\|_{L^2_x(g_0)}^2\right)^{1/2},
$$
$$
F^\top_0[\psi]
=\left(\|\partial_t\psi\|_{H^1_x(g_0)}^2
+\lambda\|\psi\|_{H^1_x(g_0)}^2\right)^{1/2},
$$
$$
F_0[\eta]=\left(F^\perp_0[\phi]^2+F^\perp_0[\partial_t\phi]^2+F^\top_0[\psi]^2\right)^{1/2}.
$$
Note that by our smallness assumption, we have
$$
\|\phi_\infty\|_{H^1_x(g_0)}^2
\leq-C\int_{S^2}\phi_\infty\cdot L(\varphi,u)\phi_\infty d\mu_0.
$$
We will then derive a Gr?nwall type inequality for $F_0[\eta]$.

We differentiate $F_0^\perp[\phi]$ first. Using (\ref{EnergyNormal0}), we obtain
$$
\begin{aligned}
\frac{1}{2}\frac{d}{dt}F_0[\phi]^2
&\leq\langle\partial_t^2\phi_\infty,\partial_t\phi_\infty\rangle_{L^2_x(g_0)}
+b\|\partial_t\phi_\infty\|_{L^2_x(g_0)}^2
-\frac{1}{2}\frac{d}{dt}\langle\phi_\infty, L(\varphi,u)\phi_\infty\rangle_{L^2_x(g_0)}\\
&\quad+\langle\partial_t^2\phi_0,\partial_t\phi_0\rangle_{L^2_x(g_0)}
+b\|\partial_t\phi_0\|_{L^2_x(g_0)}^2
+\lambda\langle\partial_t\phi_0,\phi_0\rangle_{L^2_x(g_0)}\\
&=\langle\phi_\infty,[\partial_tL(\varphi,u)]\phi\rangle_{L^2_x(g_0)}
-\langle\partial_t\phi_\infty,[L(\varphi,u),\mathcal{P}^{\varphi,0}_\infty]\phi\rangle_{L^2_x(g_0)}
+\langle\partial_t\phi,\gamma\rangle_{L^2_x(g_0)}\\
&\quad+\langle\partial_t\phi_0,L(\varphi,u)\phi_0 \rangle_{L^2_x(g_0)}^2
-\langle\partial_t\phi_0,[L(\varphi,u),\mathcal{P}^{\varphi,0}_0]\phi\rangle_{L^2_x(g_0)}^2
+\lambda\langle\partial_t\phi_0,\phi_0\rangle_{L^2_x(g_0)}^2\\
&\leq C\lambda\|\phi_\infty\|_{H^1_x(g_0)}^2
+C\lambda\|\partial_t\phi\|_{L^2_x(g_0)}\|\phi\|_{L^2_x(g_0)}
+\|\gamma\|_{L^2_x(g_0)}\|\partial_t\phi\|_{L^2_x(g_0)}.
\end{aligned}
$$
Note that in the first inequality we used $\lambda'\leq0$ so the term $\lambda'(t)\|\phi_0\|_{L^2_x(g_0)}^2$ is dropped. Strictly speaking, $\lambda'$ should be considered as a negative Borel measure, but this regularity problem does not affect what follows, since Gr?nwall's inequality remains valid for Lebesgue-Stieltjes integral. The first term of right-hand-side is controlled by $C\lambda F_0^\perp[\phi]^2$, and the second by
$$
C\sqrt{\lambda}\|\partial_t\phi\|_{L^2_x(g_0)}^2
+C\lambda^{3/2}\|\phi_\infty\|_{L^2_x(g_0)}^2
+C\lambda^{3/2}\|\phi_0\|_{L^2_x(g_0)}^2
\leq C\sqrt{\lambda}F_0^\top[\phi]^2,
$$
and the third is controlled by, using (\ref{gammaHn}) and imitating last inequality,
$$
\begin{aligned}
C\|f_\perp\|_{L^2_x(g_0)}\|\partial_t\phi\|_{L^2_x(g_0)}
&+C\lambda\left(\|\partial_t\psi\|_{L^2_x(g_0)}+\|\psi\|_{L^2_x(g_0)}\right)\|\partial_t\phi\|_{L^2_x(g_0)}\\
&\leq C\|f_\perp\|_{L^2_x(g_0)}F_0^\perp[\phi]
+C\sqrt{\lambda}F_0^\top[\psi]F_0^\perp[\phi].
\end{aligned}
$$
Thus we obtain a differential inequality: for some $M_0>0$,
\begin{equation}\label{EnergyNormalZeroTemp0}
\frac{d}{dt}F_0^\perp[\phi]^2
\leq M_0\sqrt{\lambda}F_0[\eta]^2+M_0\|f\|_{L^2_x(g_0)}F_0[\eta].
\end{equation}
We may apply this argument similarly to $\partial_t\psi$, thus obtaining
\begin{equation}\label{EnergyTangentZeroTemp0}
\frac{d}{dt}F_0^\top[\psi]^2
\leq M_0\sqrt{\lambda}F_0[\eta]^2+M_0\|f\|_{H^1_x(g_0)}^{(1)}F_0[\eta].
\end{equation}

\textbf{Step 2: Estimate of the acceleration $\partial_t^2\phi$ and the energy $E_0^{X,u}[\eta]$ in (\ref{EXun})}

Differentiating (\ref{EnergyNormal0}) with respect to $t$, we find the evolution equation for $\partial_t\phi_j$:
\begin{equation}\label{EnergyNormal0'}
\begin{aligned}
\frac{\partial^2}{\partial t^2}\frac{\partial\phi_j}{\partial t}+b\frac{\partial}{\partial t}\frac{\partial\phi_j}{\partial t}
&=L(\varphi,u)\frac{\partial\phi_j}{\partial t}
+[\partial_tL(\varphi,u)]\phi_j\\
&\quad-[{\partial_t}L(\varphi,u),\mathcal{P}^{\varphi,0}_j]\phi
-[L(\varphi,u),\mathcal{P}^{\varphi,0}_j]\frac{\partial\phi}{\partial t}
+\mathcal{P}_j^{(0)}\frac{\partial\gamma}{\partial t},
\end{aligned}
\end{equation}
where $j=0,\infty$. We may repeat the argument in step 1, with the only essential modification appearing when estimating the term $[\partial_tL(\varphi,u)]\phi_j$: it consists of second order spatial derivatives of $\phi$, and this regularity information is not guaranteed by the hyperbolic structure of the system. We thus employ the standard elliptic regularity theory: if $|X|_{C^2}+|u|_{C^1_tC^2_x}$ is small, then the modulus of continuity of the principal symbol of $L(\varphi,u)$ is controlled, so
$$
\|\phi\|_{H^2_x(g_0)}
\leq C\|L(\varphi,u)\phi\|_{L^2_x(g_0)}+C\|\phi\|_{L^2_x(g_0)}.
$$
Thus, keeping in mind $\partial_t^2\phi+b\partial_t\phi=L(\varphi,u)\phi+\gamma$, we estimate
$$
\begin{aligned}
\|[\partial_tL(\varphi,u)]\phi_j\|_{L^2_x(g_0)}
&\leq C\lambda\|\phi\|_{H^2_x(g_0)}\\
&\leq C\lambda\left(\|L(\varphi,u)\phi\|_{L^2_x(g_0)}+\|\phi\|_{L^2_x(g_0)}\right)\\
&\leq C\lambda\left(\|\partial_t^2\phi\|_{L^2_x(g_0)}+\|\partial_t\phi\|_{L^2_x(g_0)}+\|\phi\|_{L^2_x(g_0)}+\|\gamma\|_{L^2_x(g_0)}\right),
\end{aligned}
$$
so
$$
\|[\partial_tL(\varphi,u)]\phi_j\|_{L^2_x(g_0)}\|\partial_t\phi\|_{L^2_x(g_0)}
\leq C\sqrt{\lambda} F_0[\eta]^2.
$$
Imitating step 1, and adding (\ref{EnergyNormalZeroTemp0})(\ref{EnergyTangentZeroTemp0}), we finally obtain the differential inequality
\begin{equation}\label{EnergyNormalZero}
\frac{d}{dt}F_0[\eta]
\leq M_0\sqrt{\lambda}F_0[\eta]
+M_0\|f\|_{H^1_x(g_0)}^{(1)}.
\end{equation}
By Gr?nwall's inequality we find
\begin{equation}\label{EnergyNormalZeroTemp1}
F_0[\eta](t)
\leq \exp\left(M_0\int_0^t\sqrt{\lambda(s)}ds\right)
\left(F_0[\eta](0)+\int_0^t\|f(s)\|_{H^1_x(g_0)}^{(1)}ds\right).
\end{equation}
We can integrate (\ref{EnergyNormalZeroTemp1}) with respect to time and find further
\begin{equation}\label{EnergyNormalZeroTemp2}
\begin{aligned}
\|\phi(t)\|_{L^2_x(g_0)}+\|\psi(t)\|_{L^2_x(g_0)}
&\leq \|\phi(0)\|_{L^2_x(g_0)}+\|\psi(0)\|_{L^2_x(g_0)}\\
&\quad+C(1+t)Q_{u;0}(t)\left(F_0[\eta](0)+\int_0^t\|f(s)\|_{H^1_x(g_0)}^{(1)}ds\right).
\end{aligned}
\end{equation}
So adding (\ref{EnergyNormalZeroTemp1}) and (\ref{EnergyNormalZeroTemp2}), we finally obtain
\begin{equation}\label{Energyn=0}
\begin{aligned}
E_0^{X,u}[\eta](t)
&\leq CF_0[\phi](t)
+\|\phi(t)\|_{L^2_x(g_0)}
+\|\psi(t)\|_{L^2_x(g_0)}\\
&\leq M_0(1+t)Q_{u;0}(t)\left(F_0[\eta](0)+\int_0^t\|f(s)\|_{H^1_x(g_0)}^{(1)}ds\right).
\end{aligned}
\end{equation}

\textbf{Step 3: Estimate of $E_1^{X,u}[\eta]$ in (\ref{EXun}).}

We differentiate the normal equation in (\ref{LinearizationInComponents}) with the exterior differential operator $d$. This gives
$$
\frac{\partial^2}{\partial t^2}d\phi+b\frac{\partial}{\partial t}d\phi+L^1(\varphi,u)(d\phi)=-[d, L(\varphi,u)]\phi+d\gamma.
$$
A direct calculation gives
$$
\|[d, L(\varphi,u)]\phi\|_{L^2_x(g_0)}
\leq C\lambda\left(\|d\phi\|_{H^1_x(g_0)}+\|\phi\|_{L^2_x(g_0)}\right).
$$
We denote by $\mathcal{P}_0^{\varphi,1}$ the spectral projection of $d\mu_\varphi/d\mu_0(\Delta_{g_\varphi}^1-2)$ for eigenvalue zero (with multiplicity 3), and $\mathcal{P}_\infty^{\varphi,1}=\text{Id}-\mathcal{P}_0^{\varphi,1}$. We still define several weighted energy norms:
$$
F_1^\perp[\phi]
=\left(\|\partial_td\phi\|_{L^2_x(g_0)}
+\langle\mathcal{P}_\infty^{\varphi,1}(d\phi), L^1(\varphi,u)\mathcal{P}_\infty^{\varphi,1}(d\phi)\rangle_{L^2_x(g_0)}
+\lambda\|\mathcal{P}_0^{\varphi,1}(d\phi)\|_{L^2_x(g_0)}\right)^{1/2},
$$
$$
F^\top_0[\psi]
=\left(\|\partial_t\psi\|_{H^2_x(g_0)}^2
+\lambda\|\psi\|_{H^2_x(g_0)}^2\right)^{1/2},
$$
$$
F_1[\eta]=\left(F^\perp_1[\phi]^2+F^\perp_1[\partial_t\phi]^2+F^\top_1[\psi]^2\right)^{1/2}.
$$
The calculation of $D_tE_1^{X,u}[\eta]^2$ is similar as last step and uses results from last step. We first derive
$$
\begin{aligned}
\frac{\partial^2}{\partial t^2}\mathcal{P}_j^{\varphi,1}(d\phi)+b\frac{\partial}{\partial t}\mathcal{P}_j^{\varphi,1}(d\phi)+L^1(\varphi,u)\circ \mathcal{P}_j^{\varphi,1}(d\phi)
=[L^1(\varphi,u), \mathcal{P}_j^{\varphi,1}]d\phi
-\mathcal{P}_j^{\varphi,1}[d, L(\varphi,u)]\phi
+\mathcal{P}_j^{\varphi,1}(d\gamma),
\end{aligned}
$$
where $j=0,\infty$, then differentiate $F_1^\perp[\phi]^2$, use the spectral property guaranteed by Proposition \ref{SpectralProp} and employ estimates (\ref{gammaHn'}) (\ref{Energyn=0}) from last step, to derive a differential inequality:
\begin{equation}\label{EnergyNormal1}
\begin{aligned}
\frac{d}{dt}F_1^\perp[\phi]^2
&\leq C\sqrt{\lambda}F_1[\eta]^2+\|d\gamma\|_{L^2_x(g_0)}\|\partial_td\phi\|_{L^2_x(g_0)}
+\lambda\|\phi\|_{L^2_x(g_0)}\|\partial_td\phi\|_{L^2_x(g_0)}\\
&\leq C\sqrt{\lambda}F_1[\eta]^2
+C\|f\|_{H^1_x(g_0)}F_1[\eta]\\
&\quad+C\sqrt{\lambda}Q_{u;0}(t)\left(F_0[\eta](0)+\int_0^t\|f(s)\|_{H^1_x(g_0)}^{(1)}ds\right)F_1[\eta].
\end{aligned}
\end{equation}
Note that we used the monotonicity of $\lambda(t)$ and inequality $ye^{y}<e^{2y}$. Imitating the calculation of last step and employing  (\ref{Energyn=0}), we also derive differential inequalities
\begin{equation}\label{EnergyTangent1}
\begin{aligned}
\frac{d}{dt}F_1^\top[\psi]^2
&\leq C\sqrt{\lambda}F_1[\eta]^2
+C\|f\|_{H^2_x(g_0)}F_1[\eta]\\
&\quad+C\sqrt{\lambda}Q_{u;0}(t)\left(F_0[\eta](0)+\int_0^t\|f(s)\|_{H^1_x(g_0)}^{(1)}ds\right)F_1[\eta].
\end{aligned}
\end{equation}
\begin{equation}\label{EnergyNormal1'}
\begin{aligned}
\frac{d}{dt}F_1^\top[\partial_t\phi]^2
&\leq C\sqrt{\lambda}F_1[\eta]^2
+C\|f\|_{H^1_x(g_0)}^{(1)}F_1[\eta]\\
&\quad+C\sqrt{\lambda}Q_{u;0}(t)\left(F_0[\eta](0)+\int_0^t\|f(s)\|_{H^1_x(g_0)}^{(1)}ds\right)F_1[\eta].
\end{aligned}
\end{equation}
Adding (\ref{EnergyNormal1})-(\ref{EnergyNormal1'}), we get
$$
\begin{aligned}
\frac{d}{dt}F_1[\eta]
&\leq C\sqrt{\lambda}F_1[\eta]
+C\|f\|_{H^2_x(g_0)}^{(1)}\\
&\quad+C\sqrt{\lambda}Q_{u;0}(t)\left(F_0[\eta](0)+\int_0^t\|f(s)\|_{H^1_x(g_0)}^{(1)}ds\right).
\end{aligned}
$$
Again by Gr?nwall's inequality we obtain, with some $M_2>M_1$
$$
F_1[\eta](t)
\leq M_1\exp\left(M_2\int_0^t\sqrt{\lambda(s)}ds\right)\left(F_1[\eta](0)+E_0^{X,u}[\eta](0)
+\int_0^t\|f(s)\|_{H^2_x}^{(1)}ds\right).
$$
Integrating with respect to $t$, just as in the last step, we obtain
\begin{equation}\label{EnergyNormaln=1}
\begin{aligned}
E_1^{X,u}[\eta](t)
&\leq (1+t)Q_{u;1}(t)\left(E_1^{X,u}[\eta](0)
+\int_0^t\|f(s)\|_{H^2_x}^{(1)}ds\right).
\end{aligned}
\end{equation}

\textbf{Step 4: Estimate of higher derivatives.}

To obtain the energy estimate for general $n$, we differentiate the equation with respect to $x$ with some suitable differential operator. The reason that we do not localize to a coordinate patch is we need some operator that carries the information of Sobolev norm while ``approximately commutes" with $L(\varphi,u)$.

So we let $n$ be an even number and set $G_n=(1-L_\varphi)^{n/2}$. The weighted energy norms will be
$$
F_n[\eta]:=\left(F^\perp_0[G_n\phi]^2+F^\perp_0[\partial_tG_n\phi]^2+\|\partial_t\psi\|_{H^{n+1}_x(g_0)}^2
+\lambda\|\psi\|_{H^{n+1}_x(g_0)}^2\right)^{1/2},
$$
$$
F_{n+1}[\eta]:=\left(F^\perp_1[G_n\phi]^2+F^\perp_1[\partial_tG_n\phi]^2+\|\partial_t\psi\|_{H^{n+2}_x(g_0)}^2
+\lambda\|\psi\|_{H^{n+2}_x(g_0)}^2\right)^{1/2},
$$
\begin{equation}\label{EnergyNormalTemp4}
\frac{\partial^2}{\partial t^2}(G_n\phi)+b\frac{\partial}{\partial t}(G_n\phi)-L(\varphi,u)(G_n\phi)=-[L(\varphi,u),G_n]\phi+G_n\gamma.
\end{equation}
We first notice the following G?rding type inequalities, whose proof is just similar as Proposition \ref{TameEllipticEstimates1}:
$$
\|\phi\|_{H^{n+1}_x(g_0)}^2
\leq
-C\int_{S^2}\mathcal{P}^{\varphi,0}_\infty(G_n\phi)\cdot L(\varphi,u)\mathcal{P}^{\varphi,0}_\infty(G_n\phi)d\mu_0
+C_n\left(1+|X|_{C^{n+2}}\right)\|\phi\|_{L^2_x(g_0)},
$$
$$
\|\phi\|_{H^{n+2}_x(g_0)}^2
\leq
C\int_{S^2}\langle\mathcal{P}^{(1)}_\infty(dG_n\phi), L^1(\varphi,u)\mathcal{P}^{(1)}_\infty(dG_n\phi)\rangle_{g_0}d\mu_0
+C_n\left(1+|X|_{C^{n+2}}\right)\|\phi\|_{L^2_x(g_0)},
$$
implying
\begin{equation}\label{Gaarding}
E_n^{X,u}[\eta]\leq C_nF_n[\eta]+C_n\left(1+|X|_{C^{n+2}}\right)\|\phi\|_{L^2_x(g_0)},\,
E_{n+1}^{X,u}[\eta]\leq C_nF_{n+1}[\eta]+C_n\left(1+|X|_{C^{n+2}}\right)\|\phi\|_{L^2_x(g_0)}.
\end{equation}
Note that the constants do not depend on any higher norm of $X$ and $u$. We also know that $[L(\varphi,u),G_n]$ is a classical differential operator in $x$ of order $n+1$, and we further notice that $L(\varphi,u)$ is a perturbation of $L_{\varphi}$, which commutes with $G_n$. Thus, under a given local coordinate, the coefficients of $[L(\varphi,u),G_n]$ are smooth functions in
$$
\left[D_x^jX,\partial_t^kD_x^{l}u\right]_{0\leq j,l\leq n+2}^{0\leq k\leq 1},
$$
and in fact vanish when $u=0$. Using Lemma \ref{Interpolation}-\ref{PerturbOperator}, inequality (\ref{Gaarding}), and result of step 1, we find the following tame estimate under the similar assumption as step 1.:
$$
\begin{aligned}
&\|[L(\varphi,u),G_n]\phi\|_{L^2_x(g_0)}\\
&\leq C_n\lambda\|\phi\|_{H^{n+1}_x(g_0)}
+C_{n}\left(1+|X|_{C^{n+2}}+|u|_{C_x^{n+2}}^{(1)}\right)\|\phi\|_{L^2_x(g_0)}\\
&\leq C_n\lambda F_n[\eta]
+C_n(1+t)\left(1+|X|_{C^{n+2}}+|u|_{C_x^{n+2}}^{(1)}\right)Q_{u;0}(t)\left(F_0[\eta](0)+\int_0^t\|f(s)\|_{H^1_x(g_0)}^{(1)}ds\right).
\end{aligned}
$$
We then substitute $\phi$ with $G_n\phi$, $\psi$ with $G_n\psi$, $\gamma$ with $-[L(\varphi,u),G_n]\phi+G_n\gamma$ and $\xi$ with $G_n\xi$ in the calculation of last three steps, and thus obtain
\begin{equation}\label{Induction1}
\begin{aligned}
\frac{d}{dt}F_n[\eta]
&\leq C_n\sqrt{\lambda}F_n[\eta]
+C_n\left[\|f\|_{H^{n+1}_x(g_0)}^{(1)}+\left(|X|_{C^{n+4}}+|u|_{C_t^2C_x^{n+4}}\right)\|f\|_{L^2_x(g_0)}^{(1)}\right]\\
&+C_n(1+t)\left(|X|_{C^{n+2}}+|u|_{C^1_tC_x^{n+2}}\right)Q_{u;0}(t)\left(F_0[\eta](0)+\int_0^t\|f(s)\|_{H^1_x(g_0)}^{(1)}ds\right).
\end{aligned}
\end{equation}
By Gr?nwall's inequality we obtain
$$
\begin{aligned}
F_n[\eta]
&\leq Q_{u;n}(t)\left(F_n[\eta](0)+
\int_0^t\left[\|f(s)\|_{H^{n+1}_x(g_0)}^{(1)}+\left(|X|_{C^{n+4}}+|u|_{C_t^2C_x^{n+4}}\right)\|f(s)\|_{L^2_x(g_0)}^{(1)}\right]ds\right)\\
&+(1+t)^2Q_{u;n}(t)\left(|X|_{C^{n+2}}+|u|_{C^1_tC_x^{n+2}}\right)\left(F_0[\eta](0)+\int_0^t\|f(s)\|_{H^1_x(g_0)}^{(1)}ds\right).
\end{aligned}
$$
Taking into account the G?rding type inequality (\ref{Gaarding}) and estimate  (\ref{Energyn=0}), we obtain the desired result for even $n$.

We then differentiate with the exterior differential operator $d$ and obtain equations for $dG_n\phi$ and $dG_n\psi$. We can derive the desired estimate for $F_{n+1}[\eta]$, hence $E_{n+1}^{X,u}[\eta]$, similarly as above.

Finally, to obtain estimates for the frame-independent energy norm $E_n[\eta]$, it suffices to notice that by Lemma \ref{GeometricQuantities} and the equation $\Psi'(w)\eta=f$ satisfied by $\eta$,
$$
E_{n+1}[\eta]\leq C_nE_{n}^{X,u}[\eta]+C_n\left(|X|_{C^{n+3}}+|u|_{C_t^2C_x^{n+3}}\right)E_{2}[\eta],
$$
$$
E_{n}^{X,u}[\eta]\leq C_nE_{n+2}[\eta]+C_n\left(|X|_{C^{n+3}}+|u|_{C_t^2C_x^{n+3}}\right)E_{2}[\eta]+C_n\|f\|_{H^{n}_x(g_0)}.
$$
\end{proof}

\subsection{Tame Decay Estimate for the Damped Linearized Equation}
We use results from last section to solve the  equation
$$
\Phi'(X,a,u)(Y,c,v)=f,
$$
in the space $\mathbf{F}_\beta$, where the damping $b>0$ is fixed and the $\beta>0$ given by (\ref{OmegaB}). The idea is simple: energy estimates in Proposition \ref{EnergyFull} ensure that the growth of the norms $\|\phi(t)\|_{H^n_x(g_0)}$, and $\|\psi(t)\|_{H^n_x(g_0)}$ is at most exponentially fast, and if $|X|_{C^4}+|u|_{C^3_tC_x^4}$ is sufficiently small compared to $b$, then with $w=i_\varphi+a+u$ and $\phi N(w)+\psi=\Sigma_\varphi  Y+c+v$, we find that $\phi$ and $\psi$ solve equations of the form
$$
\partial_t^2\phi+b\partial_t\phi
=L_\varphi\phi+\text{Exponentially decaying terms},
$$
$$
\partial_t^2\psi+b\partial_t\psi
=\text{Exponentially decaying terms}.
$$
So using results from Section \ref{Sec3}, we find that $\phi$ and $\psi$ actually converge exponentially in $H^{n-2}_x(g_0)$, with a slower rate compared to $\beta$ (as defined in (\ref{OmegaB})); we can then use a bootstrapping argument until we arrive at a satisfactory rate of exponential decay, with a controllable loss of spatial derivatives.

We turn to the details.

\begin{proposition}\label{DecayOfLEQ}
Fix $b>0$ and $\beta$ as in (\ref{OmegaB}). Suppose $(X,a,u)\in\mathbf{F}_\beta$, and still write $\varphi=\mathcal{E}_X$. There is a $\delta_5>0$ such that if
$$
\|X\|_{H^{12}}+\|u\|_{\beta,3;12}<\delta_5,
$$
then assumptions of Proposition \ref{EnergyFull} is satisfied, and equation $\Phi'(X,a,u)(Y,c,v)=f$, with Cauchy data of $\eta:=\Sigma_\varphi Y+c+v$ given, has a unique solution $(Y,c,v)\in\mathbf{F}_\beta$, satisfying the following tame estimate with respect to the grading of $\mathbf{F}_\beta$:
$$
\|(Y,c,v)\|_{n}
\leq C_nE_{n+8}[\eta](0)+\|(X,a,u)\|_{n+12}E_3[\eta](0)
+C_n\|f\|_{\beta,1;n+7}+\|(X,a,u)\|_{n+12}\|f\|_{\beta,1;4}.
$$
\end{proposition}
\begin{proof}
Still as before, all computations below will be on a given time slice $t$, unless otherwise noted.

By the Sobolev embedding $H^2_x\hookrightarrow C_x^{1-\varepsilon}$, if $\delta_5$ is sufficiently small, then $|X|_{C^4}+|u|_{C_t^3C_x^4}$ automatically satisfies the requirement of Proposition \ref{EnergyFull}.

We shall follow the notation of Proposition \ref{EnergyFull}, i.e. we set $\varphi=\mathcal{E}_X$, $\chi=i_\varphi+a$, $\eta=\Sigma_\varphi Y+c+v$. The equation for $\eta$ is $\Psi'(\chi+u)\eta=f$, which, according to \cite{Notz2010}, has a unique solution $\eta\in C^3([0,\infty);C^\infty(S^2;\mathbb{R}^3))$, and Proposition \ref{EnergyFull} gives a tame energy estimate of the solution.
We now define $\phi=\perp_\chi\eta$, $\psi=\top_\chi\eta$, and
$$
y=f+[\Psi'(\chi)-\Psi'(\chi+u)](\eta),
$$
and thus re-write the equation in terms of $\eta$ as $\Psi'(\chi)(\eta)=y$. Writing in components, this is exactly of the form indicated in Proposition \ref{Linearized(u=0)}. Note that these symbols signify different objects as in Proposition \ref{EnergyFull}.

Let us estimate the decay of $y$. We first notice that Proposition \ref{EnergyFull}, with the aid of Lemma \ref{GeometricQuantities},  gives the following estimate:
$$
\begin{aligned}
\|\eta\|_{H^{n}_x}^{(1)}
&=E_n[\eta]\\
&\leq(1+t)^2Q_{u;n}(t)\left[E_{n+1}[\eta](0)+
\left(|X|_{C^{n+3}}+|u|_{C^3_tC^{n+3}_x}\right)E_2[\eta](0)\right]\\
&\quad+(1+t)^2Q_{u;n}(t)\left[\|f(0)\|_{H^n_x(g_0)}+\left(|X|_{C^{n+3}}+|u|_{C^3_tC^{n+3}_x}\right)\|f(0)\|_{H^{2}_x(g_0)}\right]\\
&\quad+(1+t)^2Q_{u;n}(t)\left[\int_0^t\|f(s)\|_{H^{n}_x(g_0)}^{(1)}ds
+\left(|X|_{C^{n+3}}+|u|_{C^3_tC^{n+3}_x}\right)
\int_0^t\|f(s)\|_{L^2_x(g_0)}^{(1)}ds\right].
\end{aligned}
$$
But $u\in\mathbf{E}_{\beta,3}$, $f\in\mathbf{E}_{\beta,1}$, so
$$
\begin{aligned}
Q_{u;n}(t)&=M_n\exp\left(M_n
\int_0^t\sqrt{\sup_{\tau\geq s}|u(\tau)|_{C_x^4}^{(3)}}ds
\right)\\
&\leq M_n\exp\left(M_n
\int_0^t e^{-\beta s/2}\|u\|_{\beta,3;6}^{1/2}ds
\right)\\
&\leq M_n e^{M_n\delta_5},
\end{aligned}
$$
and the integrals of norm of $f$ can be estimated similarly. So we in fact have
\begin{equation}\label{EtaHn}
\begin{aligned}
\|\eta(t)\|_{H^{n}_x}^{(1)}
&\leq C_n(1+t)^2\left[E_{n+1}[\eta](0)+
\left(|X|_{C^{n+3}}+\|u\|_{\beta,3;n+5}\right)E_2[\eta](0)\right]\\
&\quad+C_n(1+t)^2\left[\|f\|_{\beta,1;n}
+\left(|X|_{C^{n+3}}+\|u\|_{\beta,3;n+5}\right)
\|f\|_{\beta,1;2}\right].
\end{aligned}
\end{equation}

We know from (\ref{LEQ1}) that $\Psi'(\chi)-\Psi'(\chi+u)$ is a second order \emph{spatial} differential operator, whose coefficients depend smoothly on up to second order time derivatives and third order spatial derivatives of $X$ and $u$, and vanishes for $u=0$. Thus by (\ref{EtaHn}), using Lemma \ref{Interpolation}--\ref{PerturbOperator},
$$
\begin{aligned}
\|y(t)\|_{H^n_x(g_0)}^{(1)}
&\leq \|f\|_{H^n_x(g_0)}^{(1)}
+C_ne^{-\beta t}\|\eta\|_{H^{n+2}_x(g_0)}^{(1)}
+C_ne^{-\beta t}
\left(1+|X|_{C^{n+5}}+\|u\|_{\beta,3;n+7}\right)\|\eta\|_{L^2_x(g_0)}^{(1)}\\
&\leq C_ne^{-\beta t/2}\left[E_{n+3}[\eta](0)+
\left(|X|_{C^{n+5}}+\|u\|_{\beta,3;n+7}\right)E_2[\eta](0)
\right]\\
&\quad+C_ne^{-\beta t/2}\left[\|f\|_{\beta,1;n+2}
+\left(|X|_{C^{n+5}}+\|u\|_{\beta,3;n+7}\right)
\|f\|_{\beta,1;2}\right].
\end{aligned}
$$
So $\|y\|_{\beta/2,1;n}$ is estimated.

We can then start a bootstrap argument. First we apply Proposition \ref{Linearized(u=0)} to $\Psi'(\chi)\eta=y$, with $\beta/2$ in place of $\beta$. This gives a unique solution $(Y,c,v)\in\mathbf{F}_{\beta/2}$ such that
$$
\eta=\Sigma_\varphi  Y+c+v,\,v\in\mathbf{E}_{\beta/2,3},
$$
and the tame estimates for $Y$ and $c$ are already guaranteed by Proposition \ref{Linearized(u=0)}:
$$
\begin{aligned}
|c|
\leq C\left(E_1[\eta](0)+\|y\|_{\beta/2,0;0}\right)
\leq C\left(E_3[\eta](0)+\|f\|_{\beta,1;2}\right),
\end{aligned}
$$
$$
\begin{aligned}
\|Y\|_{H^n}
&\leq C_n\left(E_{n}[\eta](0)+\|y\|_{\beta/2,0;n}\right)
+C_n|X|_{C^{n+1}}\left(E_2[\eta](0)+\|y\|_{\beta/2,0;2}\right)\\
&\leq C_n\left[E_{n+3}[\eta](0)+
\left(|X|_{C^{n+5}}+\|u\|_{\beta,3;n+7}\right)E_{3}[\eta](0)
\right]\\
&\quad+C_n\left[\|f\|_{\beta,1;n+2}
+\left(|X|_{C^{n+5}}+\|u\|_{\beta,3;n+7}\right)
\|f\|_{\beta,1;4}\right],
\end{aligned}
$$
and there is a $e^{-\beta t/2}$-decay estimate for $v$:
$$
\begin{aligned}
\|v\|_{\beta/2,1;n}
&\leq C_n\left(E_{n+1}[\eta](0)+\|y\|_{\beta/2,0;n}\right)
+C_n|X|_{C^{n+3}}\left(E_1[\eta](0)+\|y\|_{\beta/2,0;2}\right)\\
&\leq C_n\left[E_{n+3}[\eta](0)+
\left(|X|_{C^{n+5}}+\|u\|_{\beta,3;n+7}\right)E_2[\eta](0)
\right]\\
&\quad+C_n\left[\|f\|_{\beta,1;n+2}
+\left(|X|_{C^{n+5}}+\|u\|_{\beta,3;n+7}\right)
\|f\|_{\beta,1;2}\right].
\end{aligned}
$$
Thus we obtain, with the aid of Lemma \ref{Interpolation}--\ref{PerturbOperator}, that on each time slice we in fact have $\|[\Psi'(\chi)-\Psi'(w+u)](\eta)\|_{H^n_x(g_0)}^{(1)}=O(e^{-\beta t})$, and to be precise,
$$
\begin{aligned}
&\|[\Psi'(\chi)-\Psi'(\chi+u)](\eta)\|_{\beta,1;n}\\
&\leq C_n\left(|Y|_{C^{n+2}}+\|v\|_{\beta/2,1;n+2}\right)
+C_n\left(|X|_{C^{n+2}}+\|u\|_{\beta,1;n+4}\right)\left(|Y|_{C^{2}}+\|v\|_{\beta/2,1;2}\right)\\
&\leq C_n\left[E_{n+7}[\eta](0)
+\left(|X|_{C^{n+9}}+\|u\|_{\beta,3;n+11}\right)E_{3}[\eta](0)\right]\\
&\quad+C_n\left[\|f\|_{\beta,1;n+6}
+\left(|X|_{C^{n+9}}+\|u\|_{\beta,3;n+11}\right)
\|f\|_{\beta,1;4}\right].
\end{aligned}
$$
Thus, applying Proposition \ref{Linearized(u=0)} to the equation $\Phi'(X,a,u)(Y,c,v)=[\Psi'(\chi)-\Psi'(w+u)](\eta)$ again, we finally obtain
$$
\begin{aligned}
\|v\|_{\beta,3;n}
&\leq C_n\left[E_{n+8}[\eta](0)
+\left(|X|_{C^{n+10}}+\|u\|_{\beta,3;n+12}\right)E_{3}[\eta](0)\right]\\
&\quad+C_n\left[\|f\|_{\beta,1;n+7}
+\left(|X|_{C^{n+10}}+\|u\|_{\beta,3;n+12}\right)
\|f\|_{\beta,1;4}\right].
\end{aligned}
$$
This completes the proof.
\end{proof}

\section{Proof of Theorem \ref{Thm1} and Theorem \ref{Thm2}}
In this section, we complete the proof of Theorem \ref{Thm1} and \ref{Thm2}. Before proceeding to the proof, we provide the details of the Nash-Moser scheme that we shall employ. The structures were first constructed by H?rmander \cite{Hormander1985}, and then refined by Baldi and Haus in \cite{BaldiHaus2017}. To avoid confusion with symbols in our paper, the notation employed in the quotation below will be different from \cite{BaldiHaus2017}.

Let $(E_a)_{a\geq0}$ be a family of decreasing Banach spaces, with continuous injection $E_b\hookrightarrow  E_a$ for $b\geq a$ satisfying
$$
\|u\|_a\leq\|u\|_b.
$$
Let $E_\infty=\cap_{a\geq0}E_a$ and equip it with the weakest topology making $E_\infty\hookrightarrow  E_a$ continuous for each $a$. Assume further the existence of a family of smoothing operators $S_j:E_a\to E_\infty$ for $j=1,2,...$, satisfying
\begin{itemize}
\item For each $a\geq0$,
$$
\|S_ju\|_a\leq C_a\|u\|_a.
$$

\item For $a<b$,
$$
\|S_ju\|_b\leq C_{a,b}2^{j(b-a)}\|u\|_a.
$$

\item For $a>b$,
$$
\|(1-S_j)u\|_b\leq C_{a,b}2^{-j(a-b)}\|u\|_a.
$$

\item Set $R_j=S_{j+1}-S_j$. Then for any $a,b$,
$$
\|u\|_a^2\leq C_a\sum_{j=0}^\infty\|R_ju\|_a^2,
$$
$$
\|R_ju\|_b\leq C_{a,b}2^{j(b-a)}\|R_ju\|_a.
$$
\end{itemize}

We now quote from \cite{BaldiHaus2017} the following Nash-Moser-H?rmander theorem:

\begin{theorem}\label{NMH}
Let $(E_a)_{a\geq0}$ and $(F_a)_{a\geq0}$ be decreasing scales of Banach spaces satisfying above requirements. Suppose  $a_0,a_1,a_2,\rho,\mu,\lambda\geq0$ satisfy
\begin{equation}\label{NMH1}
a_0\leq\mu\leq a_1,\,
a_1+\frac{\lambda}{2}<\rho<a_2+\lambda,\,
2\rho<a_1+a_2.
\end{equation}
Let $V$ be a convex neighbourhood of $0$ in $E_\mu$. Let $\Phi:V\to F_0$ be a map, such that for any $a\in[0,a_2-\mu]$, $\Phi:V\cap E_{a+\mu}\to F_a$ is $C^2$, and for all $u\in V\cap E_{a+\mu}$, $a\in[0,a_2-\mu]$ there holds
\begin{equation}\label{NMH2}
\begin{aligned}
\|\Phi''(u)[v,w]\|_a
&\leq M_1(a)\left(\|v\|_{a+\mu}\|w\|_{a_0}+\|v\|_{a_0}\|w\|_{a+\mu}\right)
+\left(M_2(a)\|u\|_{a+\mu}+M_3(a)\right)\|v\|_{a_0}\|w\|_{a_0},
\end{aligned}
\end{equation}
where $M_i(a)$ are positive increasing functions defined for $a\geq0$. Assume further there is an $\varepsilon_1>0$ such that for any $v\in V\cap E_\infty$ with $\|v\|_{a_1}\leq{\varepsilon_1}$, the linear mapping $\Phi'(v)$ has a right inverse $\Omega(v):F_\infty\to E_{a_2}$, such that for any $a\in[a_1,a_2]$,
\begin{equation}\label{NMH3}
\|\Omega(v)f\|_a
\leq Q_1(a)\|f\|_{a+\lambda-\rho}+\left(Q_2(a)\|v\|_{a+\lambda}+Q_3(a)\right)\|f\|_{0},
\end{equation}
where $Q_i(a)$ are positive increasing functions defined for $a\geq0$.

Then for any $A>0$, there exists an $\varepsilon>0$, such that for any $f\in F_\lambda$ satisfying
\begin{equation}\label{NMH4}
\sum_{j=0}^\infty\|R_jf\|_{\lambda}^2\leq A^2\|f\|_\lambda^2,\,\|f\|_\lambda\leq\varepsilon,
\end{equation}
there exists a $u\in E_\rho$ solving $\Phi(u)=\Phi(0)+f$, such that
\begin{equation}\label{NMH5}
\|u\|_\rho\leq CQ_{123}(a_2)(1+A)\|f\|_\lambda,
\end{equation}
where $Q_{123}=Q_1+Q_2+Q_3$, and $C$ depends on $a_1,a_2,\rho,\lambda$. The $\varepsilon$ is explicitly given by
\begin{equation}\label{NMH6}
\varepsilon^{-1}=C'Q_{123}(a_2)(1+A)\max\left[
1,\,\frac{1}{\varepsilon_1},\,Q_{123}(a_2)M_{123}(a_2-\mu)
\right],
\end{equation}
where $M_{123}=M_1+M_2+M_3$, and $C'$ depends on $a_1,a_2,\rho,\lambda$.

If, in addition, for some $c>0$ (\ref{NMH2}) holds for any $a\in[0,a_2+c-\mu]$, and $\Omega(v)$ maps $F_\infty$ to $E_{a_2+c}$ with (\ref{NMH3}) holding for any $a\in[a_1,a_2+c]$, and $f\in F_{\lambda+c}$ with
$$
\sum_{j=0}^\infty\|R_jf\|_{\lambda+c}^2\leq A_c^2\|f\|_{\lambda+c}^2,\,\|f\|_\lambda\leq\varepsilon,
$$
then the solution $u$ is in $E_{\rho+c}$, with
$$
\|u\|_{\rho+c}\leq C_c\left(K_1(1+A)\|f\|_\lambda+K_2(1+A)\|f\|_{\lambda+c}\right),
$$
where
$$
K_1=\bar Q_3+\bar{Q}_{12}(\bar{Q}_3\bar{M}_{12}+Q_{123}(a_2)\bar{M}_3)\sum_{j\leq N-2}z^j,
$$
$$
K_2=\bar{Q}_{12}\sum_{j\leq N-1}z^j.
$$
Here $\bar{Q}_{12}=\bar{Q}_1+\bar{Q}_2$, $\bar{Q}_i(a)=Q_i(a+c)$; $\bar{M}_{12}=\bar{M_1}+\bar{M_2}$, $\bar{M_i}=M_i(a+c-\mu)$; $C_c,N$ depend on $a_1,a_2,\rho,\lambda,c$, and
$$
z=Q_{123}(a_1)M_{123}(0)+\bar{Q}_{12}\bar{M}_{12}.
$$
\end{theorem}

This is a refinement of H?rmander's version of Nash-Moser theorem in \cite{Hormander1985}. We point out several advantages of Theorem \ref{NMH} compared to the ``simplest" version in \cite{SR1989} or the ``structuralist" version in \cite{Hamilton19821}. First of all, although the theorem is concerned with a nonlinear operator in the category of tame Fréchet spaces, the statement in fact depends only on scales in a \emph{finite} interval; this fact is better illustrated in the rougher version in \cite{SR1989}. The next advantage is that it provides an explicit linear estimate on the size of the solution in terms of the known, and avoids the presence of Besov spaces as in \cite{Hormander1985}. As pointed out by Baldi and Haus, this estimate keeps the optimal loss of regularity.

We turn to the proof of our main results. The smoothing operators $S_j$ will be constructed through (\ref{smoothing0}) and (\ref{smoothing1}), i.e.,
$$
S_j f:=\sum_{\lambda\in\sigma[-\Delta_{g_0}]:\lambda\leq2^j}\mathcal{Q}_\lambda f,
\quad
S_j X
:=\sum_{\lambda\in\sigma(D_{g_0}^*D_{g_0}),\lambda\leq2^j}\mathcal{Q}^{(1)}_\lambda X,
$$
Under the $H^n$-grading, it is easily verified that the requirements for smoothing operators are all satisfied, since there holds
$$
\|f\|_{H^s_x(g_0)}^2
\simeq_{s}\|\mathcal{Q}_0 f\|_{L^2_x(g_0)}^2
+\sum_{\lambda\in\sigma[-\Delta_{g_0}]}\lambda^{2s}\|\mathcal{Q}_\lambda f\|_{L^2_x(g_0)}^2,
$$
$$
\|X\|_{H^s_x(g_0)}^2
\simeq_{s}\|\mathcal{Q}^{(1)}_0 X\|_{L^2_x(g_0)}^2
+\sum_{\lambda\in\sigma[D_{g_0}^*D_{g_0}]}\lambda^{2s}\|\mathcal{Q}^{(1)}_\lambda X\|_{L^2_x(g_0)}^2,
$$
where the left-hand-sides are defined via local coordinates.

\begin{proof}[Proof of Theorem \ref{Thm1}]\label{ProofThm1}
Define a mapping $\mathfrak{P}:\mathbf{F}_\beta\to C^\infty(S^2;\mathbb{R}^3)\oplus C^\infty(S^2;\mathbb{R}^3)\oplus
\mathbf{E}_{\beta,1}$ by
$$
\mathfrak{P}(X,a,u)=\left(\begin{matrix}
\Xi(X,a,u(0))-i_0\\
\partial_tu(0)\\
\Phi(X,a,u)
\end{matrix}
\right)
=\left(\begin{matrix}
i_0\circ\mathcal E_X+a+u(0)-i_0\\
\partial_tu(0)\\
\Phi(X,a,u)
\end{matrix}
\right).
$$
Solving $\Phi(X,a,u)=0$ with initial data $u_0,u_1$ is equivalent to solving the following equation:
\begin{equation}\label{EQWDV}
\begin{aligned}
\mathfrak{P}(X,a,u)=\left(\begin{matrix}
u_0-i_0\\
u_1\\
0
\end{matrix}
\right).
\end{aligned}
\end{equation}

That $\mathfrak{P}$ is a smooth tame map is easily verified. The loss of regularity caused by $\mathfrak{P}$ has order 2 since $\Phi$ is a second order differential operator; note that the time differentiation here does not cause any loss of regularity since at every scale $n$, $\partial_t^2:\mathbf{E}_{\beta,3}^n\to \mathbf{E}_{\beta,1}^n$ is a continuous linear mapping between Banach spaces.

The linearization of $\mathfrak{P}$ is
$$
\mathfrak{P}'(X,a,u)(Y,c,v)
=\left(\begin{matrix}
\Sigma_\varphi  Y+c+v(0)\\
\partial_tv(0)\\
\Phi'(X,a,u)(Y,c,v)
\end{matrix}
\right).
$$
By Proposition \ref{DecayOfLEQ}, given any $f_0,f_1\in C^\infty(S^2;\mathbb{R}^3)$, $f\in\mathbf{E}_{\beta,1}$, the equation
$$
\mathfrak{P}'(X,a,u)(Y,c,v)
=\left(\begin{matrix}
f_0\\
f_1\\
f
\end{matrix}
\right)
$$
has a unique solution $(Y,c,v)\in\mathbf{F}_\beta$, satisfying the tame estimate
\begin{equation}\label{tame1}
\begin{aligned}
\|(Y,c,v)\|_n
&\leq C_n\left(\|(f_0,f_1)\|_{H^{n+8}(g_0)}+\|f\|_{\beta,1;n+7}\right)
+C_n\|(X,a,u)\|_{n+12}\left(\|(f_0,f_1)\|_{H^{4}(g_0)}+\|f\|_{\beta,1;4}\right).
\end{aligned}
\end{equation}

We may thus apply the Nash-Moser-H?rmander Theorem \ref{NMH}. Here $a$ takes value in $\mathbb{N}$ (which does not affect the argument since the Nash-Moser theorem essentially does not require the scales to vary continuously), $E_a=H^{a}_{\mathfrak{X}}\oplus\mathbb{R}^3\oplus\mathbf{E}_{\beta,3}^{a}$ (where $H^{a}_{\mathfrak{X}}$ denotes all $H^{a}$-vector fields), $F_a=H^a(S^2;\mathbb{R}^3)\oplus H^a(S^2;\mathbb{R}^3)\oplus\mathbf{E}_{\beta,1}^{a}$, the convex neighbourhood $V$ to be an open set in $E_6$ such that for any $(X,a,u)\in V$, the mapping $i_0\circ\mathcal{E}_X+a+u$ is a $C^3$ embedding of $S^2$ into $\mathbb{R}^3$. We shall then take $a_0=2$, $\mu=4$, $a_1=12$, $\lambda=41$, $\rho=33$, $a_2=55$, and $\varepsilon_1$ equal to the $\delta_5$ in Proposition \ref{DecayOfLEQ}. Using Lemma \ref{Interpolation}--\ref{GeometricQuantities}, (\ref{NMH1}) and (\ref{NMH2}) are satisfied. According to Proposition \ref{DecayOfLEQ}, (\ref{NMH3}) is satisfied, and by our choice of smoothing operators, (\ref{NMH4}) is automatically satisfied. The Nash-Moser-H?rmander theorem then ensures the existence of solution $(X,a,u)\in E_{33}$ to equation (\ref{EQWDV}) if the magnitude of $\|u_0\|_{H^{41}},\|u_1\|_{H^{41}}$ is small as indicated in (\ref{NMH6}). Uniqueness of solution is ensured by the local uniqueness results established in \cite{Notz2010}. For higher regularity, it suffices to apply the general higher regularity results in Theorem \ref{NMH}.
\end{proof}

\begin{proof}[Proof of Theorem \ref{Thm2}]\label{ProofThm2}
We leave the time scale $T=T_\varepsilon$ undetermined for the moment.

Set $E_a=C^3([0,T];H^a(S^2;\mathbb{R}^3))$, $F_a=H^a(S^2;\mathbb{R}^3)\oplus H^a(S^2;\mathbb{R}^3)\oplus C^1([0,T];H^a(S^2;\mathbb{R}^3))$.
Define a mapping $\mathfrak{Q}:E_\infty\to F_\infty$ as follows:
$$
\mathfrak{Q}(u)
=\left(\begin{matrix}
u(0)-i_0\\
\partial_tu(0)\\
\Psi(i_0+u)
\end{matrix}
\right).
$$
Then the Cauchy problem
$$
\frac{\partial^2u}{\partial t^2}
=\frac{d\mu(i_0+u)}{d\mu_0}\left(-H(i_0+u)+\frac{\kappa}{\mathrm{Vol}(i_0+u)}\right)N(i_0+u),\,
\left(
\begin{matrix}
i_0(x)+u(0,x)\\
\partial_tu(0,x)
\end{matrix}
\right)=
\left(
\begin{matrix}
u_0(x)\\
u_1(x)
\end{matrix}
\right)
$$
is equivalent to the equation
\begin{equation}\label{EQV}
\mathfrak{Q}(u)
=\left(\begin{matrix}
u_0-i_0\\
u_1\\
0
\end{matrix}
\right).
\end{equation}

For $f_0,f_1\in C^\infty(S^2;\mathbb{R}^3)$, $f\in C^3([0,T];H^a(S^2;\mathbb{R}^3))$ the linearization of equation (\ref{EQV}) is
$$
\left(\begin{matrix}
\eta(0)\\
\partial_t\eta(0)\\
\Psi'(i_0+u)\eta
\end{matrix}
\right)
=\left(\begin{matrix}
f_0\\
f_1\\
f
\end{matrix}
\right),
$$
and by Proposition \ref{EnergyFull}, there is an $\delta_5>0$ such that if $\|u\|_{C^3_tH^{7}_x}<\varepsilon_1<\delta_5$, then it has a unique solution $\eta\in E_\infty$, satisfying the energy estimate
$$
\begin{aligned}
\|\eta\|_{C_t^1H^n_x}
&\leq C_n(1+T)^2Q_n(u;T)\left[\|(f_0,f_1)\|_{H^{n+1}(g_0)}
+\|u\|_{C^3_tH^{n+5}_x}\|(f_0,f_1)\|_{L^2(g_0)}\right]\\
&\quad+C_n(1+T)^3Q_n(u;T)\left[\|f\|_{C_t^1H^n_x}+\|u\|_{C^3_tH^{n+5}_x}\|f\|_{C_t^1H^2_x}\right],
\end{aligned}
$$
where
$$
Q_n(u;T)
=\exp\left(C_nT\sqrt{\|u\|_{C^3_tH^{6}_x(g_0)}}\right)
<\exp\left(C_nT\varepsilon_1^{1/2}\right).
$$
Using the equation $\Psi'(i_0+u)\eta=0$ itself and differentiating it with respect to $t$, it is not hard to see that
$$
\begin{aligned}
\|\eta\|_{C_t^3H^n_x}
&\leq C_n(1+T)^2Q_n(u;T)\left[\|(f_0,f_1)\|_{H^{n+3}(g_0)}
+\|u\|_{C^3_tH^{n+7}_x}\|(f_0,f_1)\|_{L^2(g_0)}\right]\\
&\quad+C_n(1+T)^3Q_n(u;T)\left[\|f\|_{C_t^1H^{n+2}_x}+\|u\|_{C^3_tH^{n+5}_x}\|f\|_{C_t^1H^2_x}\right].
\end{aligned}
$$

We now apply the Nash-Moser-H?rmander theorem. Let the convex neighbourhood $V$ be an open set in $E_6$ such that $i_0+u$ is a $C^2$-embedding from $S^3$ to $\mathbb{R}^3$. We choose $a_0=2$, $\mu=2$, $a_1=7$, $\lambda=24$, $\rho=21$, $a_2=43$. By our choice of smoothing operators, (\ref{NMH4}) is automatically satisfied. Using Lemma \ref{Interpolation}--\ref{GeometricQuantities}, (\ref{NMH1}) and (\ref{NMH2}) are satisfied. According to Proposition \ref{EnergyFull}, (\ref{NMH3}) is satisfied, with the $M_i$'s being independent of $T$, and
$$
Q_i(a)=C_a(1+T)^3\exp\left(C_aT\varepsilon_1^{1/2}\right).
$$
The Nash-Moser-H?rmander theorem then ensures the existence of solution $u\in E_{21}$ to equation (\ref{EQV}) provided that $\|(u_0-i_0,u_1)\|_{H^{24}(g_0)}<\varepsilon$, where
\begin{equation}{\label{Te}}
\varepsilon^{-1}=C(1+T)^3\exp\left(CT\varepsilon_1^{1/2}\right)\max\left[1,\frac{1}{\varepsilon_1},C(1+T)^3\exp\left(CT\varepsilon_1^{1/2}\right)
\right].
\end{equation}

We can now determine the optimal $T$ from (\ref{Te}). Since we assume that $\varepsilon,\varepsilon_1$ should be small, we expect that $T$ is large. By (\ref{Te}) we find that
$$
\varepsilon^{-1}\geq C(1+T)^6\exp\left(CT\varepsilon_1^{1/2}\right),
$$
so there must hold $T=O(\varepsilon^{-1/6})$. It is easily verified that $T\simeq\varepsilon^{-1/6}$ does work if we choose $\varepsilon_1=\varepsilon^{1/3}$ since $\varepsilon^{-1/3}$ grows slower than $\varepsilon^{-1/2}$ as $\varepsilon\to0$. We thus conclude that $T\simeq\varepsilon^{-1/6}$ is a lower bound for life span if $\|(u_0-i_0,u_1)\|_{H^{24}(g_0)}<\varepsilon$ and $\varepsilon$ is sufficiently small. As for uniqueness, \cite{Notz2010} still provides the desired uniqueness result. Higher regularity results follow from Theorem \ref{NMH}, just as in the proof of Theorem \ref{Thm1}. The growth estimate indicated in the statement of Theorem \ref{Thm2} follows directly from the last half of Theorem \ref{NMH}.
\end{proof}

\appendix

\section{Weakly Hyperbolic Systems}\label{A}
In this apendix, we will sketch results on weakly hyperbolic systems and local well-posedness of (\ref{EQ0}) and (\ref{EQWD}).

Let $M$ be a compact differential manifold, and let $E_1,E_2$ be smooth vector bundles on $M$. Let $\phi$ and $\psi$ be time-dependent sections of $E_1$ and $E_2$ respectively. A \emph{weakly hyperbolic linear system in $(\phi,\psi)$} takes the following form:
\begin{equation}\label{WHLS}
\begin{aligned}
\frac{\partial^2}{\partial t^2}\phi(t)
&=L(t)\phi(t)+M(t)\partial_t\phi(t)+P(t)\psi(t)+f_1(t),\\
\frac{\partial^2}{\partial t^2}\psi(t)
&=Q(t)\phi(t)+R(t)\psi(t)+f_2(t),
\end{aligned}
\end{equation}
where given any time $t$, $L(t):\Gamma(E_1)\to\Gamma(E_1)$ is a second order elliptic operator, $P(t):\Gamma(E_2)\to\Gamma(E_2)$ and  $Q(t):\Gamma(E_1)\to\Gamma(E_1)$ are first order pseudo-differential operators, $M(t):\Gamma(E_1)\to\Gamma(E_1)$ and  $R(t):\Gamma(E_2)\to\Gamma(E_2)$ are zeroth order pseudo-differential operators, and $f_1,f_2$ are known time-dependent sections of $E_1$ and $E_2$. The notion of weakly hyperbolic linear system is just a rephrasing of Hamilton's notion of weakly parabolic linear system in \cite{Hamilton19822}.

\begin{proposition}\label{SolvingWHLS}
Fix Riemannian metrics on $E_1,E_2$ to measure function norms. Given smooth initial data $\phi[0],\psi[0]$ to (\ref{WHLS}), the equation is uniquely solvable on any time interval $[0,T]$, and the solution satisfies the tame energy estimate
$$
\begin{aligned}
E_n(t)
&\leq Ce^{Ct}E_n(0)
+C\int_0^te^{C(t-s)}\left(\|f_1(s)\|_{H^{n}_x}^{(1)}+\|f_2(s)\|_{H^{n+1}_x}\right)ds\\
&+C\int_0^te^{C(t-s)}\left([L(s)]_n^{(1)}+[M(s)]_n^{(1)}+[P(s)]_n^{(1)}+[Q(s)]_{n+1}+[R(s)]_{n+1}\right)\left(E_{n_0}(0)+\|f_1(s)\|_{H^{n_0}_x}^{(1)}\right)ds,
\end{aligned}
$$
where $\|f(s)\|^{(1)}:=\|f(s)\|+\|\partial_sf(s)\|$, the energy norm
$$
E_n(t):=\|\partial_t^2\phi(t)\|_{H^n_x}
+\|\partial_t\phi(t)\|_{H^{n+1}_x}
+\|\phi(t)\|_{H^{n+1}_x}+\|\partial_t\psi(t)\|_{H^{n+1}_x}+\|\psi(t)\|_{H^{n+1}_x},
$$
and $[\cdot]_n$ denotes the Sobolev norm on jet bundles. The integer $n_0$ is the smallest integer to ensure Sobolev embedding $H^{n_0}\hookrightarrow C^2$. The constants $C$ depend on $n$, the elliptic constants of $\{L(t)\}_{t\in[0,T]}$, and
$$
\sup_{t\in[0,T]}\left([L(t)]_{n_0}^{(1)}+[M(t)]_{n_0}^{(1)}+[P(t)]_{n_0}^{(1)}+[Q(t)]_{n_0}+[R(t)]_{n_0}\right).
$$
\end{proposition}

This is enough for establishing the local well-posedness result for (\ref{EQ0}) and (\ref{EQWD}), defined on a general compact surface $M$. In fact, these two equations both fall into the class of ``evolutionary problems with an integrability condition", treated by Hamilton in \cite{Hamilton19822}. The linearization of both gives a weakly hyperbolic system, whose unique solution satisfies a tame estimate. By the Nash-Moser theorem, the following is true:
\begin{theorem}
Let $M$ be a compact oriented surface, and let $u_0:M\hookrightarrow \mathbb{R}^3$ be a smooth embedding. Fix constants $\kappa>0$, $b\geq0$. Let $u_1:M\to\mathbb{R}^3$ be any smooth mapping. Then there exists a $T>0$ depending on $u_0,u_1$ such that the Cauchy problem
$$
\frac{\partial^2u}{\partial t^2}+b\frac{\partial u}{\partial t}
=\frac{d\mu(u)}{d\mu_0}\left(-H(u)+\frac{\kappa}{\mathrm{Vol}(u)}\right)N(u),\,
\left(
\begin{matrix}
u(0,x)\\
\partial_tu(0,x)
\end{matrix}
\right)=
\left(
\begin{matrix}
u_0(x)\\
u_1(x)
\end{matrix}
\right)
$$
has a unique smooth solution $u:[0,T]\times M\to\mathbb{R}^3$.
\end{theorem}

\section{Discussion on ``Elementary" Methods in Establishing Long-time Results}\label{B}
In this appendix, we shall sketch the method of estimating lifespan for (\ref{EQ0}) in \cite{Notz2010}, and explain why it cannot be improved, thus why the Nash-Moser method in the proof of Theorem \ref{Thm2} is unavoidable.

In \cite{Notz2010}, the lifespan of equation (\ref{EQ0}) with initial data $\varepsilon$-close to a static solution was estimated using a standard continuous induction argument, based on energy estimates for the equation satisfied by geometric quantities (the second fundamental form and components of the velocity etc.). The equation is a complicated quasilinear weakly hyperbolic system. The lifespan so obtained was $\sim\log1/\varepsilon$.

Following \cite{Notz2010}, we consider the evolution of various geometric quantities. We shall write $w=i_0+u$, where $u$ is a perturbation, and suppose $w$ solves the  Cauchy problem of (\ref{EQ0}) on a time interval $[0,T]$:
\begin{equation}\label{EQ0'}
\frac{\partial^2u}{\partial t^2}
=\frac{d\mu(w)}{d\mu_0}\left(-H(w)+\frac{\kappa}{\mathrm{Vol}(w)}\right)N(w),\,
\left(
\begin{matrix}
w(0,x)\\
\partial_tu(0,x)
\end{matrix}
\right)=
\left(
\begin{matrix}
i_0(x)+u_0(x)\\
u_1(x)
\end{matrix}
\right).\tag{EQ0'}
\end{equation}
Let
$$
\sigma=\perp_w\partial_tw,\,
S=\top_w\partial_tu,\,
B_{ij}=(\bar{\nabla}_{\partial_iw}\partial_tu,\partial_jw),
$$
where $\bar{\nabla}$ denotes the connection in $\mathbb{R}^3$. Then $\sigma$ is a scalar function, $S$ is a vector field, and $B=B_{ij}$ is a symmetric second order tensor. Given any local coordinate on $S^2$, the difference between the Christoffel symbols of $g(w)$ and $g_0$ is denoted by $\Gamma_{ij}^k(w)-\Gamma_{ij}^k(g_0)$, which is a tensor.

Now define the following tensors:
$$
\zeta=\left(\begin{matrix}
\sigma\\
h_{ij}(w)-h_{ij}(g_0)
\end{matrix}\right),\,
\xi=\left(\begin{matrix}
u\\
du\\
S\\
B_i^j(w)\\
\Gamma_{ij}^k(w)-\Gamma_{ij}^k(g_0)
\end{matrix}\right).
$$
Here $B_i^j(w)=g^{ik}(w)B_k^j(w)$. We shall directly quote from \cite{Notz2010} the evolution equations of $\zeta$ and $\xi$. For convenience we omit the dependence on $w$, and add an upper circle for geometric quantities induced by $g_0$.

The evolution of $\zeta$ is given by
$$
\begin{aligned}
\partial_t^2\sigma
&=\frac{d\mu}{d\mu_0}\left[\Delta\sigma
+|h|^2\sigma-\frac{\kappa}{\mathrm{Vol}(w)^2}\int_{S^2}\sigma d\mu
+\left(-H+\frac{\kappa}{\mathrm{Vol}(w)}\right)
(\text{div}S+H)\sigma
\right]\\
&\quad+\sigma\left(|\nabla\sigma|^2+S^kS^lh^j_kh_{jl}-2h_k^i\partial_i\sigma S^k\right)
-2\partial_t S^k\partial_k\sigma+2h_{ik}S^i\partial_tS^k,
\end{aligned}
$$
$$
\begin{aligned}
\partial_t^2h_{ij}
&=\frac{d\mu}{d\mu_0}\left\{\Delta h_{ij}
+|h|^2h_{ij}-Hh_{il}h^{l}_j
+\left(-H+\frac{\kappa}{\mathrm{Vol}(w)}\right)
\left[\nabla_i(\Gamma_{jl}^l-\mathring\Gamma_{jl}^l
)+(\Gamma_{ik}^k-\mathring\Gamma_{ik}^k
)(\Gamma_{jl}^l-\mathring\Gamma_{jl}^l)\right]
\right\}\\
&\quad+\frac{d\mu}{d\mu_0}\nabla_jH(\Gamma_{ik}^k-\mathring\Gamma_{ik}^k)
+\frac{d\mu}{d\mu_0}\nabla_iH(\Gamma_{jl}^l-\mathring\Gamma_{jl}^l)
+\frac{d\mu}{d\mu_0}\left(-H+\frac{\kappa}{\mathrm{Vol}(w)}\right)h_{ik}h^{k}_j\\
&\quad+h_{ij}\left(|\nabla\sigma|^2+h^j_kh_{jl}S^kS^l-2h^i_k\partial_i\sigma S^k\right)
+2\partial_t\Gamma^k_{ij}(\partial_k\sigma-h_{kl}S^l).
\end{aligned}
$$
Here the operator $\Delta$ acting on $h_{ij}$ is the trace Laplacian $g_{ij}\nabla^i\nabla^j$.

The evolution of $\xi$ is given by
$$
\partial_t^2u
=\frac{d\mu}{d\mu_0}\left(-H+\frac{\kappa}{\mathrm{Vol}(w)}\right)N(w),
$$
$$
\partial_t^2\partial_ku
=\partial_k\left[\frac{d\mu}{d\mu_0}\left(-H+\frac{\kappa}{\mathrm{Vol}(w)}\right)\right]N(w)
+\frac{d\mu}{d\mu_0}\left(-H+\frac{\kappa}{\mathrm{Vol}(w)}\right)h_k^l\partial_lu,
$$
$$
\begin{aligned}
\partial_t^2S^m
&=\frac{d\mu}{d\mu_0}\left(-H+\frac{\kappa}{\mathrm{Vol}(w)}\right)\left[-\nabla^m\sigma
+\sigma\left(\Gamma_{ij}^k-\mathring\Gamma_{ij}^k\right)g^{km}\right]\\
&\quad+2\partial_t\sigma\left(\nabla^m\sigma-h^m_kS^k\right)-\sigma\frac{d\mu}{d\mu_0}\nabla^mH-2\sigma B_i^m(\nabla^i\sigma-h^i_jS^j)
-2\partial_tS^kB_k^m,
\end{aligned}
$$
$$
\begin{aligned}
\partial_t^2B_i^k
&=-\frac{d\mu}{d\mu_0}(\nabla^k\sigma-h^k_lS^l)\partial_iH
+\frac{d\mu}{d\mu_0}\left(-H+\frac{\kappa}{\mathrm{Vol}(w)}\right)(\Gamma_{il}^l-\mathring\Gamma_{il}^l)(\nabla^k\sigma-h^k_mS^m)\\
&\quad+\frac{d\mu}{d\mu_0}\left[
-\frac{d\mu}{d\mu_0}\partial_tH-\frac{\kappa}{\mathrm{Vol}(w)^2}\int_{S^2}\sigma d\mu
+\left(-H+\frac{\kappa}{\mathrm{Vol}(w)}\right)g^{ml}B_{ml}\right]h_i^k\\
&\quad+\frac{d\mu}{d\mu_0}\left(-H+\frac{\kappa}{\mathrm{Vol}(w)}\right)(\partial_th^k_i+h_i^lB_l^k-B_i^lh_l^k)
-2\partial_tB_i^lB_l^k+2\left[-B_i^l(\partial_l\sigma-h_{lm}S^m)\right](\nabla^k\sigma-h_l^kS^l)\\
&\quad-(\partial_i\sigma-h_{il}S^l)\left[\frac{d\mu}{d\mu_0}\nabla^kH+2(\nabla^m\sigma-h_j^mS^j)B_m^k\right]\\
&\quad+\frac{d\mu}{d\mu_0}\left(-H+\frac{\kappa}{\mathrm{Vol}(w)}\right)g^{km}(\Gamma_{ml}^l-\mathring\Gamma_{ml}^l)(\partial_i\sigma-h_{ij}S^j),
\end{aligned}
$$
$$
\begin{aligned}
\partial_t^2\Gamma_{ij}^k
&=\frac{d\mu}{d\mu_0}(-\nabla_iHh_j^k-\nabla_jHh_i^k+\nabla^kHh_{ij})
+\frac{d\mu}{d\mu_0}\left(-H+\frac{\kappa}{\mathrm{Vol}(w)}\right)\nabla_ih_j^k\\
&\quad+\frac{d\mu}{d\mu_0}\left(-H+\frac{\kappa}{\mathrm{Vol}(w)}\right)\left[
h_j^k(\Gamma_{il}^l-\mathring\Gamma_{il}^l)+h_i^k(\Gamma_{jl}^l-\mathring\Gamma_{jl}^l)-h_{ij}g^{km}(\Gamma_{ml}^l-\mathring\Gamma_{ml}^l)\right]\\
&\quad-2\partial_t\Gamma_{ij}^lB_l^k-2\partial_th_{ij}(\nabla^k\sigma-h_m^kS^m)+2h_{ij}(\nabla^l\sigma-h_m^lS^m)B_l^k.
\end{aligned}
$$

Observe that all geometric quantities are smooth in $(\zeta,\xi)$; for example, the induced metric $g(w)$ and the Radon-Nikodym derivative $d\mu/d\mu_0$ are smooth functions of $du$, and the (scalar) mean curvature $H$ is the contraction of $g(w)$ with $h(w)$, and
$$
\partial_k\frac{d\mu}{d\mu_0}
=(\Gamma_{ik}^i-\mathring\Gamma_{ik}^i)\frac{d\mu}{d\mu_0}.
$$
We thus reduce these evolution equations to much terser forms:
\begin{proposition}\label{P5.1}
There is a $\delta_6>0$ such that if $|u|_{C_t^2C_x^2}<\delta_6$, then the system of $(\zeta,\xi)$ can be re-arranged to a terser form as
\begin{equation}\label{Ev.Geom}
\begin{aligned}
\partial_t^2\zeta
&=A\zeta
+I_0^1(\zeta,\xi,\mathring{\nabla}\zeta,\partial_t\xi,\mathring{\nabla}\xi)\cdot\zeta
+I_1^1(\zeta,\xi,\mathring{\nabla}\zeta,\partial_t\xi,\mathring{\nabla}\xi)\cdot\mathring{\nabla}\zeta
+I_2^1(\zeta,\xi,\mathring{\nabla}\xi)\cdot\mathring{\nabla}^2\zeta,\\
\partial_t^2\xi
&=J^2(\zeta,\partial_t\zeta,\mathring{\nabla}\zeta)
+Q^2(\zeta,\xi,\partial_t\zeta,\mathring{\nabla}\zeta,\partial_t\xi),
\end{aligned}
\end{equation}
where $A$ is the elliptic operator given by
$$
A\left(
\begin{matrix}
\sigma\\
h_{ij}
\end{matrix}
\right)
=\left(
\begin{matrix}
\displaystyle{\mathring\Delta\sigma+2\sigma-\frac{6}{4\pi}\int_{S^2}\sigma d\mu_0}\\
\mathring\Delta h_{ij}
\end{matrix}
\right),
$$
and the $I,J$'s are tensors that depend smoothly on their arguments, and vanish linearly when the arguments tend to zero; the $Q^2$ tensor depends smoothly on its arguments, and vanishes quadratically when the arguments tend to zero.
\end{proposition}
\begin{proof}
This is a direct calculation using Taylor's formula. For the $\zeta$ component, the only difficulty is the equation satisfied by $h_{ij}$. However, note that e.g. for any symmetric section $W$ of $T^*(S^2)\otimes T^*(S^2)$, we have
$$
|\mathring{h}|^2W_{ij}-\mathring{H}\mathring{h}^l_jW_{il}\equiv0,
$$
so
$$
\frac{d\mu}{d\mu_0}(\Delta h_{ij}
+|h|^2h_{ij}-Hh_{il}h^{l}_j)
$$
is reduced to the form
$$
\mathring{\Delta}(h_{ij}-\mathring{h}_{ij})
+(\text{function of }du,\,\mathring{\nabla}du)\cdot\mathring{\nabla}^2h_{ij}
+(\text{function of }h)\cdot(h-\mathring{h})\cdot h_{ij}.
$$
Similarly, the only difficulty for the $\xi$ component are terms involving $dH$. However, we compute
$$
\begin{aligned}
dH
&=d\left[(g^{ij}-\mathring g^{ij})h_{ij}+\mathring g^{ij}(h_{ij}-\mathring h_{ij})-2\right]\\
&=-h_{ij}\left(\Gamma_{lk}^ig^{lj}-\mathring\Gamma_{lk}^i\mathring g^{lj}
+\Gamma_{lk}^jg^{li}-\mathring\Gamma_{lk}^j\mathring g^{li}\right)dx^k+(g^{ij}-\mathring g^{ij})\partial_kh_{ij}dx^k\\
&\quad+\partial_k\mathring g^{ij}(h_{ij}-\mathring h_{ij})dx^k
+\mathring g^{ij}\partial_k(h_{ij}-\mathring h_{ij})dx^k,
\end{aligned}
$$
where we used the compatibility condition
$$
dg_{ij}=g_{il}\Gamma_{jk}^ldx^k+g_{jl}\Gamma_{ik}^ldx^k.
$$
So $dH$ can be expressed linearly in first order derivatives of $\zeta$ and no derivative of $\xi$.
\end{proof}

Once the original equation (\ref{EQ0'}) is reduced to the weakly hyperbolic quasilinear system (\ref{Ev.Geom}), a standard fixed-point type argument with the aid of Proposition \ref{SolvingWHLS} will give the local well-posedness result of that system. Nevertheless, we should note that (\ref{Ev.Geom}) is derived from (\ref{EQ0'}). Going back from the solution of (\ref{Ev.Geom}) to the original unknown $u$ in (\ref{EQ0'}) shall encounter obstructions resulting from geometric constraints, namely the Gauss-Codazzi equations for evolving submanifolds. We do not yet know if the verification of these geometric-dynamical constraints is possible; should it be possible, it is certainly as lengthy as the Nash-Moser iteration scheme. From the analysis above, it is better illustrated why the Nash-Moser scheme is unavoidable in solving (\ref{EQ0}).

However, if (\ref{EQ0'}) is already proved to be locally well-posed, system (\ref{Ev.Geom}) shall provide an estimate for the lifespan if some sufficiently high Sobolev norm of $u[0]$ is $\varepsilon$-small. This is the method employed in \cite{Notz2010} to estimate the lifespan, and it is applicable to generic closed constant mean curvature hypersurfaces in a generic ambient manifold. Nevertheless, the result produced by this approach for perturbation of $S^2$ in $\mathbb{R}^3$ is not optimal compared to our proof of Theorem \ref{Thm2}. Let us briefly explain the reason below. The idea is to estimate the energy norm
$$
E_n(t):=\|\partial_t^2\zeta(t)\|_{H^n_x}
+\|\partial_t\zeta(t)\|_{H^{n+1}_x}
+\|\zeta(t)\|_{H^{n+1}_x}+\|\partial_t\xi(t)\|_{H^{n+1}_x}+\|\xi(t)\|_{H^{n+1}_x}
$$
for some large $n$ by proving inequalities of the form
$$
\frac{d}{dt}E_n(t)\leq C(\varepsilon)E_n(t),
$$
and use a continuous induction argument to make the quantity $\varepsilon E_n(t)$ bounded. A problem then arises: in the evolution equations of $u$ and $du$, we find
$$
-H+\frac{\kappa}{\text{Vol}(w)}
=(g^{ij}-\mathring{g}^{ij})h_{ij}+\mathring{g}^{ij}(h_{ij}-\mathring{h}_{ij})+\frac{\kappa}{\text{Vol}(w)}-2.
$$
Thus, as $\zeta,\xi\to0$, the difference $-H+\kappa/\text{Vol}(w)$ only vanishes linearly. Consequently, in establishing energy estimate for (\ref{Ev.Geom}), we find that the information of (perturbed) non-growing modes is lost compared to the linearized problem (\ref{LinearizationInComponents}), and the best to expect for energy estimate is
$$
E_n(t)\lesssim e^{Ct},
$$
where no smallness for $C$ can be guaranteed even if $u[0]$ is small. Thus the standard continuous induction argument, employed in e.g. \cite{Klainerman1985} or section 6.4 of \cite{Hormander1997}, will only give the lifespan estimate
$$
T\simeq\log\frac{1}{\varepsilon},
$$
as obtained in \cite{Notz2010}. This loss of information on slow growth of modes near zero is because that this method does not employ the stability of $S^2$, and thus the estimate $T\simeq\log1/\varepsilon$ for lifespan is not optimal in this case.

\end{spacing}
\end{document}